\documentclass[12pt]{article} 
\usepackage{amsmath,amssymb,latexsym,amsthm,tikz, enumerate,amscd,hyperref} \usepackage[abbrev]{amsrefs} 
\usepackage{graphicx}
\usepackage{todonotes}
\newtheoremstyle{dotless}{}{}{\itshape}{}{\bfseries}{}{ }{}
\usetikzlibrary{decorations.markings}

\newtheorem{Theorem}{Theorem}[section]

\newtheorem{Prop}[Theorem]{Proposition} 
\newtheorem{corollary}[Theorem]{Corollary} 
\newtheorem{Question}[Theorem]{Question} 
\newtheorem{lemma}[Theorem]{Lemma}
\newtheorem{conjecture}[Theorem]{Conjecture}
\newtheorem*{theorem*}{Theorem}
\newtheorem*{Kquestion}{The $K(\pi,1)$-Question}

\newtheorem*{thmA}{Theorem A}
\newtheorem*{thmB}{Theorem B}
\newtheorem*{thmC}{Theorem C}

\theoremstyle{definition} 
\newtheorem{definition}[Theorem]{Definition} 
\newtheorem{remark}[Theorem]{Remark}
\newtheorem{remarks}[Theorem]{Remarks}
\newtheorem{example}[Theorem]{Example}

\newtheorem*{acknowledgements}{Acknowledgements}

\newtheorem*{Remark}{Remark}

\newtheorem*{Conjecture2}{The Action Dimension Conjecture}

\DeclareMathOperator{\aster}{\text{\LARGE{\textasteriskcentered}}}

\DeclareMathOperator{\piprod}{\raisebox{-0.1em}{\huge{$\pi$}}\kern -0.2em}

\newcommand{\ca}{\mathcal {A}} 
 
\newcommand{\cac}{\mathcal {C}}

\newcommand{\cg}{\mathcal {G}} 
\newcommand{\ch}{\mathcal {H}}

\newcommand{\cp}{\mathcal {P}} 
\newcommand{\cq}{\mathcal {Q}}
 
\newcommand{\cs}{\mathcal {S}} 
 
\newcommand{\cu}{\mathcal {U}}

\newcommand{\gq}{\lim{G\cq}}
\newcommand{\gsl}{\lim{G\cs(L)}}

 \newcommand{\cc}{\mathbb{C}}

\newcommand{\nn}{{\mathbb N}} 
 
\newcommand{\QQ}{{\mathbb{Q}}} 
\newcommand{\rr}{{\mathbb R}} 
\newcommand{\hh}{{\mathbb H}}

\newcommand{\zz}{{\mathbb Z}}

\newcommand{\ux}{{\underline{X}}}
\newcommand{\uy}{{\underline{Y}}}

\newcommand{\minus}{{-1}}

\newcommand{\ti}{\tilde}

\newcommand{\h}{\hat}

\newcommand{\ol}{\overline}

\newcommand{\fiq}{FI\cq(\ca)} 
\newcommand{\viq}{VI\cq(\ca)} 
\newcommand{\iq}{I\cq(\ca)}

\newcommand{\cd}{\operatorname{cd}} 
 
\newcommand{\gdim}{\operatorname{gdim}}

\newcommand{\obdim}{\operatorname{obdim}}
\newcommand{\pobdim}{\operatorname{pobdim}}
\newcommand{\actdim}{\operatorname{actdim}}
\newcommand{\embdim}{\operatorname{embdim}}

\newcommand{\vk}{\operatorname{vk}}
\newcommand{\vktwo}{\vk}
\newcommand{\cvk}{\operatorname{\nu\gk}}

\newcommand{\rk}{\mathrm{rk}}

\newcommand{\raag}{\mathrm {RAAG}} 
\newcommand{\racg}{\mathrm {RACG}}

\def\clap#1{\hbox to 0pt{\hss#1\hss}}

\newcommand{\comment}[1]{} 
 
\newcommand{\ga}{\alpha} 
\newcommand{\gb}{\beta} 
\newcommand{\gd}{\delta} 

\newcommand{\geps}{\varepsilon} 
\newcommand{\gf}{\varphi} 
\newcommand{\gi}{\iota} 
\newcommand{\gr}{\rho} 
\newcommand{\gs}{\sigma} 
 
\newcommand{\gt}{\tau} 
 
\newcommand{\gk}{\kappa}

\newcommand{\T}{T}

\newcommand{\gG}{\Gamma} 
\newcommand{\gD}{\Delta} 
\newcommand{\gO}{\Omega}

\newcommand{\gS}{\Sigma}

\newcommand{\vertex}{\operatorname{Vert}}

\newcommand{\Ima}{\operatorname{Im}} 
\newcommand{\Lk}{\operatorname{Lk}} 
 
\newcommand{\St}{\operatorname{St}}

\newcommand{\cat}{\operatorname{CAT}} 
\newcommand{\cone}{\operatorname{Cone}} 
\newcommand{\codim}{\operatorname{codim}}

\newcommand\mapright[1]{\smash{\mathop{\longrightarrow}\limits^{#1}}}

\newcommand{\itemEq}[1]{%
         \begingroup%
         \setlength{\abovedisplayskip}{0pt}%
         \setlength{\belowdisplayskip}{0pt}%
         \parbox[c]{\linewidth}{\begin{flalign}#1&&\end{flalign}}%
         \endgroup}

\newenvironment{enumerate1}{ 
\begin{enumerate}[\upshape (1)]}
	{ 
\end{enumerate}
} 

\newenvironment{enumeratei}{ 
\begin{enumerate}[\upshape (i)]}
	{ 
\end{enumerate}
} 

\newenvironment{enumeratei'}{ 
\begin{enumerate}[\upshape (i)$'$]}
	{ 
\end{enumerate}
} 

\newenvironment{enumerate1'}{ 
\begin{enumerate}[\upshape (1)$'$]}
	{ 
\end{enumerate}
} 

\newenvironment{enumerateI}{ 
\begin{enumerate}[\upshape (I)]}{ 
\end{enumerate}
} 
\newenvironment{enumeratea}{ 
\begin{enumerate}[\upshape (a)]}{ 
\end{enumerate}
}

\newenvironment{enumeratea'}{ 
\begin{enumerate}[\upshape (a)$'$]}{ 
\end{enumerate}
}

\usepackage{color}
\usepackage[outerbars,color]{changebar}
\usepackage{ifthen}

\newboolean{usecolor}
\setboolean{usecolor}{true}

\ifthenelse{\boolean{usecolor}}
{
 \definecolor{colore}{cmyk}{0,1,0.6,0}
 \definecolor{coloregen}{cmyk}{0.7,0,1,0}
 \definecolor{coloresimo}{cmyk}{1,0.6,0,0}
}
{
 \definecolor{colore}{cmyk}{0,0,0,1}
 \definecolor{coloregen}{cmyk}{0,0,0,1}
 \definecolor{coloresimo}{cmyk}{0,0,0,1}
}

\newcommand{\nota}[1]
{\-\marginpar[\raggedleft\tiny \textcolor{colore}{#1}]{\raggedright\tiny \textcolor{colore}{#1}}}

\numberwithin{equation}{section} 
\begin{document}
\title{Action dimensions of some simple complexes of groups}

\author{Michael W. Davis 
\and {Giang Le}  
\and{Kevin Schreve}
}
\date{\today} \maketitle

\begin{abstract} 

\noindent
The action dimension of a discrete group $G$ is the minimum dimension of contractible manifold that admits a proper $G$-action.
We compute the action dimension of the direct limit of a simple complex of groups for several classes of examples including: 1) Artin groups, 2) graph products of groups, and 3) fundamental groups of aspherical complements of arrangements of affine hyperplanes.
\medskip 

\noindent
\textbf{AMS classification numbers}. Primary: 20F36, 20F55, 20F65, 57S30, 57Q35, 
Secondary: 20J06, 32S22
\smallskip

\end{abstract}

\section*{Introduction}\label{s:intro}
Suppose $G$ is a discrete, torsion-free group with classifying space $BG$.
Its \emph{geometric dimension}, $\gdim(G)$, is the smallest dimension of a model for its classifying space $BG$ by a CW complex. 
This number is equal to $\cd G$, the cohomological dimension of $G$ (provided $\cd G\neq 2$). 
Its \emph{action dimension}, $\actdim G$, is the smallest dimension of a model for $BG$ by a manifold.  
In other words, $\actdim G$ is the minimum dimension of a thickening of a CW model for $BG$ to a manifold, possibly with boundary. 
It follows that $\gdim G\le \actdim G$ with equality if and only if $BG$ is homotopy equivalent to a closed manifold. 
From general principles, $\actdim G\le 2\gdim G$.

A common method of constructing  groups and their classifying spaces is to use the notion of a ``complex of groups'' (cf.\ \cite{bh}).  Here we will only use the easier notion of a \emph{simple complex of groups} over a poset $\cq$. By definition, this means a functor, $G\cq$, from $\cq$ to the category of groups and monomorphisms. So, $G\cq$ is the following data: a collection of groups $\{G_\gs\}_{\gs\in \cq}$ and monomorphisms $\phi_{\gs\gt}:G_\gt\to G_\gs$, defined when $\gt<\gs$ and satisfying $\phi_{\gs\gt}\phi_{\gt\mu}=\phi_{\gs\mu}$ when $\mu<\gt<\gs$. Suppose we have models for the $BG_\gs$ and realizations for the homomorphisms $\phi_{\gs\gt}$ by maps
$\ol{\phi}_{\gs\gt}:BG_\gt\to BG_\gs$. One then can glue together the $\{BG_\gs\}_{\gs\in \cq}$ (or more precisely iterated mapping cylinders of the $\ol{\phi}_{\gs\gt}$) to form a space $BG\cq$ called the \emph{aspherical realization} of $G\cq$. If the geometric realization of $\cq$ is simply connected, then it follows from van Kampen's Theorem that $\pi_1(BG\cq)$ is the direct limit, $G$, of the system of groups $\{G_\gs, \phi_{\gs\gt}\}$. The space $BG\cq$ may or may not be aspherical. This is the ``$K(\pi,1)$-Question'' for $G\cq$. When the answer is affirmative, $BG\cq$ is a model for $BG$. In this case, we get an upper bound for $\gdim G$ in terms of the geometric dimensions of the $G_\gs$. Similarly, if each $BG_\gs$ is modeled by a manifold with boundary, $M_\gs$, and each $\ol{\phi}_{\gs\gt}$ is homotopic to an embedding $f_{\gs\gt}:M_\gt\to \partial M_\gs$, then we can glue together suitably thickened versions of the $M_\gs$ to get a thickening of $BG\cq$ to a manifold with boundary $M$. So, provided the $K(\pi,1)$-Question for $G\cq$ has a positive answer,  one gets an upper bound for the action dimension of $G$ in terms of the action dimensions of the $\{G_\gs\}_{\gs \in \cq}$  

In \cite{bkk}, Bestvina, Kapovich, and Kleiner define a number, $\obdim G$, called the ``obstructor dimension" of $G$. It is a lower bound for the action dimension of $G$. It is based on the classical van Kampen obstruction, $\vk^n(K)$, for embedding a finite simplicial complex $K$ into $\rr^n$. 
This obstruction is a cohomology class with $\zz_2$ coefficients in the configuration space of unordered pairs of distinct points in $K$.
Suppose $EG$ denotes the universal cover of $BG$. The idea of \cite{bkk} is to find a complex $K$ and a coarse embedding of $\cone_\infty K$ into $EG$, where $\cone_\infty K$ means the cone of infinite radius on $K$. It is proved in \cite{bkk} that $\vk^n(K)$ is also an obstruction to a coarse embedding of $\cone_\infty K$ into any contractible $(n+1)$-manifold; hence, when $\vk^n(K)\neq 0$, $\actdim G\ge n+2$. In \cite{bkk} the \emph{obstructor dimension} of $G$  is defined to be $n+2$, where $n$ is the largest integer so that there exists a complex $K$ with $\vk^n(K)\neq 0$, together with a coarse embedding of $\cone_\infty K$ into $EG$. We also shall have occasion to use a variation of this notion due to Yoon \cite{yoon}, called the ``proper obstructor dimension of $G$.'' His idea is to consider coarse embeddings $T\to EG$ where $T$ is a contractible simplicial complex, not necessarily of the form $T=\cone_\infty K$. 

In this paper we use the following two techniques to compute the action dimension for certain groups which are direct limits of  simple complexes of groups.

\begin{enumerateI}
\item
(\emph{Gluing}). Construct a thickening of $BG\cq$ by gluing together manifolds with boundary that are models for the $BG_\gs$, hence, establishing an upper bound for $\actdim G$.  The pieces that are to be glued together will have the form $M_\gs\times D_\gs$ where $M_\gs\sim BG_\gs$ and where $D_\gs$ is a ``dual disk.'' Two such pieces will be glued together along a piece that is a common codimension-$0$ submanifold of both boundaries. (The details of this method are described in Section~\ref{s:gluing}.)
\item
(\emph{Obstructors}). Lower bounds for $\obdim G$ (and hence, for $\actdim G$) are established by finding obstructors for $G$. In most of our examples the coarse obstructor will be a finite union of contractible manifolds, containing a common basepoint, each of which is the universal cover of some closed aspherical manifold. These contractible manifolds are called \emph{sheets}. When a sheet can be compactified to a disk, it is homeomorphic to the cone on a sphere and the coarse obstructor has the form $\cone_\infty K$ where $K$ is some configuration of spheres ($K$ will often be a ``polyhedral join'' of spheres, cf.\ Definition~\ref{d:polyjoin}). So,  such  cases reduce to calculating van Kampen obstructions $\vk^n(K)$.
\end{enumerateI}

A necessary condition for $BG\cq$ to be aspherical is that the geometric realization $\vert \cq\vert$ of the poset $\cq$ is contractible (see Remark~\ref{r:retract}). Often this will be automatic, since if $\cq$ has a minimum element for which the corresponding local group is the trivial group, then $\vert\cq\vert$ is a cone. If $\cq$ has such a minimum element, then at the final stage of (I) we will need to glue a disk onto the result of previous gluings.  This will entail that each of previous manifolds $M_\gs$ has nonempty boundary. However, if $\cq$ has no such minimum, then for any minimal element $\gs$ of $\cq$ one can allow $M_\gs$ to be a closed manifold. We shall return to this point later in the Introduction. 

In most applications $\cq$ will be the poset $\cs(L)$ of simplices in some simplicial complex $L$, including the empty simplex (so  $\vert\cs(L)\vert$ will be a cone).

A prototypical example is the case where $G=A_L$, the right-angled Artin group (or ``$\raag$'') associated to a $d$-dimensional flag complex $L$. Then the simple complex of groups is the \emph{Artin complex}, $A\cs(L)= \{A_\gs\}_{\gs\in\cs(L)}$ of spherical Artin subgroups, $BA_L$ is the standard model for its classifying space as a union of tori, and $\gdim A_L =d+1$. It turns out that the appropriate obstructor is the polyhedral join of $0$-spheres $O_1L$, called the \emph{octahedralization} of $L$. The following two results are proved in \cite{ados}*{Theorems 5.1 and 5.2}.
\begin{enumeratei}
\item
If $H_d(L,\zz_2)\neq 0$, then $\vk^{2d}(O_1L)\neq 0$. Hence, $\actdim A_L=\obdim A_L= 2d+2$.
\item
If $H_d(L,\zz_2)=0$ and $d\neq 2$, then $\actdim A_L\le 2d+1$.
\end{enumeratei}

We discuss below  three generalizations of $\raag$s: 1) general Artin groups, 2) graph products of fundamental groups of closed aspherical manifolds, and 3) fundamental groups of aspherical complements of affine hyperplane arrangements. In each case we prove results similar to (i) and (ii) above. In all three cases the relevant obstructor will be a polyhedral join $O_mL$ of $(m-1)$-spheres. 
Although we conjecture that a result similar to (ii) holds in both cases 1) and 2), the proof of (ii) for $\raag$s uses a special argument. ($A_L$ is a subgroup of the right-angled Coxeter group corresponding to $O_1L$; if $\vk^{2d}(O_1L)=0$, then $O_1L$ embeds in a flag triangulation of $S^{2d}$ and hence, the Coxeter group is a subgroup of a Coxeter group that acts cocompactly on  a contractible $(2d+1)$-manifold.) When the group is not a $\raag$, we instead use gluing methods to prove the analog of (ii) in the case where $L$ embeds in a contractible simplicial complex of the same dimension $d$. We abbreviate this condition by saying that $L$ is EDCE. (Note: if $L$ is  EDCE, then $H_d(L;\zz_2)=0$.)
\medskip

\noindent
1) (\emph{General Artin groups}). Suppose $A_L$ is the Artin group, where the simplicial complex $L$ is the nerve of the associated Coxeter system. The \emph{Artin complex} $A\cs(L)$ is the simple complex of groups $\{A_\gs\}_{\gs \in \cs(L)}$ of spherical Artin subgroups. The group $A_L$ is the direct limit of the system of groups defined by $A\cs(L)$. Conjecturally, the $K(\pi,1)$-Question has an affirmative answer for all Artin complexes. This is known in many cases (see \cite{cd1}). When the answer is affirmative, $BA\cs(L)$ is covered by the ``Salvetti complex'' of $A_L$ and $\gdim A_L =\dim BA\cs(L) = d+1$. The next result is proved in \cite{dh16} and \cite{giang}. We shall give the details of the arguments in Sections \ref{s:gluing} and \ref{s:obstructors}.

\begin{thmA}[cf.\ \cite{dh16} and Proposition~\ref{p:artinEDCE}]\label{A}
Suppose the $K(\pi,1)$-Question has a positive answer for $A\cs(L)$. 
\begin{enumeratei}
\item
If $H_d(L;\zz_2)\neq 0$, then $\actdim A_L=\obdim A_L=2d+2=2\gdim A_L$.
\item
If $L$ is EDCE, then $\actdim A_L\le 2d+1$.
\end{enumeratei}
\end{thmA}

In \cite{dh16} it is proved that the relevant obstructor for the Artin complex is a polyhedral join of $0$-spheres, $O_1L_\oslash$, where $L_\oslash$ is a certain subdivision of $L$ whose simplices index the ``standard free abelian subgroups'' of $A_L$. As before, $\cone(O_1L_\oslash)$ coarsely embeds in $EA_L$. The calculation of the van Kampen obstruction in part (i) of Theorem A is then the same as its calculation for $O_1L$ in the case of a $\raag$.
\medskip

\noindent
2) (\emph{Graph products}). Suppose $\{G_v\}_{v\in V}$ is a collection of groups indexed by the vertex set $V$ of a simplicial graph $L^1$. The \emph{graph product} $G$ is the quotient of the free product of the $G_v$ by the relations that $G_v$ and $G_w$ commute whenever $\{v,w\}$ is an edge of $L^1$. Let $L$ be the flag complex associated to $L^1$. For each simplex $\gs\in \cs(L)$, $G_\gs$ denotes the direct product of the $G_v$ over the vertex set of $\gs$. The \emph{graph product complex} $G\cs(L)$ is the simple complex of groups $\{G_\gs\}_{\gs\in \cs(L)}$ where the monomorphisms $\phi_{\gs\gt}$ are the natural inclusions. Obviously, $G$ is the direct limit of the $G_\gs$. Moreover, the $K(\pi,1)$-Question for $G\cs(L)$ always has a positive answer (see \ref{ss:flag}). One can essentially determine $\actdim G$ in terms of the $\actdim G_v$ and we do so in Sections \ref{s:gluing} and \ref{s:obstructors}. The most interesting case is when each $G_v$ is the fundamental group of a closed aspherical manifold $M_v$. A $\raag$ is the special case where each $M_v=S^1$. To simplify the statements of our results, assume each vertex manifold $M_v$ has the same dimension $m$. Then $\gdim G=m(d+1)$. Using appropriate thickenings of the $M_\gs$ we can use gluing technique (I) to construct a manifold $M$ of dimension $(m+1)(d+1)$ which is a model for $BG$. By gluing together standard lifts of the $M_\gs$, we get a coarsely embedded $\cone_\infty(O_{m-1}L)$ in $EG$. We then get the following analog of Theorem A.

\begin{thmB}[cf.\ Corollaries~\ref{c:closedgprdt2}, \ref{c:gp} and \ref{c:gprdt2EDCE}]\label{B}
Suppose, as above, that $G$ is the graph product of fundamental groups of closed, aspherical $m$-manifolds $M_v$ over a $d$-dimensional flag complex $L$. 
\begin{enumeratei}
\item
If $H_d(L,\zz_2)\neq 0$, then $\actdim G=\obdim G= (m+1)(d+1)$.
\item
If $L$ is EDCE, then $\actdim G\le (m+1)(d+1)-1$.
\end{enumeratei}
\end{thmB}

A mild generalization of this is the case where each local group $G_\gs$ is the fundamental group of a closed aspherical manifold $M_\gs$ of dimension $m(\dim(\gs)+1)$ and where the $\phi_{\gs\gt}$ are realized by embeddings $f_{\gs\gt}: M_\gt\hookrightarrow M_\gs$. In Subsection~\ref{ss:gprdts}, we call such a system $\{M_\gs\}_{\gs\in \cs(L)}$ a \emph{simple complex of closed aspherical manifolds}. When $L$ is a flag complex, the $K(\pi,1)$-Question for $G\cs(L)$ has a positive answer (cf. Theorem~\ref{t:flag}) and the conclusion of Theorem B holds without change.
\medskip

\noindent
3) (\emph{Aspherical complements of hyperplane arrangements}). Suppose $\ca$ is an arrangement of affine hyperplanes in $\cc^{n}$. Let $M(\ca)$ denote the complement $\cc^{n} - \bigcup\ca$. The relevant poset $\cq$ is the intersection poset of $\ca$. Its elements are the proper subspaces $\gs$ of $\cc^{n}$ which are intersections of hyperplanes in $\ca$, ordered by reverse inclusion. The minimal elements of $\cq$ are a family of parallel subspaces, and the arrangement is essential if these are zero dimensional. For each $\gs\in \cq$, there is a central arrangement $\ca_\gs$ in the subspace normal to $\gs$ in $\cc^{n}$. Put $G=\pi_1(M(\ca))$ and $G_\gs=\pi_1(M(\ca_\gs))$. This is the data for a simple complex of groups $G\cq$. Of course, $M(\ca)$ need not be aspherical; however, if it is, then so are the $M(\ca_\gs)$. Let us assume that $M(\ca)$ is aspherical. 

Since $M(\ca)$ is a $2n$-manifold, $\actdim G\le 2n$. If $\ca$ is central (meaning that the hyperplanes are linear), then $M(\ca)$ deformation retracts onto the complement of the hyperplane arrangement in the unit sphere $S^{2n-1}$ in $\cc^{n}$ and hence, $\actdim G\le 2n-1$. If $\ca$ is not essential, then the parallel subspaces can be deformation retracted, and again  $\actdim G\le 2n-1$. On the other hand, for essential aspherical arrangements we have the following theorem.

\begin{thmC}[cf.\ Theorem \ref{t:nothe}]\label{C}
Let $\ca$ be an essential arrangement of affine hyperplanes in $\cc^n$. Suppose $M(\ca)$ is aspherical. Let $G = \pi_1(M(\ca))$.
If $\ca$ does not decompose as a product with a factor equivalent to a nontrivial central arrangement, then $\actdim G=\obdim G=2n$.
\end{thmC} 
This theorem implies that if $\ca$ decomposes as a product of irreducibles and $k$ is the number of factors which are irreducible central arrangements, then $\actdim G=2n-k$.

To understand why this is an analog of the previous theorems, two points require explanation. First, $\vert\cq(\ca)\vert$ is homotopy equivalent to a wedge of $(n-1)$-spheres and, when $\ca$ is not central, then there is at least one sphere in the wedge. In particular, $H_{n-1}(\vert\cq\vert,\zz_2)\neq 0$. Secondly, there is a simplicial complex $\iq$, called the ``irreducible complex,'' such that $\iq$ is homotopy equivalent to $\vert\cq(\ca)\vert$, and so that the simplices of $\iq$ index the standard free abelian subgroups in $G$. So, $\iq$ is the analog to $L_\oslash$. As before, the relevant obstructor is the polyhedral join $O_1(\iq)$. 

The computations are evidence for the following conjecture connecting action dimension to $L^2$-cohomology. 

\begin{Conjecture2}\textup{(Davis-Okun \cite{do01})}.
If the $i^{th}$ $L^2$-Betti number of $G$ is nonzero, then $\actdim G \ge 2i$.
\end{Conjecture2}

For example, if $\ca$ is an irreducible, essential, affine hyperplane arrangement in $\cc^n$ such that $M(\ca)$ is aspherical, then the $L^2$-Betti number of $\pi_1(M(\ca))$ are zero in degrees $\neq n$ and if the arrangement is irreducible and not central, then the $n^{th}$ $L^2$-Betti number of $\pi_1(M(\ca))$ is nonzero,  cf.\ \cite{djl}. 

In a forthcoming paper, we will determine the action dimensions for some more examples of simple complexes of groups. These computations provide further evidence for the Action Dimension Conjecture.

This paper is organized as follows. In Section \ref{s:sc}, we review simple complexes of groups and explain our main examples. In Section \ref{s:gluing}, we discuss thickenings of simplicial complexes and explain our method of gluing together manifolds. In Section \ref{s:examples}, we perform this operation to give upper bounds on action dimension for Artin groups, graph products and hyperplane complements. In Section \ref{s:obstructors}, we review the van Kampen obstruction, introduce our main example of an obstructor complex, and use it to give lower bounds for the action dimension of graph products and other simple complexes of groups. In Section \ref{hypcomp}, we again use obstructor complexes to compute the action dimension of fundamental groups of hyperplane complements.

\begin{acknowledgements}
This material is based upon work supported by the National Science Foundation under Grant No. DMS-1440140 while the authors were in residence at the Mathematical Sciences Research Institute in Berkeley, California, during the fall semester of 2016. The first author was partially supported by an NSA grant and he received support in the spring of 2017 as a Simons Fellow at the Isaac Newton Institute.  The third author is partially supported by NSF grant DMS-1045119. This material is based upon work supported by the National Science Foundation under Award No. 1704364. 
\end{acknowledgements}

\section{Simple complexes of groups}\label{s:sc}

\subsection{The basic construction}\label{ss:posetgps}
We begin by reviewing the theory of simple complexes of groups as developed in Bridson-Haefliger \cite[II.12]{bh}.

Let $\cq$ be a poset. As in the Introduction, a \emph{simple complex of groups $G\cq$ over $\cq$} is a collection of groups $\{G_\gs\}_{\gs \in \cq}$ and monomorphisms $\phi_{\gs\gt}:G_\gt\to G_\gs$ defined whenever $\gt<\gs$. The $G_\gs$ are the \emph{local groups}. Furthermore, $G\cq$ must be a functor from $\cq$ to the category of groups and monomorphisms in the sense that 
$\phi_{\gs\gt}\phi_{\gt\mu}= \phi_{\gs\mu}$ whenever $\mu<\gt<\gs$. 
(In \cite[II.12.11, p.\,375]{bh} the order relation on $\cq$ is reversed.)
A \emph{simple morphism} $\psi =(\psi_\gs)$ from $G\cq$ to a group $G$ is a function which assigns to each $\gs\in \cq$ a homomorphism $\psi_\gs:G_\gs\to G$ such that $\psi_\gt= \psi_\gs \phi_{\gs\gt}$ whenever $\gt<\gs$. The simple morphism $\psi$ is \emph{injective on local groups} if each $\psi_\gs$ is injective.
Such a simple complex of groups $G\cq=\{G_\gs, \phi_{\gs\gt}\}$;  has a direct limit, denoted $\gq$.
	
For each $\gs\in \cq$, there is a canonical homomorphism $\gi_\gs: G_\gs\to \gq$, hence, a \emph{canonical simple morphism} $\gi:G\cq \to \gq$. The simple complex of groups $G\cq$ is \emph{developable} if $\gi$ is injective on local groups.
The direct limit has the universal property that for any group $H$ and simple morphism $\psi: G\cq\to H$, there is a unique homomorphism $\hat{\psi}: \gq\to H$ such that $\psi_\gs=\hat{\psi} \gi_\gs$ (see \cite[II.12.13, p.\,376]{bh}). 

The \emph{order complex} of a poset $\cp$ is the simplicial complex whose simplices are the totally ordered finite subsets $\{\gt_0,\dots, \gt_k\}$ of elements of $\cp$ (where $\gt_0<\cdots < \gt_k$). 
The underlying topological space of  order complex of $\cp$ is denoted $\vert\cp\vert$ and is called the \emph{geometric realization of $\cp$}. If $\cp^{opp}$ denotes the opposite poset where the order relations are reversed, then $\vert \cp^{opp}\vert$ is isomorphic to $\vert\cp\vert$. There are two natural stratifications of $\vert\cp\vert$ both indexed by $\cp$:
	\begin{equation}\label{e:posetstrata}
	\vert\cp\vert_\gs := \vert\cp_{\le \gs}\vert\quad \text{and}\quad \vert\cp\vert^\gs := \vert\cp_{\ge \gs}\vert .
	\end{equation}
In the first stratification, inclusion of one stratum into another corresponds to the original order relation on $\cp$; in the second one, the order relation is reversed. So, we should regard $\{\vert\cp\vert^\gs \}$ as being indexed by $\cp^{opp}$. Given $x\in \vert\cp\vert$, let $\gs(x)$ be the index of the smallest stratum $\vert\cp\vert^\gs$ containing $x$.

Given a simple complex of groups $G\cq$ and a simple morphism $\psi:G\cq\to G$ that is injective on local groups, one can define a poset $D(\cq,\psi)$ equipped with a $G$-action, as well as, a space $D(\vert\cq\vert,\psi)$ also equipped with a $G$-action. The poset $D(\cq,\psi)$ is the disjoint union $\coprod_\cq G/G_\gs$, i.e., it is the set of pairs $(gG_\gs, \gs)$, where $\gs\in \cq$ and $g\in G$. (Here we are identifying $\psi_\gs(G_\gs)$ with $G_\gs$.) The order relation is defined by inclusion of cosets, i.e.,
	\[
	(gG_\gs,\gs)<(g'G_{\gs'}, \gs') \iff \gs < \gs' \ \text{and} \ (g')^\minus g \in G_{\gs'}
	\]
The space $D(\vert\cq\vert, \psi)$ is called the \emph{basic construction} or the \emph{development} of $(\vert\cq\vert, \psi)$.  It is the quotient of $G\times \vert\cq\vert$ by the equivalence relation $\sim$ defined by
	\begin{equation}\label{e:equiv}
	(g,x)\sim (g', x') \iff x=x' \ \text{and}\ gG_{\gs(x)}=g'G_{\gs(x)} . 
	\end{equation}
One checks that $D(\vert\cq\vert, \psi)$ is the geometric realization of $D(\cq,\psi)$.
Let $[g,x]$ denote the equivalence class of $(g,x)$. The natural $G$-action on $D(\vert\cq\vert, \psi)$ is induced from the action of $G$ on itself by left translation; the isotropy subgroup at $[g,x]$ is $gG_{\gs(x)}g^\minus$, and projection onto the second factor identifies the orbit space with $\vert\cq\vert$. 
The orbit projection $D(\vert\cq\vert, \psi)\to \vert\cq\vert$ has a a section $i:\vert\cq\vert\to D(\vert\cq\vert,\psi)$ defined by $i:x\mapsto [1,x]$. Thus, $i(\vert\cq\vert)$ is a strict fundamental domain for the $G$-action on $D(\vert\cq\vert,\psi)$. (To say that a closed subspace of a $G$-space is a \emph{strict fundamental domain} means that it intersects each orbit in exactly one point.) Note that this implies that $\vert\cq\vert$ is a retract of $D(\vert\cq\vert,\psi)$ (the orbit projection is the retraction).

\begin{remark}\label{r:strict}
Suppose $G$ acts on a space $D$ with a strict fundamental domain $Y$ so that the stratification of $Y$ by closures of points of the same orbit type $\{Y^\gs\}_{\gs\in \cq}$ is indexed by some poset $\cq^{opp}$ and so that for any point $y\in Y_\gs$, the isotropy subgroup $G_y$ contains $G_\gs$. This gives the data for a simple complex of groups $G\cq$ over $\cq$, where the local group $G_\gs$ is the isotropy subgroup at a generic point of $Y^\gs$. There is a canonical simple morphism $\psi= (\psi_\gs)$ from $G\cq$ to $G$ corresponding to the inclusions $G_\gs\hookrightarrow G$. As before, one defines the \emph{development} of $Y$ with respect to $\psi$ by $D(Y,\psi)=(G\times Y)/\sim$, where the equivalence relation $\sim$ is defined as in \eqref{e:equiv}. Moreover, the inclusion $Y\hookrightarrow D$ extends to a $G$-equivariant homeomorphism $D(Y,\psi) \to D$ (see \cite{bh}*{Prop.\,12.20\,(1), II.12}). 
In other words, $D$ is determined by the strict fundamental domain $Y$, the group $G$ and the simple morphism from $G\cq$ to $G$. (In previous work of the first author, e.g., \cite{dbook}*{\S 5.1}, the basic construction is denoted by $\cu(G,Y)$ rather than by $D(Y,\psi)$.)
\end{remark}

Next, we recall a basic lemma, which can be found as \cite{bh}*{Prop.\,12.20\,(4), II.12} or \cite{serre}*{Thms.\,6, 10, pp.~32, 39}. 

\begin{lemma}\label{l:1connected}
Let $G\cq$ be a developable simple complex of groups over $\cq$. Let $\psi$ be a simple morphism from $G\cq$ to a group $G$ and let $\h{\psi}:\gq \to G$ be the induced homomorphism. Suppose $\vert\cq\vert$ is simply connected. Then $D(\vert\cq\vert,\psi)$ is simply connected if and only if $\h{\psi}$ is an isomorphism. (In particular, if $G=\gq$, then $D(\vert\cq\vert,\gi)$ is simply connected.)
\end{lemma}

\begin{remark}\label{r:fund}
Although we will not define the concepts of the ``universal cover'' or the ``fundamental group'' for a general complex of groups, we will give definitions for simple complexes of groups which agree with the more general definitions in \cite{bh}. Put $G=\gq$ and let $\gi:G\cq\to G$ be the canonical simple morphism. First suppose $\vert\cq\vert$ is simply connected. Then, since $D(\vert \cq\vert, \gi)$ is simply connected (by Lemma~\ref{l:1connected}),  $D(\vert\cq\vert,\gi)$ is the universal cover of $G\cq$ and $\pi_1(G\cq)=G$ (cf.\,{bh} \cite{bh}*{III.$\cac$.3.11(1), p.\,551}). 
If $\vert\cq\vert$ is not connected, then each component of $\vert\cq\vert$ gives its own simple complex of groups and can be treated separately. So, suppose $\vert\cq\vert$ is connected, but not necessarily simply connected. Put $\pi=\pi_1(\vert\cq\vert)$. The poset structure on $\vert\cq\vert$ lifts to a poset structure $\cp$ on the universal cover of $\vert\cq\vert$ (so that the universal cover of $\vert\cq\vert$ is $\vert\cp\vert$). The group of deck transformations $\pi$ acts on $\vert\cp\vert$ and on $\cp$. Let $p:\cp\to \cq$ be the projection. We have a simple complex of groups $G\cp$ defined by $G_{\tilde{\gs}}=G_{\gs}$ where $\tilde{\gs}\in \cp$ lies above $\gs$. Put $\ti{G}=\lim G\cp$ and let $\ti{\gi}:G\cp\to \ti{G}$ be the canonical simple morphism. By Lemma~\ref{l:1connected}, the basic construction $D(\vert\cp\vert,\ti{\gi})$ is simply connected. The group $\pi$ acts on $G\cp$ and hence, on $\ti{G}$. So, the semidirect product $\ti{G} \rtimes \pi$ acts on $D(\vert\cp\vert,\ti{\gi})$ with orbit space $\vert\cp\vert/\pi= \vert\cq\vert$. Therefore, $D(\vert\cp\vert,\ti{\gi})$ is the universal cover of $G\cq$ and $\pi_1(G\cq)=\ti{G} \rtimes \pi$.
\end{remark}

\paragraph{Aspherical realizations.} The classifying space of a discrete group $H$ is denoted $BH$. The universal cover of $BH$ is denoted  $EH$. A simple complex of groups $G\cq$ gives the data for a poset of spaces $\{BG_\gs, \ol{\phi}_{\gs\gt}\}$, where $\gt<\gs\in \cq$ and where $\ol{\phi}_{\gs\gt}:BG_\gt\to BG_\gs$ is the map induced by the monomorphism $\phi_{\gs\gt}:G_\gt \to G_\gs$. Using this data we can glue together the disjoint union of spaces $\coprod \vert\cq\vert^\gs \times BG_\gs$ by using iterated mapping cylinders. For each $\gt<\gs$, $\vert\cq\vert^\gs$ is a subcomplex of $\vert\cq\vert^\gt$ and one glues the subspace $\vert\cq\vert^\gs \times BG_\gt$ of $\vert\cq\vert^\gt \times BG_\gt$ to $\vert\cq\vert^\gs \times BG_\gs$ via the map:
	\[
	\mathbb{I}\times\ol{\phi}_{\gs\gt}:\vert\cq\vert^\gs \times BG_\gt\to \vert\cq\vert^\gs \times BG_\gs,
	\]
where $\mathbb{I}$ denotes the identity map on $\vert\cq\vert^\gs$. The resulting space $BG\cq$ is the \emph{aspherical realization} of $G\cq$. It is well-defined up to homotopy equivalence. If $\vert \cq\vert$ is connected and simply connected, then it follows from van Kampen's Theorem that $\pi_1(BG\cq)= \gq$. We note that   
	\begin{equation}\label{e:dim}
  	 \dim BG\cq=\sup\{ (\dim BG_\gs + \dim \vert \cq \vert^\gs ) \mid \gs\in \cq\}
  	 \end{equation}

\begin{Prop}\label{p:developable}
Suppose $G\cq$ is a simple complex of groups over $\cq$ with $\vert\cq\vert$ simply connected. Let $G=\gq$ and let $D=D(\vert\cq\vert, \gi)$ be the basic construction. When $G\cq$ is developable, $BG\cq$ is homotopy equivalent to the Borel construction $EG\times_G D$.
\end{Prop}

\begin{proof}
Projection on the second factor induces a projection $p:EG\times_G D \to D/G = \vert\cq\vert$ so that the inverse image of the vertex $\gs$ is homotopy equivalent to $BG_\gs$. This uses the fact that $G\cq\to G$ is injective on local groups, otherwise, $p^\minus(\gs)$ is homotopy equivalent to $B\ol{G}_\gs$, where $\ol{G}_\gs$ means the image of $G_\gs$ in $G$. It follows that $EG\times_G D$ is an aspherical realization of $G\cq$.
\end{proof}

Recall that the $K(\pi,1)$-Question for $G\cq$ is the following.
\begin{Kquestion}
Is $BG\cq$ aspherical?
\end{Kquestion}

\begin{corollary}\label{p:aspherical} Suppose $G\cq$ is a developable simple complex of groups with $\vert\cq\vert$ simply connected.  Then the $K(\pi,1)$-Question for $G\cq$ has a positive answer if and only if $D$ is contractible.
\end{corollary} 

\begin{proof}
By Proposition~\ref{p:developable}, $BG\cq \sim EG\times_G D$. Hence, $BG\cq$ is aspherical if and only if the universal cover of $EG\times_G D$ is contractible. This universal cover is $EG\times D$, which yields the corollary.
\end{proof}

\begin{remark}\label{r:retract}
Since $\vert\cq\vert$ is a retract of $D$, a necessary condition for $D$ to be contractible is that $\vert\cq\vert$ is contractible.
\end{remark}
 
In the construction of $BG\cq$ we can assume that for each $\gs$, $\dim BG_\gs=\gdim G_\gs$. 
There is a simple condition which implies that the maximum value of the quantity $\dim BG_\gs + \dim \vert\cq\vert^\gs$ in \eqref{e:dim} occurs when $\gs$ is a maximal element of $\cq$, i.e., when $\vert\cq\vert^\gs$ is a point. It is:
  \begin{equation}\label{e:gdim}
  \gdim G_\gs >\gdim G_\gt,\quad \text{whenever $\gs>\gt$}.
  \end{equation}
Since a $k$-simplex in $\vert\cq\vert^\gt$ corresponds to a chain $\gt=\gt_0<\cdots <\gt_\gk$, condition \eqref{e:gdim} implies $\gdim G_{\gt_k}\ge \gdim G_{\gt}+k \ge \gdim G_{\gt} +\dim \vert \cq\vert^{\gt}$; so the maximum value occurs when $\gs$ is maximal.

\begin{Prop}\label{p:gdim}
Suppose  $G\cq$ is a developable  simple complex of groups and that the $K(\pi,1)$-Question for $G\cq$ has a positive answer. Then
  \[
  \gdim G\le \dim BG\cq =\sup_{\gs\in \cq} \{ \gdim G_\gs + \dim \vert\cq\vert^\gs\}.
  \]
If \eqref{e:gdim} holds, then
  \(
  \gdim G =\sup_{\gs\in \cq} \{\gdim G_\gs\}.
  \)
\end{Prop}

\begin{proof}
Since $\gdim G\le \dim BG\cq$, the first formula follows from \eqref{e:dim}. So, if condition \eqref{e:gdim} holds, $\gdim G \le\sup \{\gdim G_\gs\mid \gs\in \cq\}$. Since $EG/G_\gs$ is a model for $BG_\gs$, $\gdim G\ge \gdim G_\gs$, so the previous inequality must be an equality. 
\end{proof}

\subsection{Examples where $\cq$ is a poset of simplices}\label{ss:SL}
Given a simplicial complex $L$, its poset of simplices (including the empty simplex) is denoted by $\cs(L)$. Its geometric realization $\vert\cs(L)\vert$ is the cone on the barycentric subdivision of $L$. (The vertex corresponding to the empty simplex is the cone point.)
Most of this paper concerns simple complexes of groups over posets $\cq$ of the form $\cs(L)$. Given such a simple complex of groups $G\cs(L)$, we often will write $G$ for $\gsl$ and $D$ for the basic construction $D(\vert\cs(L)\vert, \gi)$. In this subsection we introduce  examples of main interest of such complexes of groups coming from Coxeter groups, Artin groups, and graph products.

\begin{example}\label{ex:cox}(\emph{Coxeter groups}) 
Suppose $(W,S)$ is a Coxeter system (cf. \cite{bourbaki} or \cite{dbook}). This means that $W$ is a group, that $S$ is a distinguished set of generators and that $W$ has a presentation of the form $$W := \langle s_i \in S| s_i^2 = (s_is_j)^{m_{ij}} = 1 \rangle$$ For any subset $T$ of $S$, the subgroup generated by $T$ is denoted $W_T$ and called the \emph{special subgroup} corresponding to $T$. It is a standard fact that $(W_T,T)$ also is a Coxeter system (cf. \cite[pp.12-13]{bourbaki}). The subset $T$ is \emph{spherical} if $W_T$ is finite, in this case $W_T$ is a \emph{spherical special subgroup}. Let $\cs(W,S)$ denote the poset of spherical subsets of $S$. There is a simplicial complex $L$ ($=L(W,S)$), called the \emph{nerve} of $(W,S)$. Its vertex set is $S$ and a subset $T\le S$ spans a simplex of $L$ if and only if $T$ is spherical. Thus, $\cs(L)=\cs(W,S)$, the poset of spherical subsets. This gives a simple complex of groups over $\cs(L)$, denoted $W\cs(L)$, and called the \emph{complex of spherical Coxeter groups} (or simply the \emph{Coxeter complex}). The local group $W_\gs$ at a simplex $\gs\in \cs(L)$ is the spherical subgroup generated by the vertices of $\gs$. The direct limit of $W\cs(L)$ is the Coxeter group $W$. 
\end{example}

\begin{example}\label{ex:art}(\emph{Artin groups}). Given a Coxeter system $(W,S)$ there is an associated \emph{Artin group} $A$. 
This group has one generator $x_s$ for each $s\in S$; its relations are the \emph{braid relations}:
	\[
	\underbrace{x_sx_t\ \cdots}_{m_{st} \text{ terms}} = \underbrace{x_tx_s\ \cdots}_{m_{st} \text{ terms}},
	\]
where both sides of the equation are alternating words in $x_s$ and $x_t$, where 
$\{s,t\}$ ranges over the edges of the nerve $L$, and where $m_{st}$ denotes the order of $st$ in $W$. The matrix $(m_{st})$ is called the \emph{Coxeter matrix}. For any $T\le S$, let $A_T$ denote the Artin group corresponding to the Coxeter system $(W_T,T)$. As with Coxeter groups, $A_T$ can be identified with the special subgroup of $A$ generated by $\{x_t\}_{t\in T}$. (When $T$ is a spherical subset of $S$, this is proved in \cite[Th\`eor\'eme 4.13 (iii)]{deligne} and in general in \cite[Theorem 4.14]{lek}.) The Artin group $A_T$ is \emph{spherical} if $T$ is a spherical subset of $S$. As with Coxeter groups, let $\cs(L)$ be the poset of spherical subsets of $S$. The spherical special subgroups of $A$ give a simple complex of groups $A\cs(L)$ called the \emph{Artin complex}. When $\gs$ is a simplex of $L$ and $T=\vertex \gs$, we shall often write $A_\gs$ instead of $A_T$. It is clear that the direct limit $\lim{A\cs(L)}$ is $A$. Since $A_\gs$ is isomorphic to a subgroup of $A$, $A\cs(L)$ is developable. Note that if $L'$ is any subcomplex of the full simplex on $S$ with the same $1$-skeleton as $L$, then $\lim A\cs(L')=A$.

By Deligne's Theorem in \cite{deligne}, any spherical Artin group $A_\gs$ has a classifying space $BA_\gs$ which is a finite CW complex of dimension one greater than $\dim \gs$. (There is a specific model for $BA_\gs$ called the \emph{Salvetti complex}, see \cite{cd2} or Subsection \ref{ss:hyp}). The ``$K(\pi,1)$-Conjecture'' for Artin groups is the conjecture that the $K(\pi,1)$-Question has a positive answer for any Artin poset $A\cs(L)$. By Proposition~\ref{p:aspherical} this is equivalent to the conjecture that $D(\vert A\cs(L)\vert, A)$ is contractible. It is the most important unsolved problem concerning general Artin groups. A detailed discussion can be found in \cite{cd1}, where the conjecture is proved whenever $L$ is a flag complex (we generalize this in Theorem~\ref{t:flag} below). 
\end{example}

The discussion in Example~\ref{ex:art} yields the following proposition.

\begin{Prop}\label{p:gdartin}
Let $A$ be an Artin group such that the nerve $L$ of its associated Coxeter system is $d$-dimensional. If the $K(\pi,1)$-Question for the associated Artin poset has a positive answer, then $\gdim A= d+1$. 
\end{Prop}

\begin{proof}
Since the aspherical realization $BA\cs(L)$ is formed by gluing together the $BA_\gs\times \vert\cs(L)\vert^\gs$ with $\gs\in \cs(L)$ and since $\dim BA_\gs=\dim \gs +1$, we have $\dim BA\cs(L)=d+1$. Since $BA\cs(L)$ is a model for $BA$, $\gdim A\le d+1$. On the other hand, for any $d$-simplex $\gs\in \cs(L)$, the spherical Artin group $A_\gs$ contains a free abelian subgroup of rank $d+1$ (e.g., see \cite{dh16}). So, $\cd A_L\ge d+1$.  Hence, $d+1\ge \gdim A_L\ge \cd A_L \ge d+1$; so all inequalities are equalities.
\end{proof}

\begin{definition}A simplicial complex $L$ is a \emph{flag complex} if it satisfies the following: if $T$ is any finite set of vertices of $L$ which are pairwise connected by edges, then $T$ spans a simplex of $L$. A simplicial graph $L^1$ determines a flag complex $L$: the simplices of $L$ are the cliques in $L^1$. (This is also called the ``clique complex'' of $L^1$.) 
\end{definition}

\begin{definition}\label{d:graphprdt}
Suppose $L^1$ is a simplicial graph with vertex set $V$ and edge set $E$. Let $\{G_v\}_{v\in V}$ be a collection of groups indexed by $V$. 
The \emph{graph product} of the $G_v$, denoted $\prod_{L^1} G_v$, is the quotient of the free product of the $G_v$, $v\in V$, by the normal subgroup generated by all commutators of the form, $[g_v,g_w]$, where $\{v,w\}\in E$, $g_v\in G_v$ and $g_w\in G_w$.
\end{definition}

\begin{example}\label{ex:gprdt} (\emph{The graph product complex}).
Suppose $\prod_{L^1} G_v$ is a graph product and that $L$ is the flag complex determined by $L^1$. There is a simple complex of groups $G\cs(L)$ over $\cs(L)$ called the \emph{graph product complex}. It is defined by putting $G_\gs$ equal to the direct product,
	\[
	\prod_{v\in \vertex \gs} G_v,
	\]
for each simplex $\gs \in L$ and letting $\phi_{\gs\gt}:G_\gt\to G_\gs$ be the natural inclusion whenever $\gt < \gs$. It is immediate that $\prod_{L^1} G_v=\gsl$.
 
\end{example}

\begin{remark}\label{r:poly}
Another approach to the graph product complex is given in \cite{d98} or \cite{davis12} using the notion of a ``polyhedral product''.
Suppose $L$ is a simplicial complex with vertex set $V$ and that we are given a collection of pairs of spaces $(\ux,\uy)=\{(X_v,Y_v)\}_{v\in V}$ together with a choice of basepoint $*_v \in Y_v$.  Let $\prod_{v\in V} X_v$ denote the subspace of the Cartesian product consisting of all $V$-tuples $(x_v)_{v\in V}$ such that $x_v=*_v$ for all but finitely many $v$.  Given a $V$-tuple, $x=(x_v)$, put
	\(
	\gs(x)=\{v\in V\mid x_v\in X_v-Y_v\}.
	\)
The \emph{polyhedral product} $(\ux, \uy)^L$ is the subspace of $\prod_{v\in V} X_v$ consisting of all $x$ such that $\gs(x)$ is a simplex in $L$. (Usually, we only shall be concerned with the case where each $Y_v$ equals the basepoint $*_v$, in which case the notation will be simplified to $\ux^L$.)
For any $\gt\in \cs(L)$, the set of all $x$ with $\gs(x)\ge\gt$, can be identified with the product, $X^\gt :=\prod_{v\in \vertex\gt} X_v$. Thus, $\ux^L$ is the union of the $X^\gt$, with $\gt\in \cs(L)$. 

Next let $L^1$ be any simplicial graph and $L' $ any simplicial complex with $1$-skeleton $=L^1$. 
Let $G$ be the graph product $\prod_{L^1} G_v$ and $H=\prod_{v\in V} G_v$ be the direct product. Let $\gi:G\cs(L^1)\to G$ and $\psi:G\cs(L^1)\to G$ be the natural simple morphisms.
Consider the polyhedral product:
	\[
	\mathcal{Z}(L')=(\underline{\cone G_v}, \underline{G_v})^{L'}.
	\] 
Then $\mathcal{Z}(L')$ is locally isomorphic to a product of cones on discrete sets. 
Moreover, $\mathcal{Z}(L')$ can be identified with the development $D(\vert\cs(L')\vert,\psi)$. The fundamental group of $\mathcal{Z}(L')$ can be identified with the kernel of the natural epimorphism $G\to H$. So, the $H$-action on $\mathcal{Z}(L')$ lifts to a $G$-action on the universal cover $\tilde{\mathcal{Z}}(L')$ with the same strict fundamental domain. It follows from Remark~\ref{r:strict} that $\tilde{\mathcal{Z}}(L')= D(\vert\cs(L')\vert,\gi)$.  When $L'=L$, $\tilde{\mathcal{Z}}(L)$ can be identified with the standard realization of a right-angled building (cf.\ \cite{davis12}*{\S 2.2}).  Therefore, $\tilde{\mathcal{Z}}(L)$ has a $\cat(0)$ structure and hence, is contractible.  A corollary to these observations is the following result.
\end{remark}

\begin{Prop}\label{p:graphprdt}\textup{(cf. \cite{davis12}*{Theorem~2.22} and \cite{dk})}.
Suppose $G=\prod_{L^1} G_v$ is the graph product of nontrivial groups over a simplicial graph $L^1$, and let $L$ be the  flag complex determined by $L^1$..  Then the $K(\pi,1)$-Question for the graph product complex $G\cs(L)$ has a positive answer.
\end{Prop}

Since $BG_1\times BG_2$ is a model for $B(G_1\times G_2)$, it is obvious that
	\[
	\gdim(G_1\times G_2)\le \gdim G_1 +\gdim G_2.
	\]
On the other hand, by using certain torsion-free subgroups of $\racg$s, Dranishnikov \cite{dran} showed that there are groups $G_1$ and $G_2$ for which the inequality is strict (cf.\ \cite[Example 8.5.9]{dbook}). Hence, for $G_\gs=\prod_{v\in \gs}G_v$, we have that $\gdim G_\gs\le \sum_{v\in \gs} \gdim G_v$ and the inequality can be strict. 

A corollary to Proposition~\ref{p:graphprdt} is the following calculation of the geometric dimension of any graph product of groups.

\begin{corollary}\label{c:gprdt}
Suppose $G=\prod_{L^1} G_v$ is the graph product of nontrivial groups over a simplicial graph $L^1$. Let $L$ be the flag complex determined by $L^1$. Then $\gdim G = \sup\{\gdim G_\gs \mid \gs\in \cs(L)\}$.
\end{corollary}

\subsection{Flag complexes and the $K(\pi,1)$-Question}\label{ss:flag}

We will use the following result from \cite{serre}.
\begin{Prop}\label{p:serre}\textup{(\cite{serre}*{Prop.~3, p.~6}).} Let $G_i$ be a collection of groups with common subgroup $A$ and let
 $\aster_A\, G_i$ denote the amalgamated product. Let $H_i\subset G_i$ be subgroups and suppose the intersection $B=H_i\cap A$ is independent of $i$. Then the natural homomorphism $\aster_B\ H_i \to \aster_A\, G_i$ is injective.
\end{Prop}
As usual, $G\cs(L)$ is a simple complex of groups. For any full subcomplex $L'$ of $L$, put $G_{L'}=\lim G\cs(L')$ and let $\ol{G}_{L'}$ denote the image of $G_{L'}$ in $G_L$. For any $\gs\in \cs(L)$, let $\ol{G_\gs}$ denote the image of $G_\gs$ in $G_L$. (If $G\cs(L)$ is developable, then $G_\gs\to \ol{G}_\gs$ is an isomorphism.) 
The \emph{intersection of local groups condition} for simplices $\gs$ and $\gt$ of $L$ is the following:
	\begin{equation}\label{e:simplexinter}
	\ol{G_{\gs}}\cap\ol{G_\gt} =\ol{G_{\gs\cap\gt}}. 
	\end{equation}
	
\begin{lemma}\label{l:inject}
Suppose $G\cs(L)$ is developable and that \eqref{e:simplexinter} holds for all $\gs,\gt \in \cs(L)$. If $L$ is a flag complex, then for any full subcomplex $L'$ of $L$, the natural map $G_{L'}\to G_L$ is an injection.
\end{lemma}

We first prove this in a special case where $L$ is the cone on another flag complex. We will use the notation $\St(u)$ for the cone and $\Lk(u)$ for the base of the cone; regarding $u$ as the cone vertex.
\smallskip

\noindent
\textbf{Special Case:} We will prove the following:
\begin{enumeratei}
\item
The natural map $G_{\Lk(u)}\to G_{\St(u)}$ is injective.
\item
Suppose $L_0$, $L_1$ are full subcomplexes of $\Lk(u)$ with $L_0\subset L_1$. Then
$\displaystyle{G_{L_1} \cap G_{L_0 * u}= G_{L_0}}$.
\end{enumeratei}

\begin{proof}[Proof of the Special Case]
The proof is by induction on the number of vertices of $\St(u)$. When $\St(u)$ is a simplex, (i) and (Ii) follow from \eqref{e:simplexinter}. If $\St(u)$ is not a simplex, then, since $\Lk(u)$ is a flag complex, there is a vertex $v$ in $\Lk(u)$ which is not connected by an edge to some other vertex of $\Lk(u)$. Denote the star of $v$ in $\St(u)$ simply by $\St(v)$. So, $\St(v)$ is a proper subcomplex of $\St(u)$. By induction on the number of vertices in the cone, (i) holds for $\St(v)$, i.e., 
\(
G_{\Lk(v)}\to G_{\St(v)}
\)
is injective, where $\Lk(v):=\Lk(v,\St(v))$. By induction on the number of vertices, the natural maps,
\[
G_{\Lk(u)-v}\to G_{\St(u)-v},\ \ G_{\St(v,\Lk(u))} \to G_{\St(v)}, \ \ \text{and}\ \ G_{\Lk(v, \Lk(u))} \to G_{\Lk(v)},
\]
are injective. We have decompositions as amalgamated products:
	\begin{align*}
	G_{\St(u)}&=G_{\St(u)-v} *_{G_{\Lk(v)}} G_{\St(v)}\\
	G_{\Lk(u)}&=G_{\Lk(u,\St(u)-v)} *_{G_{\Lk(v,\Lk(u))}} G_{\Lk(u,\St(v))}
	\end{align*}
By induction on the number of vertices, we have that (ii) holds in two cases: first with $L_1 = \Lk(u,\St(u)-v) = \Lk(u) - v$ and $L_0=\Lk(v,\Lk(u))$ and second with $L_1 = \St(v)$ and $L_0=\Lk(v,\Lk(u))$. This yields
	\begin{align*}
	G_{\Lk(u,\St(u)-v)} \cap G_{\Lk(v,\Lk(u)) * u} &= G_{\Lk(v,\Lk(u))},\\
	G_{\St(v)} \cap G_{\Lk(v,\Lk(u)) * u} &= G_{\Lk(v,\Lk(u))}.
	\end{align*}
Applying Proposition~\ref{p:serre}, we see that 
\(
G_{\Lk(u)} \to G_{\St(u)}
\)
is injective.
\end{proof}

\begin{proof}[Completion of Proof]
We can now complete the proof of Lemma~\ref{l:inject}. The argument is again by induction on the number of vertices. Let $L'$ be a full subcomplex of $L$. If $L = L' \ast v$, we are done by the special case. Otherwise, there is a vertex $v$ such that $L\neq \St(v)$ and $v \notin L'$. There is an amalgamated product decomposition:
	\[
	G_L=G_{L-v}*_{G_{\Lk(v)}} G_{\St(v)}.
	\]
By inductive hypothesis for the proper subcomplexes $L-v$ and $\St(v)$, we have inclusions $G_{\Lk(v)}\to G_{L-v}$ and $G_{\Lk(v)}\to G_{\St(v)}$. Since both factors inject into the amalgamated product, $G_{L-v}\to G_L$ is an injection. Continuing by deleting one vertex at a time, we see the same holds for $G_{L'}$. 
\end{proof}

\begin{figure}
   \centering
   \begin{tikzpicture}[scale = .7]
   
   \draw[thick, fill = gray, opacity = .3] (0,0) circle (4cm);
   \draw[thick] (0,0) circle (4cm);
   \draw[fill=gray!60] (0,0) -- +(60:4) arc (60:-60:4);
   
   \node[left] at (0,0) {$u$};
   \draw[thick] (0,0) -- (4,0);
  \node[right] at (4,0) {$v$};
  \node[left] at (-4,4) {$\Lk(u)$};
  \node at (2,-1.7) {$\St(v)$};
  \node at (-2,0) {$\St(u)$};
  \draw[->][thick] (-4.1,4) .. controls (-3.4,3.6) .. (-2.7,3.1);
  
   \draw[thick] (0,0) -- (2,3.42);
    \draw[thick] (0,0) -- (2,-3.42);
    
      \node[right] at (4,2.82) {$\Lk(u,\St(v)) = \St(v, \Lk(u))$};
      
      \draw[->][thick] (4,2.82) .. controls (3.5,2.82) .. (2.95,2.82);
      
   \node[above] at (1,5) {$\Lk(u, \Lk(v))$};  
   
    \draw[->][thick] (1,5.2) .. controls (1.3,4.9) .. (2,3.6); 
   
   \end{tikzpicture}
  \caption{Setup for Lemma 1.19}
 \end{figure}
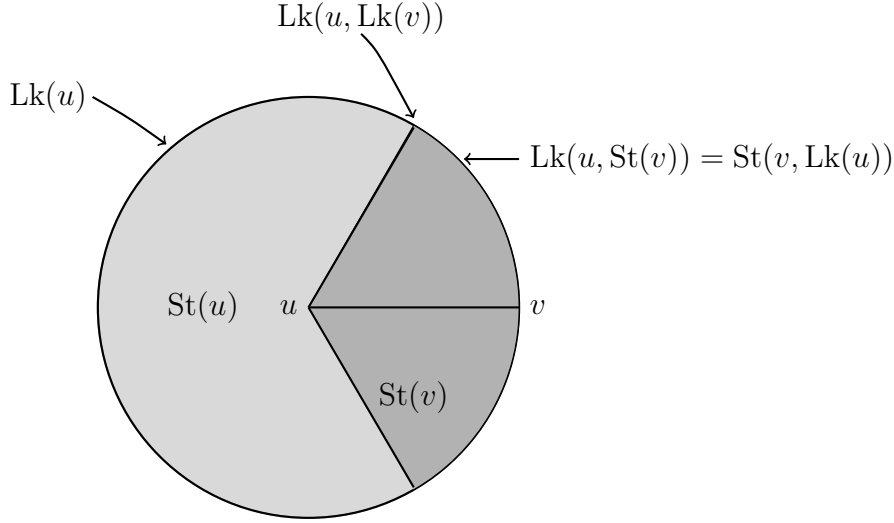
 
\begin{Theorem}\label{t:flag}\textup{(cf.~\cite{cd1}*{Remark on p.~619}).} Suppose $L$ is a flag complex, that $G\cs(L)$ is developable, and that condition \eqref{e:simplexinter} holds for all $\gt$, $\gs$ in $\cs(L)$.
Then the $K(\pi,1)$-Question for $G\cs(L)$ has a positive answer.
\end{Theorem}

\begin{proof}
 The proof is by induction on the number of maximal simplices in $L$. The base case is when $L$ is a single simplex $\gs$. Then $G_\gs=\lim G\cs(\gs)$ and $BG\cs(\gs)\sim BG_\gs$. So, the answer to the $K(\pi,1)$-Question for $G\cs(\gs)$ is positive.  

Suppose $L$ is not a simplex. Since $L$ is a flag complex, there are distinct vertices $v_1, v_2$ in $V=\vertex L$ so that $v_1$ and $v_2$ are not joined by an edge of $L$. Let $L_1$, $L_2$ and $L_0$ be the full subcomplexes of $L$ spanned by $V-\{v_1\}$, $V-\{v_2\}$ and $V-\{v_1,v_2\}$, respectively, and let $G_1$, $G_2$ and $G_0$ be the respective direct limits.  Since $G_L$ is a direct limit, it is the amalgamated product $G_L=G_1*_{G_{0}} G_2$. By Lemma~\ref{l:inject}, the natural maps from $G_0\to G_i$ and $G_i\to G_L$ are injections for $i=1,2$.
A lemma of Whitehead states that if two aspherical complexes are glued together along an aspherical subcomplex so that the fundamental group injects into either side, then the result of the gluing is also aspherical (for example, see \cite[Thm.~E.1.15]{dbook}). (This is a special case of theorem that we are proving.) By inductive hypothesis, $BG\cs(L_i) \sim BG_i$ and, by its definition, $BG\cs(L)$ is the union of $BG\cs(L_1)$ and $BG\cs(L_2)$ along $BG\cs(L_0)$. Hence, Whitehead's Lemma shows that $BG\cs(L)\sim BG_L$.
\end{proof}

\begin{remark}
Theorem~\ref{t:flag} shows that if $G\cs(L)$ is developable, \eqref{e:simplexinter} holds and $L$ is a flag complex, then $D(\vert\cs(L)\vert,G_L)$ is contractible. Moreover, $D(\vert\cs(L')\vert,G_{L'})$ is contractible for any full subcomplex $L'\le L$. For graph product complexes, these developments are right-angled buildings and hence, are $\cat(0)$. This raises the question of nonpositive curvature in the context of Theorem~\ref{t:flag}. This is equivalent to the question of whether the link of each simplex in $D(\vert\cs(L)\vert,G_L)$ is $\cat(1)$, (see \cite{bh}). For bounded curvature to make sense, we should first put a piecewise spherical structure on $L$. Since $L$ is already assumed to be a flag complex we might as well assume that each simplex in $L$ is an all right spherical simplex so that $D(\vert\cs(L)\vert,G_L)$ becomes a cube complex, see \cite{dbook}. The link, $\Lk_\gs$, corresponding to a nonempty simplex $\gs$ in $L$ is the development of the simplex of groups $G\cs(\partial \gs)$ with respect to the natural simple morphism $\psi:G\cs(L) \to G_\gs$, so the question becomes whether $\Lk_\gs=D(\vert\cs(\partial \gs)\vert, \psi)$ is $\cat(1)$. By Gromov's Lemma (cf.\ \cite{dbook}*{Appendix I.6, pp.516-517}), this is equivalent to the question of whether it is a flag complex. Note that this is implied by the condition that the development of the simplex of groups, $D(\vert\cs(\partial \gs)\vert,G_{\partial \gs})$, is a flag complex. (Here $G_{\partial \gs}=\lim{G\cs(\partial\gs)}$).

Similar considerations led R.~Charney and the first author to conjecture in \cite{cd1} that the development of any Artin complex in Example~\ref{ex:art} can be given $\cat(0)$ structure. This would imply the $K(\pi,1)$-Conjecture for all Artin groups. The piecewise spherical structure on $L$ should be the natural one in which $\gs$ is isometric to a fundamental spherical simplex in the round sphere on which $W_\gs$ acts as reflection group; in other words, the spherical simplex which makes the spherical realization of the Coxeter complex for $W_\gs$ into a round sphere. In the Artin complex, $\Lk_\gs$ is the Deligne complex for $A_\gs$. So, the conjecture of \cite{cd1} is that the natural piecewise spherical metric on the Deligne complex for $A_\gs$ is $\cat(1)$. 
\end{remark}

\subsection{$\cq$ is the poset of nonempty simplices in $L$}\label{ss:cq=cs}
Suppose $L$ is a full subcomplex of another simplicial complex $\bar L$. Then every simplex $\bar \gs \in \cs(\bar L)$ can be decomposed as a join $\bar \gs=\gs*\gt$, where $\gs=\bar \gs \cap L$ and $\gt$ is a simplex whose vertices lie in $\bar L-L$. Suppose $G\cs(L)$ and $G\cs(\bar L)$ are simple complexes of groups such that $G\cs(L)$ is the restriction of $G\cs(\bar L)$ to $\cs(L)$. Then $G\cs(\bar L)$ is a \emph{trivial extension} of $G\cs(L)$ if $G_{\bar \gs}= G_\gs$, whenever $\bar \gs=\gs *\gt$ decomposes as a join as above. 

When $L$ is a subcomplex of $\bar L$, we can replace $\bar L$ by a subdivision relative to $L$ so that $L$ becomes a full subcomplex of $\bar L$. If $L$ is a flag complex, then the subdivision $\bar L$ can be assumed to be flag. So, any simple complex of groups $G\cs(L)$ over $\cs(L)$ has a trivial extension to a simple complex of groups over $\cs(\bar L)$ as above. The following lemma is clear. 

\begin{lemma}\label{l:trivialext}
Suppose $G\cs(\bar L)$ is a trivial extension of $G\cs(L)$.  Then
\begin{enumeratei}
\item
$\lim{G\cs(\bar L)}=\gsl$.
\item
$G\cs(\bar L)$ is developable if and only if $G\cs(L)$ is developable.
\item
$BG\cs(\bar L)$ is homotopy equivalent to $BG\cs(L)$.
\end{enumeratei}
\end{lemma}

For a given simplicial complex $L$, let $\cs^o(L)$ denote its poset of nonempty simplices, i.e., $\cs^o(L)=\cs(L)_{>\emptyset}$. The geometric realization of $\cs^o(L)$ is the barycentric subdivision of $L$. Let $G\cs(L)$ be a simple complex of groups over $\cs(L)$ and let $G\cs^o(L)$ denote its restriction to $\cs^o(L)$. For example, $\cs^o(\cone L)$ can be identified with $\cs(L)$ and any simple complex of groups $G\cs(L)$ can be identified with a simple complex of groups over $\cs^o(\cone L)$. Since we shall always assume that the local group $G_\emptyset$ attached to the empty simplex is the trivial group,  then simple complexes of groups $G\cs^o(L)$ and $G\cs(L)$ have the same direct limit, which we denote by $G$. The next lemma is immediate.

\begin{lemma}\label{l:SvsSo}
The complex of groups $G\cs^o(L)$ is developable if and only if $G\cs(L)$ is developable.
\end{lemma}

\begin{lemma}\label{l:BGS}
Suppose $L$ is connected. Then the following statements are equivalent
\begin{enumeratei}
\item
$L$ is contractible.
\item
The inclusion $D(\vert\cs^o(L)\vert, G)\hookrightarrow D(\vert\cs(L)\vert, G)$ of developments is a $G$-equivariant homotopy equivalence.
\item 
The inclusion $BG\cs^o(L)\hookrightarrow BG\cs(L)$ of aspherical realizations is a homotopy equivalence.
\item
The $K(\pi,1)$-Questions for $G\cs^o(L)$ and for $G\cs(L)$ have the same answers.
\end{enumeratei}
\end{lemma}

\begin{proof}
By Remark~\ref{r:fund}, the following three conditions are equivalent:
\begin{itemize}
\item
$L$ is simply connected.
\item
$D(\vert\cs^o(L)\vert, G)$ is simply connected.
\item
$\pi_1(BG\cs^o(L))=G$.
\end{itemize}
To simplify notation, put $D=D(\vert\cs(L)\vert, G)$ and $D^o=D(\vert\cs^o(L)\vert, G)$. 
Let $C=\vert\cs(L)\vert -\vert\cs^o(L)\vert$ be the open dual cone of the vertex in $\vert\cs(L)\vert$ corresponding to $\emptyset$.  The inverse image of $C$ in $D(\vert\cs(L)\vert, G)$ consists of $G$ copies of $C$.  Hence,
	\[
	H_*(D, D^o)\cong  \bigoplus_G H_*(\cone L,L), \text{ and} 	
	\]
	\[
	H_*(EG\cs(L), EG\cs^o(L))\cong  \bigoplus_G H_*(\cone L,L). 	
	\]
The equivalence of conditions (i), (ii) and (iii) follows.  The equivalence of condition (iv) follows from Corollary~\ref{p:aspherical}. 
\end{proof}

We will use the above lemma in the following way. Given $G\cs(L)$, first embed $L$ as a full subcomplex of a contractible simplicial complex $L_c$ and then take a trivial extension of $G\cs(L)$ to $G\cs(L_c)$. By Lemmas \ref{l:trivialext} and \ref{l:BGS}, there are homotopy equivalences $BG\cs(L)\sim BG\cs(L_c)\sim BG\cs^o(L_c)$, i.e., $BG\cs^o(L_c)$ is another model for $BG\cs(L)$. If $\dim L_c=\dim L$, then  $\dim(\vert\cs^o(L)\vert) = \dim (\vert\cs^o(L)\vert) -1$ and we will sometimes be able to use this to reduce upper bounds on the action dimension or the geometric dimension by $1$. 

\begin{definition}\label{d:EDCE}
A $d$-dimensional simplicial complex $L$ is \emph{equidimensionally, contractibly embeddable} (abbreviated EDCE) if it can be embedded in a contractible complex $L_c$ of the same dimension $d$.
\end{definition}

\begin{remark}
For $d\neq 2$ the condition that $L$ be EDCE is equivalent to the following two conditions:
\[
H_{d-1}(L,\zz)\ \ \text{is free abelian and}\ \ H_d(L,\zz)=0.
\]
By the Universal Coefficient Theorem, the above conditions are equivalent to the condition that $H^d(L,\zz)=0$. Indeed, when these conditions hold, one can use standard methods to attach cells of dimension $\le d$ to kill all the homology of $L$. When $d=2$, we will only end up with an $L_c$ which is acyclic. A conjecture of Kervaire asserts that one cannot kill a nontrivial group by adding the same number of generators and relations. So, in the many situations where Kervaire's Conjecture is known to hold, it is not possible to obtain a contractible complex by adding $1$ and $2$-cells to an acyclic complex.
\end{remark}

\section{Gluing}\label{s:gluing} 
Suppose we are given a collection of manifolds with boundary $\{M_\gt\}_{\gt \in \cq}$ indexed by a poset $\cq$. Further suppose that whenever $\gt<\gs$, we have an embedding $i_{\gt\gs}: M_\gt\hookrightarrow \partial M_\gs$ with trivial normal bundle.  This means, in particular, that whenever $\gs$ is not a minimal element of the poset, $\partial M_\gs$ must be nonempty. 
The poset of manifolds $\{M_\gt\}_{\gt\in \cq}$ is \emph{$n$-dimensional} if for each maximal $\gs$, $\dim M_\gs=n$. Henceforth, we assume this. Let $c(\gt)=n-\dim M_\gt$ and let $D_\gt$ be a disk of dimension $c(\gt)$. The basic idea in this section is that we can glue together the $M_\gt\times D_\gt$ along codimension-zero submanifolds of their boundaries to obtain $M$, an $n$-manifold with boundary. Let $G_\gt:=\pi_1(M_{\gt})$. 
We assume each inclusion $M_{\gt}\hookrightarrow M_{\gs}$ is $\pi_1$-injective, and that $G\cq=\{G_\gt\}_{\gt\in \cq}$ is a simple complex of groups. If $\vert\cq\vert$ is simply connected, then $\pi_1(M)=\gq$ is the direct limit of the $\pi_1(M_\gt)$. So, if each $M_\gt$ is aspherical, $M$ will be a model for $BG\cq$. In practice $\cq$ will be $\cs(L)$ or $\cs^o(L)$ for some $L$.

The main work in Subsection~\ref{ss:thick} is to describe the ``dual disk'' $D_\mu$, for each $\mu\in \cs(L)$ and a decomposition of its boundary sphere into pieces $\T_\mu(\gt)$. For $\gt>\mu$, $M_\mu\times D_\mu$ will be glued onto $M_\gt\times D_\gt$ along $M_\mu\times\T_\mu(\gt)$.
The dual disk is a thickening of the dual cone, $\cone (L_\mu)$, where $L_\mu$ means the normal link of $\mu \in L$ (i.e., $L_\mu$ is the geometric realization of $\cs(L)_{>\mu}$). If $D_\mu$ is a $k$-disk, then it should be attached to other pieces along a submanifold $\T_\mu$ of codimension $0$ in $\partial D_\mu$ ($= S^{k-1}$). The manifold $\T_\mu$ is a thickening of $L_\mu$ in $S^{k-1}$ ($\T_\mu$ is the ``thick link''). The thick link $\T_\mu$ is further decomposed into pieces $\T_\mu(\gt)$ on which the piece corresponding to $\gt$ is to be attached.

\begin{figure}\label{f:thick}
\centering
\begin{tikzpicture}[scale = .8]

 \draw (-4,2) -- (-2,0) -- (2,0) -- (4,2);
  \draw (4,-2) -- (2,0);
    \draw (-4,-2) -- (-2,0);

    \node[below] at (-1.8,0) {$\mu$};
    \node[below] at (-2.8,-1) {$\gs$};
  \node[below] at (0,0) {$\gt$};

         \draw[fill =black] (-2,0) circle (.3mm);
         
         \begin{scope}[yshift = -5cm]
          \draw (-4,2) -- (-2,0) -- (2,0) -- (4,2);
  \draw (4,-2) -- (2,0);
    \draw (-4,-2) -- (-2,0);
    
     \draw[thick] (-3.55,.45) -- (-2.45, 1.55);
        \draw[thick] (3.55,.45) -- (2.45, 1.55);
         \draw[thick]  (2.45, -1.55) -- (3.55,-.45);   
          \draw[thick]  (-2.45, -1.55) -- (-3.55,-.45);

    \node at (-5,0) {$T_\mu(\gs)$};

\draw[->][thick] (-4.3,-0.3) .. controls (-4,-.7) .. (-3.25,-.84);

    \node[below] at (-1.8,0) {$[\mu]$};
        \node[right] at (-2.15,-2) {$[\gs]$};
    
    \draw[->][thick] (-2,-2) .. controls (-2.5,-1.7) .. (-3,-1.1);

  \node[below] at (0,0) {$[\gt]$};

    \draw[fill =black] (0,0) circle (.5mm);
    
         \draw[fill =black] (-3,-1) circle (.5mm);
         \draw[fill =black] (-2,0) circle (.5mm);
 
  \draw[fill = gray, opacity = .3] (0cm,0cm) ellipse[x radius=.5cm,y radius=1cm];
  \draw[thick,dotted] (0,-1) arc[start angle=-90,end angle=90, x radius=0.5, y radius=1];
   \draw[thick] (0,1) arc[start angle=90,end angle=270, x radius=0.5, y radius=1];

 \draw[fill = gray, opacity = .3, rotate around={45:(-3,1)}] (-3,1) ellipse [x radius=.8cm,y radius=.3cm];
   \draw[thick, rotate around={45:(-3,1)}] (-2.2,1) arc[start angle=0,end angle=-180, x radius=0.8, y radius=.3];
   \draw[thick, dotted, rotate around={45:(-3,1)}] (-2.2,1) arc[start angle=0,end angle=180, x radius=0.8, y radius=.3];

  \draw[fill = gray, opacity = .3, rotate around={-45:(-3,-1)}] (-3,-1) ellipse [x radius=.8cm,y radius=.3cm];
   \draw[thick,dotted, rotate around={-45:(-3,-1)}] (-2.2,-1) arc[start angle=0,end angle=-180, x radius=0.8, y radius=.3];
   \draw[thick, rotate around={-45:(-3,-1)}] (-2.2,-1) arc[start angle=0,end angle=180, x radius=0.8, y radius=.3];

   \draw[fill = gray, opacity = .3, rotate around={135:(3,1)}] (3,1) ellipse [x radius=.8cm,y radius=.3cm];
     \draw[thick,dotted, rotate around={135:(3,1)}] (3.8,1) arc[start angle=0,end angle=-180, x radius=0.8, y radius=.3];
   \draw[thick, rotate around={135:(3,1)}] (3.8,1) arc[start angle=0,end angle=180, x radius=0.8, y radius=.3];

  \draw[fill = gray, opacity = .3, rotate around={-135:(3,-1)}] (3,-1) ellipse [x radius=.8cm,y radius=.3cm];
   \draw[thick, rotate around={-135:(3,-1)}] (3.8,-1) arc[start angle=0,end angle=-180, x radius=0.8, y radius=.3];
   \draw[thick, dotted, rotate around={-135:(3,-1)}] (3.8,-1) arc[start angle=0,end angle=180, x radius=0.8, y radius=.3];

\draw plot [smooth] coordinates {(0,1) (1.5,1) (2.43,1.57)};

\draw plot [smooth] coordinates {(0,-1) (1.5,-1) (2.43,-1.57)};

\draw plot [smooth] coordinates {(0,1) (-1.5,1) (-2.43,1.57)};

\draw plot [smooth] coordinates {(0,-1) (-1.5,-1) (-2.43,-1.57)};

\draw[->][thick] (0,2) .. controls (.25,1) .. (0,.5);
  \node[above] at (0,2) {$T(\gt)$};
  
  \draw[->][thick] (-2,2) .. controls (-1.5,1.3) .. (-1.5,.5);
  \node[above] at (-2,2) {$T(\mu)$};
  
    \draw[->][thick] (-2,2) .. controls (-1.5,1.3) .. (-1.5,.5);

  \draw plot [smooth] coordinates {(-3,-2.5) (-2.8, -2) (-2.43,-1.57)};
  
    \draw plot [smooth] coordinates {(3,2.5) (2.8, 2) (2.43,1.57)};
      \draw plot [smooth] coordinates {(-3,2.5) (-2.8, 2) (-2.43,1.57)};
            \draw plot [smooth] coordinates {(3,-2.5) (2.8, -2) (2.43,-1.57)};

\draw plot [smooth] coordinates {(-4.5, 1.2) (-3.55, .42) (-3.4,0) (-3.55, -.42) (-4.5,-1.2)};
\draw plot [smooth] coordinates {(4.5, 1.2) (3.55, .42) (3.4,0) (3.55, -.42) (4.5,-1.2)};
 \end{scope}

\end{tikzpicture}
\caption{A thickening of a graph in $S^3$, and some of the associated dual disks.}
\end{figure}
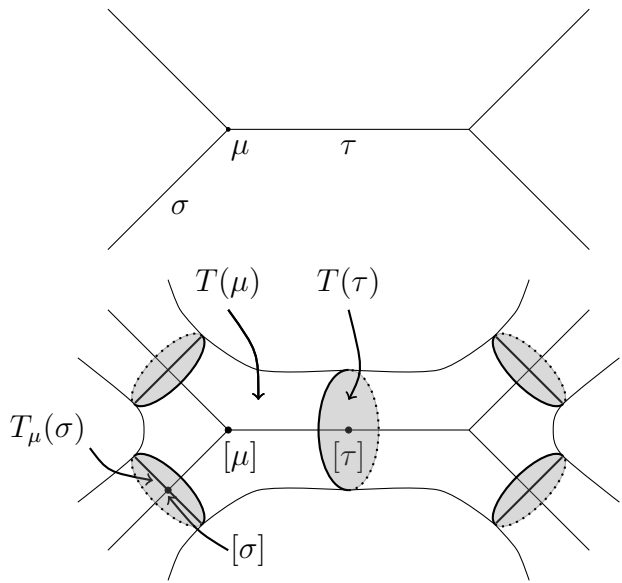

Here is a picture to keep in mind. Suppose $L$ is a simplicial graph, so that $L$ embeds in $S^3$. The dual disk $D_\emptyset$ is $D^4$ and $\T_\emptyset$ is a thickening of $L$ in $S^3$ to a $3$-manifold with boundary. For each vertex $\mu$, $\T_\emptyset(\mu)$ is a $3$-ball neighborhood of $\mu$ and for each edge $\gt$, $\T_\emptyset(\gt)$ is a tubular neighborhood of the edge as a solid cylinder. (See Figure~\ref{f:thick}.)
The dual disk to a vertex $\mu$ is a $2$-disk. If the degree of $\mu$ is $p$, then $L_\mu$ consists of $p$ points and $\T_\mu$ is a thickening to $p$ intervals in $S^1$. The dual cell to an edge (i.e. to a maximal simplex) is a $0$-disk.

The simplicial complex $\vert\cs(L)\vert$ is equal to  $\cone L'$, where $L'$ means the barycentric subdivision of $L$.
So, a $k$-simplex in $\cone L'$ is a chain $\{\gt_0,\dots,\gt_k\}$ of length $k+1$ in $\cs(L)$. We denote the geometric realization of this simplex by $[\gt_0,\dots, \gt_k]$. (We will always write such chains in increasing order, i.e., $\gt_0<\cdots <\gt_k$.) A chain of length $1$, $\{\gt\}$,  is either the cone point $[\emptyset]$ or it corresponds to a vertex $[\gt]$ of $L'$ (thought of as the barycenter of $\gt$). Thus, $[\gt_0,\dots, \gt_k]$ is the $k$-simplex in $\cone L'$ spanned by the vertices $[\gt_0],\dots, [\gt_k]$.  Given a simplex $\ga=[\gt_0,\dots \gt_k]$ in $\cone L'$, its \emph{minimum vertex} is defined by 
	 \(
	 \min\ga=\gt_0.
	 \)	 

To understand the decompositions of dual disks, one first needs to understand the stratification of $\vert\cs(L)\vert$ into ``dual cones.''   In this section we shall often use the notation $K$ (or $K_\emptyset$) instead of $\vert\cs(L)\vert$. Similarly, $K_\gs$ is $\vert\cs(L)\vert_{\ge \gs}$ and called the \emph{dual cone of $\gs$}.  (Here we are reversing the use of superscripts and subscripts from notation in \eqref{e:posetstrata}.)  We will use $\partial K_\gs$ or $L'_\gs$ for $\vert\cs(L)\vert_{>\gs}$, the barycentric subdivision of the link of $\gs$.  Thus, a simplex $\ga=[\gt_0,\dots \gt_k]$ of $K$ lies in $K_\gs$ (resp., $\partial K_\gs$) if and only if $\min\ga\le \gs$ (resp., $\min\ga< \gs$).

\subsection{Thick links}\label{ss:thick}
Suppose $L$ is a finite simplicial complex of dimension $d$.  We recall a method for  thickening  $L$ to a manifold with boundary $\T$.  First, piecewise linearly embed $L$ into a sphere $S^{n-1}$, where $n\ge 2d+2$.  Thus, $L$ is a full subcomplex of some PL triangulation $S$ of $S^{n-1}$.  Denote  the barycentric subdivisions of $L$ and $S$ by $L'$ and $S'$, respectively.  Let $\T$ denote the first derived neighborhood of $L$ in $S$.  In other words, $\T$ is the union of simplices in $S'$ which have nonempty intersection with $L'$.

If $\gr_0<\dots<\gr_k$ is a chain of simplices in $\cs(S)$, then $[\gr_0,\dots,\gr_k]$ denotes the simplex in $S'$ spanned by their barycenters.

For each vertex $\nu$ in $L$, let $\T(\nu)$ denote the closed star of $[\nu]$ in $S'$, i.e.,  
	\[
	\T(\nu)=\bigcup\{\gamma \in S'\mid \nu\le \min \gamma\}.
	\]
Then $\T(\nu)$ is an $(n-1)$-disk.  Moreover, $\T=\bigcup \T(\nu)$, where the union is over all vertices $\nu\in L$.  If $\gt$ is a simplex of $L$, then let $T(\gt)$ denote the normal star of $\gt$ in $S$, i.e.,
	\[
	\T(\gt)=\bigcup\{\gamma \in S'\mid \gt\le \min \gamma\} = [\gt] *\Lk(\gt,S)'.
	\]
So, $T(\gt)$ is the cone on $\Lk(\gt,S)'$. 
Since $\Lk(\gt,S)$ is a sphere of dimension $n-k-2$, where $k=\dim\gt$, we see that $\T(\gt)$ is PL homeomorphic to a $(n-k-1)$-disk.
Moreover, if $\vertex \gt=\{\nu_0,\dots, \nu_k\}$, then $\T(\gt)=\T(\nu_0)\cap \cdots \cap \T(\nu_k)$.  So, $\T$ is an $(n-1)$-manifold with boundary embedded in $S^{n-1}=\partial D^n$

Next we want to apply this construction to links in $L$.  For each simplex $\mu\in \cs(L)$, let $L_\mu=\Lk(\mu,L)$.  We want to thicken $L_\mu$ in a sphere $S_\mu:=S^{c(\mu)-1}$ of an appropriate dimension $c(\mu)-1$ to obtain a manifold with boundary $\T_\mu$ called a \emph{thick link}, embedded in a disk $D_\mu=D^{c(\mu)}$, called the \emph{dual disk}.  Similarly, whenever $\mu<\gt$, we will have $\T_\mu(\gt)$, the normal star of $\gt$ in $S_\mu$.  Thus, $\T_\mu(\gt)$ is a thickening of $\partial K_\mu(\gt)$.   In particular, when $\mu$ is the empty simplex we will have $\T_\emptyset = T$, $D_\emptyset= D^n$ and $\T_\emptyset(\gt)=\T(\gt)$.  

To specify the dimensions of the dual disks
suppose we are given a function $c:\cs(L)\to \nn$ so that 
	\begin{itemize}
	\item \itemEq{\text{For any maximal simplex } \gs \text{ of } L, c(\gs)=0.} 
	\item \itemEq{\text {If } \gt<\gs, \text{then } c(\gt)-c(\gs)\ge 2\codim(\gt,\gs)}.
	\end{itemize}

Put $n=c(\emptyset)$. (If we are using $\cs^o(L)$ instead of $\cs(L)$, then we don't require the second condition for $\gt=\emptyset$.) 
Thus, the dual disk $D_\mu$ will be a thickening of the dual cone $K_\mu$ and the thick link $\T_\mu$ will be a thickening of $\partial K_\mu$ ($=L'_\mu$).

The link $L'_\mu$ naturally is a subcomplex of $L'_\emptyset$ and, in fact, whenever $\mu<\gt$, $L'_\gt$ is a subcomplex of $L'_\mu$.  Similarly, if $S'_\mu$ means the barycentric subdivision of $S^{c(\mu)-1}$, then we  want to arrange that $S'_\gt$ is a subcomplex of $S'_\mu$. This amounts to requiring that $S'_\gt$ is PL embedded in $\Lk(\gt,S^{c(\mu)-1})' \subset S'_\mu$.  For example, if $\mu=\emptyset$ and $L_\emptyset$ is a graph with a vertex $\gt$, then $\dim L'_\emptyset =1$ and $\dim L'_\gt=0$.  We can thicken $L'_\emptyset$ in $S^3=S_\emptyset$, while $L'_\gt$ should be thickened in $S'_\gt=S^1$, see Figure 2. 

Next we want to explain how  dual disks can be regarded as manifolds with corners.
Recall that an $n$-manifold with boundary, $P$, is a \emph{smooth manifold with corners} if it is locally differentiably modeled on the simplicial cone $[0,\infty)^n$. This can be extended to a definition for topological manifolds by requiring that the overlap maps preserve the stratification of $[0,\infty)^n$ by intersections with coordinate subspaces.  The stratification of $[0,\infty)^n$ by its faces then induces a stratification of $P$.   A codimension-one stratum of $P$ is called a \emph{facet}.

\begin{lemma}\label{l:corners}For each $\mu\in \cs(L)$, the dual disk $D_\mu$ is a $c(\mu)$-manifold with corners. The facets are $\{\T_\mu(\gt)\}$ where $\mu$ is a codimension one face of $\gt$ (i.e., where $[\gt]$ is a vertex of $L'_\mu$) together with $\ol{\partial D_\mu - \T_\mu}$.
\end{lemma}

The facet $\ol{\partial D_\mu - \T_\mu}$ is called a \emph{boundary piece}; the other facets are \emph{ordinary facets}.  

\begin{proof}
If $\mu$ is a codimension one face of $\gt$, then $\T_\mu(\gt)\subset \partial D_\mu$ is a disk of codimension one in $D_\mu$. In general, if $\mu$ is a codimension-$k$ face of a simplex $\gs$ in $L$ and $\{[\gt_0],\dots, [\gt_k]\}$ are the vertices of $\gs$ in $L'_\mu$, then
$\T_\mu(\gs)=\T_\mu(\gt_0) \cap \cdots \cap \T_\mu(\gt_0)$ is a disk of codimension $k$ in $D_\mu$. Hence, the $\T_\mu(\gt)$ intersect in the same fashion as the facets of the simplicial cone $[0,\infty)^{c(\mu)}$.
\end{proof}

\subsection{Gluing complexes of manifolds with boundary}\label{ss:systems}
A \emph{complex of manifolds with boundary over $\cs(L)$} is a collection of manifolds with boundary, $\{M_\gt\}_{\gt\in \cs(L)}$, together with embeddings $i_{\gs\gt}:M_\gt\hookrightarrow \partial M_\gs$ defined whenever $\gt<\gs$. In particular, we must have that $\partial M_\gs$ is nonempty whenever $\gs>\emptyset$; however, the minimum manifold $M_\emptyset$ can have empty boundary, and we usually assume this.  In other words, the manifolds are indexed by the vertices of $\cone(L')$ (the cone on the barycentric subdivision of $L$), while the embeddings $i_{\gs\gt}$ are indexed by the edges of $\cone(L')$. In addition, we require that there are certain $(k-1)$-parameter families of isotopies between the $i_{\gs\gt}$ which are indexed by the $k$-simplices of $\cone(L')$.  We eventually will want to require that the codimension of $M_\gt$ in $M_\gs$ is $c(\gt)-c(\gs)$ where $c:\cs(L)\to \nn$ is a function as in Subsection~\ref{ss:thick}, and that the image of $M_\gt$ in $\partial M_\gs$ has trivial normal bundle. 
Before giving the precise requirements let us mention that we also will be using the notion of a complex of manifolds over $\cs^0(L)$, where the empty simplex is not needed and where $\cone(L')$ is replaced by $L'$. 

 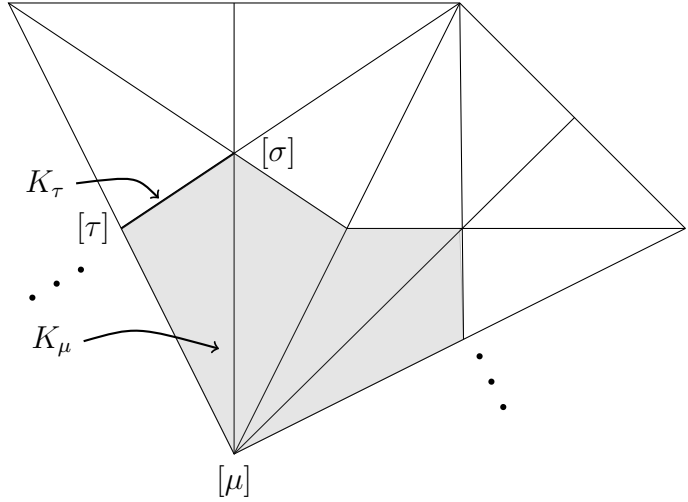
\begin{figure}
   \centering
   \begin{tikzpicture}[scale = 1]
   
\draw[thick] (-1.5,3) -- (0,4);
  
   \draw (3,6) -- (-3,6) -- (0,0) -- (3,6);
   \draw (0,0) -- (0,6);
   \draw (3,6) -- (-1.5,3);
   \draw (-3,6) -- (1.5,3);
      \draw (6,3) -- (1.5,3);
    \draw (0,0) -- (6,3);
        \draw (3,6) -- (6,3);
        \draw (0,0) -- (4.53,4.48);
       \draw (3,6) -- (3.05,1.525);

\draw[fill = gray, opacity = .2] (0,0) -- (-1.5, 3) -- (0,4) -- (1.5,3) -- (3,3) -- (3.05,1.525) -- (0,0);

    
    
    %

 \node[below] at (0,0) {$[\mu]$};
  \node[left] at (-1.5,3) {$[\gt]$};
  \node[right] at (.2,4) {$[\gs]$};
  \node[scale = 2,rotate=-65] at (3.45,.9) {$\dots$};
  \node[scale = 2,rotate=30] at (-2.3,2.3) {$\dots$};

    \draw[->][thick] (-2.1,3.6) .. controls (-1.3, 3.7) .. (-1,3.4);
  \node[left] at (-2.1,3.6) {$K_{\gt}$};
  
      \draw[->][thick] (-2,1.5) .. controls (-1.3, 1.7) .. (-.2,1.4);
  \node[left] at (-2,1.5) {$K_{\mu}$};
  

      \end{tikzpicture}
   \caption{Neighborhoods in the barycentric subdivision of $L$, the link of $\emptyset$ in $\cs(L)$.}
 \end{figure}

We want to describe how to glue together the $M_\gt \times D_\gt$, where $D_\gt$ is the dual disk defined in Subsection~\ref{ss:thick}.  It is easier to first describe how to glue together the $M_\gt\times K_\gt$ where $K_\gt$ is the dual cone of $\gt$ (i.e., $K_\gt =[\gt]*L'_\gt$ is a subcomplex of the $\vert\cs(L)\vert=[\emptyset]*L'$).  The gluing is accomplished using various embeddings:
	\[
	h_{\ga*\mu}:M_\mu \times \ga \hookrightarrow \partial M_{\min\ga} \times \ga,
	\]
where $\ga$ is a simplex in $L'$ and $\mu\in \cs(L)$ is a proper face of $\min \ga$.  The embedding $h_{\ga*\mu}$ will be used to glue $M_\mu\times (\ga*[\mu])$ onto $M_{\min \ga}\times K_{\min\ga}$.  For the gluing to be well-defined the $h_{\ga*\mu}$ must satisfy certain compatiblity relations which we will now describe. 

First suppose $\ga$ is a vertex of $L'$, i.e., $\ga=[\gt]$, where $\gt$ is a simplex of $L$.  Then $h_{[\gt]*[\mu]}:= i_{\gt\mu}\times \mathbb{I}_{[\gt]}: M_\mu\times [\gt]\to \partial M_\gt\times [\gt]$.  Essentially, $h_{[\gt]*[\mu]}$ is  $i_{\gt\mu}$.   Next, suppose $\ga=[\gt,\gs]$ is an edge in $L'$ so that $h_{[\gt,\gs]*[\mu]}: M_\mu\times [\gt,\gs] \to M_\gt\times [\gt,\gs]$ is an isotopy.  The restriction of this isotopy to $M_\mu\times [\gt]$ is $h_{[\gt]*[\mu]}:M_\mu\times  [\gt] \to M_\gt\times [\gt]$.  Its restriction to the other end when precomposed with $h_{[\gs]*[\gt]}$ should equal $h_{[\gs]*[\mu]}$, that is,
	\begin{equation}\label{e:iso1}
	(h_{[\gs]*[\gt]}) (h_{[\gt,\gs]*[\mu]}\vert_{M_\mu\times[\gs]}) = h_{[\gs]*[\mu]}.
	\end{equation}
In other words, the composition of the two gluing maps defined on the left hand side of \eqref{e:iso1} is equal to the gluing map on the right hand side.

Given a $k$-simplex $\ga=[\gt_0,\dots ,\gt_k]$ in $L'$, let $\ga_i=[\gt_0,\dots,\hat{\gt_i}, \dots, \gt_k]$ denote the face opposite $[\gt_i]$.  Then $h_{\ga*[\mu]}:M_\mu\times \ga\to \partial M_{\gt_0}\times \ga$ is a $k$-parameter isotopy. For $i\neq 0$, we require:
	\begin{equation}\label{e:iso2}
	h_{\ga*[\mu]}\vert_{M_\mu\times \ga_i}= h_{\ga_i*[\mu]},
	\end{equation}
while for $i=0$,
	\begin{equation}\label{e:iso3}
	(h_{\ga_0*[\gt_0]})(h_{\ga*[\mu]}\vert_{M_\mu\times \ga_0})=h_{\ga_0*[\mu]}
	\end{equation}

There is one further condition which  our isotopies should satisfy.  If $\gt=\min\ga$ and if $\mu$ and $\mu'$ are two faces of $\gt$ and $\gr=\mu\cap \mu'$, then we require that
	\begin{equation}\label{e:inter}
	\Ima (h_{\ga*[\mu]})\, \cap\, \Ima (h_{\ga*[\mu']}) = \Ima (h_{\ga*[\gr]})
	\end{equation}
In other words, on the complement of $\Ima (h_{\ga*[\gr]})$,  the embeddings $h_{\ga*[\mu]}:M_\mu\to M_\gt$ and $h_{\ga*[\mu']}:M_{\mu'}\to M_\gt$
must have disjoint images.

Next we describe how to glue together  $\{M_\gt \times K_\gt\}_{\gt\in \cs(L)}$ to obtain a space $X$ together with a projection map $p:X\to K$. Start with the disjoint union $\coprod M_\gt\times K_\gt$.  We will construct $X(0)\subset\cdots \subset X(k) \subset \cdots\subset X(d+1)=X$ so that $X(k)$ will be the inverse image of the $k$-skeleton of $K$ in $X$. The space $X(0)$ is defined to be the disjoint union $\coprod M_\gs \times [\gs]$.   Next, given an edge $[\gt,\gs]$ in $[\emptyset]*L'$, we glue $M_\gt \times [\gt,\gs]$ to $M_\gs\times [\gs]$ via $h_{[\gs]*[\gt]}:M_\gt\times [\gs]\to \partial M_\gs\times [\gs]$. (Recall that $h_{[\gs]*[\gt]}=i_{\gs\gt}$.)  After doing this gluing for each edge $[\gt,\gs]$, we obtain $X(1)$.  Notice that if $\gt$ is a $(d-1)$-simplex of $L$, then the link, $L'_\gt$, is the disjoint union of the vertices $[\gs]$ where $\gt<\gs$ and $K_\gt=[\gt]*\coprod\, [\gs] =\bigcup\, [\gt,\gs]$; hence, after building $X(1)$ we will have glued $M_\gt\times K_\gt$ onto $X(0)$.  Next, consider a $2$-simplex $[\mu,\gt,\gs]$.  Glue $M_\mu\times [\mu,\gt,\gs]$ onto $M_\gt\times [\gt,\gs]$ using $h_{[\gt,\gs]*[\mu]}:M_\mu \times [\gt,\gs]\to \partial M_\gt \times [\gt,\gs]$.  By \eqref{e:iso1} this is compatible with the previously defined gluing map $M_\mu\times [\gs] \to \partial M_\gs\times [\gs]$.  Hence, the union of the two gluing maps gives a well-defined map $M_\mu \times [\gt,\gs]\to (\partial M_\gt \times [\gt,\gs]) \cup (\partial M_\gs\times [\gs])$ which we can use to glue $M_\mu\times [\mu,\gt,\gs]$ onto $X(1)$.

\begin{figure}\label{f:gluing}
\centering
\begin{tikzpicture}[scale = .65]

\draw (-.5,1) -- (-.5,-1);
 \draw [fill] (-.5,1) circle [radius=0.04];
  \draw [fill] (-.5,-1) circle [radius=0.04];
\draw [fill] (0.8,-2) circle [radius=0.05];
   
 \draw[ fill =  gray, opacity = .2] (3,-1) -- (3,1) -- (1,2) -- (1,0) -- (3,-1);
  \draw (3,-1) -- (3,1) -- (1,2) -- (1,0) -- (3,-1);
\draw(3, -1) -- (4, -.5);
\draw(3, 1) -- (4, 1.5);
\draw(1, 2) -- (2, 2.5);
   \node[left] at (-.5,0) {$M_\gt$};
   \node[below] at (0.8,-2) {$M_\mu$};
    \node at (3,1.75) {$M_\gs$};
 \node[below] at (0.8,-3) {$X(0)$};
  \node[below] at (-.5,-1) {$\partial M_\gt$};
  
    \node at (2,.5) {$\partial M_\gs$};
 
\begin{scope}[xshift = 12cm]
\draw (-.5,1) -- (-.5,-1);
 \draw [fill] (-.5,1) circle [radius=0.04];
  \draw [fill] (-.5,-1) circle [radius=0.04];
\draw [fill] (0.8,-2) circle [radius=0.05];
   
 \draw[ fill =  gray, opacity = .2] (3,-1) -- (3,1) -- (1,2) -- (1,0) -- (3,-1);
  \draw (3,-1) -- (3,1) -- (1,2) -- (1,0) -- (3,-1);
\draw(3, -1) -- (4, -.5);
\draw(3, 1) -- (4, 1.5);
\draw(1, 2) -- (2, 2.5);

\draw[fill = gray, opacity = .2](-.5, 1) -- (2, 1.5) -- (2,-.5) -- (-.5, -1);

 \draw [fill] (2,1.5) circle [radius=0.04];
  \draw [fill] (2,-.5) circle [radius=0.04];
   \draw [fill] (-.5,-1) circle [radius=0.04];
  \draw [fill] (2,-.5) circle [radius=0.04];

\draw(-.5, 1) -- (2, 1.5) -- (2,-.5) -- (-.5, -1);
\draw(.8, -2) -- (-.5, -1);
\draw(.8, -2) -- (2, -.5);

   \node[left] at (-.5,0) {$M_\gt$};
   \node[below] at (0.8,-2) {$M_\mu$};
    \node at (3,1.75) {$M_\gs$};
     \node[below] at (0.8,-3) {$X(1)$};
     
      \draw[->][thick] (-1.6,-1) .. controls (-1,-1.1) .. (-.6,-1);
       \node[left, scale = .9] at (-1.6,-1) {$h_{[\gt]\ast[\mu]}(M_\mu \times [\tau])$};
         \draw[->][thick] (3,-1.5) .. controls (2.5,-1.2) .. (2,-.6);
          \node[right, scale = .9] at (3,-1.5) {$h_{[\gs]\ast[\mu]}(M_\mu \times [\gs])$};
          
                   \draw[->][thick] (0,2) .. controls (1.5,1) .. (2,.5);
          \node[left, scale = .9] at (0,2) {$h_{[\gs]\ast[\gt]}(M_\gt \times [\gs])$};
     
    \end{scope}
    
    \begin{scope}[yshift = -6cm, xshift = 6cm]
\draw (-.5,1) -- (-.5,-1);
 \draw [fill] (-.5,1) circle [radius=0.04];
  \draw [fill] (-.5,-1) circle [radius=0.04];
\draw [fill] (0.8,-2) circle [radius=0.05];
   
    \draw [fill] (2,1.5) circle [radius=0.04];
  \draw [fill] (2,-.5) circle [radius=0.04];
   \draw [fill] (-.5,-1) circle [radius=0.04];
  \draw [fill] (2,-.5) circle [radius=0.04];
   
 \draw[ fill =  gray, opacity = .25] (3,-1) -- (3,1) -- (1,2) -- (1,0) -- (3,-1);
  \draw (3,-1) -- (3,1) -- (1,2) -- (1,0) -- (3,-1);
\draw(3, -1) -- (4, -.5);
\draw(3, 1) -- (4, 1.5);
\draw(1, 2) -- (2, 2.5);

\draw[fill = gray, opacity = .25](-.5, 1) -- (2, 1.5) -- (2,-.5) -- (-.5, -1);
\draw(-.5, 1) -- (2, 1.5) -- (2,-.5) -- (-.5, -1);
\draw(.8, -2) -- (-.5, -1);
\draw(.8, -2) -- (2, -.5);

   \node[left] at (-.5,0) {$M_\gt$};
   \node[below] at (0.8,-2) {$M_\mu$};
    \node at (3,1.75) {$M_\gs$};
    \draw[fill = gray, opacity = .25](.8, -2) -- (2, -.5) -- (-.5, -1);
     \node[below] at (0.8,-3) {$X(2)$};
     
             \draw[->][thick] (3,-1.5) .. controls (2.5,-1.4) .. (.7,-.85);
          \node[right, scale = .9] at (3,-1.5) {$h_{[\gt,\gs]\ast[\mu]}(M_\mu \times [\gt,\gs])$};
    \end{scope}

\end{tikzpicture}
\caption{The first stages of our gluing procedure before thickening}
\end{figure}
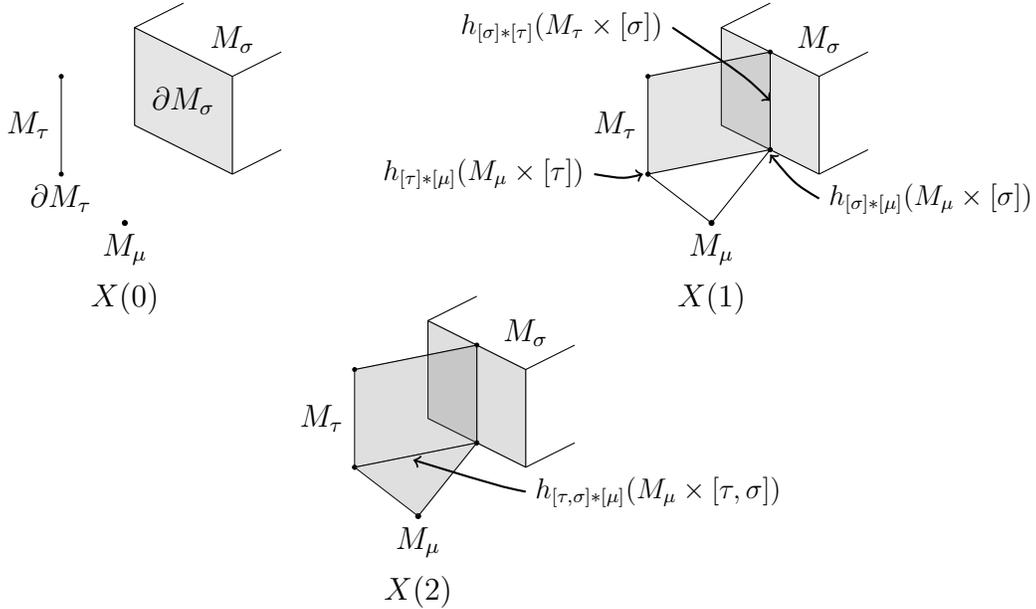

Continue by induction.  Suppose $X(k)$ has been defined over the $k$-skeleton of $K$.  For each $k$-simplex $\ga$ in $L'$ and $\mu\in \cs(L)$ with $\mu < \min \ga$, we have $h_{\ga*\mu}:M_\mu \times \ga \to \partial M_{\min\ga}\times \ga$.  By \eqref{e:iso2} and \eqref{e:iso3}, this is compatible with previously defined gluing maps $h_{\gb*[\mu]}$, where $\gb$ is a face of $\ga$.  So, for  $\ga=[\gt_0,\dots, \gt_k]$, we get a well-defined map from $M_\mu\times \ga$ onto the image of 
\[
(\partial M_{\gt_0}\times [\gt_0,\dots, \gt_k]) \cup (\partial M_{\gt_1}\times [\gt_1,\dots, \gt_k] )\cup \cdots  \cup (\partial M_{\gt_k}\times [\gt_k] )
\]
in $X(k)$.

We can summarize the above as follows.  The union of the  $h_{\ga*[\mu]}$, with $\min \ga=\gt$, define an embedding:
	\begin{equation}\label{e:H}
	H_{\gt\mu}:M_\mu\times \partial K_\mu(\gt)\to \partial M_\gt \times K_\gt.
	\end{equation}
(Recall $\partial K_\mu(\gt)= K_\gt$.)  Formulas \eqref{e:iso1}, \eqref{e:iso2} and \eqref{e:iso3} imply that whenever $\mu<\gt<\gs$:
	\begin{equation}\label{e:iso4}
	(H_{\gs\gt} \vert_{\partial M_\gt\times \partial K_\gt(\gs)}) (H_{\gt\mu}\vert_{\partial K_\mu(\gs)})=H_{\gs\mu} .
	\end{equation}
Hence, if $X_{>\mu}=p^{-1}(\partial K_\mu)$, then  the union of the $H_{\gt\mu}$ fit together to give a well-defined map $H_\mu:M_\mu\times \partial K_\mu\to X_{>\mu}$,  that specifies the gluing of $M_\mu\times \partial K_\mu$ onto $X_{>\mu}$.  By \eqref{e:inter} $H_\mu$ is an embedding.

\begin{remark}
If each $M_\gt$ is aspherical, then $X$ is a model for the aspherical realization $BG\cs(L)$.
\end{remark}

Next we want to replace dual cones by dual disks.  We suppose that
\begin{enumeratea}
\item
$M_\emptyset$ is a point.
\item 
For all maximal simplices $\gs$ in $L$, the manifolds $M_\gs$ all have the same dimension $n$.
\item \label{i:codim2}
For each $\gt<\gs$, $\dim M_\gs -\dim M_\gt \ge 2\codim (\gt,\gs)$. 

\end{enumeratea}

The conditions in \eqref{i:codim2} are called the \emph{codimension $\geq 2$ conditions}.
The dimensions of the $M_\gt$ give us the data for a function $c:\cs(L)\to \nn$ as in Subsection \ref{ss:thick}, defined by
	\[
	c(\gt) = n- \dim M_\gt.
	\]
As in Subsection \ref{ss:thick}, we can use the function $c$ to define the thick link $\T_\mu$ for each $\mu\in\cs(L)$.  

Whenever $\mu$ is a codimension-one face of $\gt$,  $\T_\mu(\gt)$ is a regular neighborhood of $[\gt]$ in $S'_\mu$.  

Moreover, $\T_\mu(\gt)$ is homeomorphic to $E\times D_\gt$ where $E$ is a disk of  dimension equal to the codimension of $M_\mu$ in $\partial M_\gt$.  Since the normal bundle of $M_\mu$ in $\partial M_\gt$ is trivial,  $M_\mu\times E$ can be embedded as a codimension-zero submanifold of $\partial M_\gt$ and
hence, the embedding $H_{\gt\mu}$ defined by \eqref{e:H} extends to an embedding:
	\begin{equation}\label{e:J}
	J_{\gt\mu}:M_\mu\times \T_\mu(\gt)\hookrightarrow \partial M_\gt \times D_\gt.
	\end{equation}
When $\codim (\mu,\gt) >1$, $\T_\mu(\gt)$ remains a thickening of $\partial K_\mu(\gt)$ (as well as a thickening of $D_\gt$) and the embedding $H_{\gt\mu}$ again extends to an embedding $J_{\gt\mu}:M_\mu\times \T_\mu(\gt)\to \partial M_\gt \times D_\gt$.  Moreover, the $J_{\gt\mu}$ satisfy the analogous formulas to \eqref{e:iso4}.

Finally, we build an $n$-manifold with boundary $M$ in exactly the same fashion we constucted the space $X$ above, namely,
	\begin{equation}\label{e:M}
	M:=\left(\coprod _{\gt\in \cs(L)} M_\gt \times D_\gt \right)\, / \sim
	\end{equation}
where the equivalence relation $\sim$ is defined as before except that we use as gluing maps the $J_{\gt\mu}$ rather than the  $H_{\gt\mu}$.

At this point we should explain why $M$ is a $n$-manifold with boundary.  Since $D_\mu$ is a $c(\mu)$-manifold with corners (cf.\,Lemma~\ref{l:corners}) and $M_\mu$ is a $(n-c(\mu))$-manifold with boundary, $M_\mu\times D_\mu$ is a $n$-manifold with corners.  Each facet is either the product of $M_\mu$ with a facet of $D_\mu$ or it has the form $\partial M_\mu\times D_\mu$. If $\mu<\gt$ is a codimension-one face, then $M_\mu\times D_\mu$ is glued onto $M_\gt\times D_\gt$ by the embedding $J_{\gt\mu}:M_\mu\times \T_\mu(\gt)\to\partial M_\gt \times D_\gt$ from a facet of $M_\mu\times D_\mu$ to a facet of $M_\gt\times D_\gt$.  Similarly, if $\{[\gt_0], \dots, [\gt_k]\}$ is the vertex set of $\gs$ in $L'_\mu$, then $M_\mu\times D_\mu$ is glued onto $M_\gs\times D_\gs$ via a map from the stratum $M_\mu\times T_\mu(\gs)$ to an intersection of strata in $\partial M_\gs\to D_\gs$.

Let $\{M_\gt\}_{\gt\in \cs(L)}$ be a complex of manifolds with boundary over $\cs(L)$. Let $G_\gt = \pi_1(M_\gt)$.  As usual suppose the induced homomorphisms $\phi_{\gs\gt}=(i_{\gs\gt})^*:G_\gt\to G_\gs$ are injective and the system $G\cs(L)=\{G_\gt\}_{\gt\in\cs(L)}$ is a simple complex of groups. Let $G=\gsl$ be the direct limit. Also suppose
	\begin{itemize}\label{e:conditions}
	\item
	$G\cs(L)$ is developable, and
	\item
	each $M_\gt$ is homotopy equivalent to $BG_\gt$,
	\end{itemize}
 
\begin{Theorem}\label{t:gluing}
Suppose $\{M_\gt\}_{\gt\in \cs(L)}$ is a complex of manifolds with boundary over $\cs(L)$ satisfying the above conditions and let $M$ be the manifold with boundary defined by \eqref{e:M}.  Then $M$ is a thickening of $BG\cs(L)$.  If the $K(\pi,1)$-Question for $G\cs(L)$ has a positive answer, then $M\sim BG$, so
	\[
	\actdim G\le \dim M.
	\]
\end{Theorem}

\subsection{Complexes of manifolds with boundary over $\cs^o(L)$}\label{ss:systems0}

We turn to the case where the manifolds $M_v$ associated to the vertices $v$ of $L$ are allowed to have empty boundary. 
When this happens our algorithm does not tell us how to glue on the final disk $M_\emptyset \times D_\emptyset$. (If $\partial M_v=\emptyset$, there is no place to embed $M_\emptyset \times \T_\emptyset$.) One can still try to accomplish the gluings without  the final disk.  For each higher-dimensional simplex $\gt$, we still require $\partial M_\gt$  to be nonempty. By excluding the empty simplex, we consider a complex of manifolds with boundary  
over $\cs^o(L)$ where the embeddings $i_{\gs\gt}$ and gluing maps $H_{\gt\mu}$, $J_{\gt\mu}$ satisfy the same conditions as in the previous subsection.  We can glue together the pieces together as before; however, the resulting manifold will usually not be a model for $BG$.

So, let us consider a complex of manifolds with boundary $\{M^o_\gt\}_{\gt\in \cs^o(L)}$, where $\partial M^o_v$ is allowed to be empty for $v\in \vertex L$, but $\partial M^o_\gt\neq \emptyset$ if $\dim \gt>0$.  The other conditions in the previous section are satisfied mutatis mutandis.
The complex is \emph{$n$-dimensional} if $\dim M^o_\gs=n$ for each maximal simplex $\gs$.  As in \eqref{e:M}, we can glue together the $\{M^o_\gt\times D_\gt\}_{\gt\in \cs^o(L)}$ to obtain an $n$-manifold with boundary $M^o$.  As before, $M^o$ will a model for $BG\cs^o(L)$; however, as explained in Subsection~\ref{ss:cq=cs},
$BG\cs^o(L)$ will not be homotopy equivalent to $BG$ unless $L$ is contractible (see Lemma~\ref{l:BGS}).  When $L$ is EDCE (see Definition~\ref{d:EDCE}), we can use the methods of Subsection~\ref{ss:cq=cs} to conclude that $\actdim G\le n$.

Given the $n$-dimensional system $\{M^o_\gt\}_{\gt\in \cs^o(L)}$ over $\cs^o(L)$, there is an easy way to extend it to an $(n+1)$-dimensional system over $\cs(L)$: simply take the product of each manifold $M^o_\gt$ with the unit interval. In other words, put $M_\gt= M^o_\gt \times [0,1]$. So, the vertex manifold $M^o_v$ is replaced by the manifold with boundary $M_v=M^o_v\times [0,1]$. Let $M_\emptyset$ be a point and $i_{\gt\emptyset}:M_\emptyset \to \partial M_\gt$ an appropriate embedding.

Since $\pi_1(M^o_\gt)=\pi_1(M_\gt)$, the systems of fundamental groups $\{\pi_1(M_\gt)\}$ and $\{\pi_1(M^o_\gt)\}$ define the same simple complex of groups $G\cs(L)$ over $\cs(L)$ given by $G_\gt=\pi_1(M_\gt)$.  Let $G=\lim G\cs(L)$.
Let $M^o$ and $M$ denote, respectively, the results of gluing together the systems $\{M^o_\gt\}_{\gt\in \cs^o(L)}$ and $\{M_\gt\}_{\gt\in \cs(L)}$. 
We assume each $M_\gt$ is connected, so that $\pi_1(M) = G$.
If $L$ is simply connected, then $\pi_1(M^o)$ is also equal to $G$.  (In general, $\pi_1(M^o)$ is the semidirect product described in Remark~\ref{r:fund}.)

From Lemma~\ref{l:BGS} we get the following.
\begin{lemma}\label{l:coneL}
The inclusion $M^o\hookrightarrow M$ is a homotopy equivalence if and only if $L$ is contractible.
\end{lemma}

Suppose $L$ is a full subcomplex of another simplicial complex $\bar{L}$. In Subsection~\ref{ss:cq=cs} we defined the notion of a trivial extension of a simple complex of groups over $\cs(L)$ to one over $\cs(\bar{L})$.  Similarly, we can define the notion of a trivial extension of a complex of manifolds with boundary over $\cs(L)$.  Let $\{M_\gt\}_{\gt\in \cs(L)}$ be such a system. For simplicity, suppose that the dimension of $M_\gt$ depends only on the dimension of the simplex $\gt$, i.e., if $\dim M_\gt=n(k)$ for each $k$-simplex $\gt$. Recall any simplex $\bar{\gs}$ can be decomposed as a join $\bar{\gs}=\gs*\gt$, where $\gs\in \cs(L)$ and where vertices of $\gt$ lie in $\bar{L}-L$.  
A complex of manifolds with boundary $\{N_\gt\}_{\gt\in \cs(\bar{L})}$ is called a \emph{trivial extension} of $\{M_\gt\}_{\gt\in \cs(L)}$ if 
$N_{\bar{\gs}}=M_{\gs}\times D^{n(k)}$, whenever $\bar{\gs}=\gs*\gt$ as above and $D^{n(k)}$ is the disk $N_\gt$.
 A \emph{trivial extension} of a complex over $\cs^o(L)$ is defined similarly.

The next lemma follows immediately from Lemma~\ref{l:trivialext}.
\begin{lemma}\label{l:extension}
Suppose $\{M_\gs\}_{\gs\in \cs(L)}$ is a complex of manifolds and $\{N_{\bar{\gs}}\}_{\bar{\gs}\in \cs(\bar{L})}$ is a trivial extension of it over $\bar{L}$. Let $M$ and $N$ be the result of gluing these systems together. Then $M$ is homotopy equivalent to $N$.
\end{lemma}

\begin{Prop}\label{p:ext}
Suppose $L$ is a full subcomplex of a contractible complex $L_c$ of the same dimension as $L$. Let $\{M^o_\gt\}_{\gt\in \cs^o(L)}$ be an $n$-dimensional complex of manifolds with boundary over $\cs^o(L)$ and $\{N^o_\gt\}_{\gt\in \cs^o(L_c)}$ a trivial extension of it to $L_c$.  Let $N^o$ be the result of gluing together $\{N^o_\gt\}$ and let $M$ and $N$ be the result of gluing the corresponding $(n+1)$-dimensional complexes over $\cs(L)$ and $\cs(L_c)$, respectively. Then $M$ is homotopy equivalent to $N^o$ 
\end{Prop}

\begin{proof}
By Lemma~\ref{l:extension}, $M$ is homotopy equivalent to $N$ and by Lemma~\ref{l:coneL}, $N$ is homotopy equivalent to $N^o$.
\end{proof}

For the next theorem we suppose that $\{M^o_\gs\}$ is a $n$-dimensional complex of manifolds with boundary over $\cs^o(L)$ where the vertex manifolds $M^o_v$ are allowed to have empty boundary.  Further suppose that the conditions for Theorem~\ref{t:gluing} are satisfied, i.e.,  each $M^o_\gs$ is aspherical, the associated simple complex of groups $G\cs(L)$ is developable and  $BG\cs(L)\sim BG$ for $G=\gsl$  (so the $K(\pi,1)$-Question has a  positive answer).

\begin{Theorem}\label{t:gluing2}
Let $\{M^o_\gs\}_{\gs\in \cs^o(L)}$ be an $n$-dimensional complex of aspherical manifolds with boundary.  Then, with hypotheses as above:
\begin{enumeratei}
\item
$\actdim G\le n+1$.
\item
If $L$ is EDCE, then $\actdim G\le n$.
\end{enumeratei}
\end{Theorem}

\begin{proof}
Statement (i) follows from Theorem~\ref{t:gluing} applied to the $(n+1)$-dimensional system $\{M_\gs\}$, where $M_\gs=M^o_\gs\times [0,1]$.  When $L$ is EDCE, $L$ embeds in a contractible complex $L_c$ of the same dimension; so, by Proposition~\ref{p:ext}, $M$ is homotopy equivalent to the $n$-manifold $N^o$.  Hence, $\actdim G\le n$.
\end{proof}

\begin{remark}
Most  of the constructions in Subsections \ref{ss:systems} and \ref{ss:systems0} work if one replaces the poset of simplices $\cs(L)$ by a more arbitrary poset $\cq$.  Specifically, suppose $\cq$ is a poset with a minimum element $m$ and that $\cq^o=\cq-\{m\}$.  As before, one can define the notion of a posets of manifolds with boundary over $\cq$ and $\cq^o$, respectively,  and  prove versions of Theorems \ref{t:gluing} and \ref{t:gluing2}.
\end{remark}

\section{Examples}\label{s:examples} 
\subsection{Simple complexes of closed aspherical manifolds}\label{ss:gprdts}
We begin by discussing graph products. 
Notation is continued from Definition~\ref{d:graphprdt} and Example~\ref{ex:gprdt}: $L^1$ is a simplicial graph with vertex set $V$, $L$ is the associated flag complex, $\{G_v\}_{v\in V}$ is a collection of groups of type $F$, $G=\prod_{L^1} G_v$ denotes their graph product and, for each $\gs\in \cs(L)$, $G_\gs=
\prod_{v\in \vertex\gs} G_v$ is the direct product. As in Example~\ref{ex:gprdt}, this defines a simple complex of groups, $G\cs(L)=\{G_\gs\}_{\gs\in \cs(L)}$. By Corollary~\ref{c:gprdt}, the $K(\pi,1)$-Question has a positive answer for $G\cs(L)$. 
 For each $v\in V$, let $M_v$ be a model for $BG_v$ by a manifold of minimum dimension $m_v=\actdim G_v$. Let $M_\gs=\prod_{v\in \vertex \gs} M_v$.  
 
Suppose each $M_v$ has nonempty boundary and has dimension at least $2$.  Let $m_\gs=\dim M_\gs=\sum_{v\in \vertex \gs}\dim M_v$ and put 
	\begin{equation}\label{e:dimn}
	n=\sup\{m_\gs\mid \gs\in \cs(L)\}
	\end{equation}
By taking products with disks of suitable dimensions, we can assume that for each maximal simplex $\gs$ of $L$, $\dim M_\gs = n$.  It is   straightforward to give $\{M_\gs\}_{\gs\in \cs(L)}$ the structure of an $n$-dimensional complex of manifolds with boundary satisfying the conditions in Subsection~\ref{ss:systems}.  Here  are the details.  Choose basepoints $x_v\in \partial M_v$.  Whenever  $\mu$ is a face of a simplex $\gt\in \cs(L)$,
put $V(\gt\mu)=\vertex \gt-\vertex \mu$ and $x_{\gt\mu}=(x_v)_{v\in V(\gt\mu)}$ be a basepoint in  $\prod _{v\in V(\gt\mu)} \partial M_v$.  This gives an inclusion
	\[
	i_{\gt\mu}:M_\mu\hookrightarrow M_\mu\times x_{\gt\mu}\subset\partial M_\gt
	\]
with trivial normal bundle.  For a $k$-simplex $\ga = [\gt_0, \gt_1, \dots, \gt_k] \in L'$, define 
	\[
	h_{[\ga \ast \gt_0]}: M_\mu \times \alpha \rightarrow  \partial M_{\gt_0} \times \alpha
	\]
by $h_{[\ga \ast \gt_0]}(z,x) = (i_{\gt\mu}(z),x)$.
The $h_{[\ga \ast \gt_0]}(z,x)$ obviously satisfy the conditions in Subsection \ref{ss:systems}; so, $\{M'_\gs\}_{\gs\in \cs(L)}$ is a complex of manifolds with boundary. 
Using Theorem~\ref{t:gluing}, we get an $n$-dimensional manifold $M$ which is a model for $BG$.  This gives the following.

 \begin{Prop}\label{p:gprdt1}\textup{(cf.\ Theorem~\ref{t:gluing} and Corollary~\ref{c:gprdt}).}
 Suppose each $G_v$ is the fundamental group of an aspherical manifold $M_v$ with nonempty boundary.  Let
$G=\prod_{L^1} G_v$ be the graph product. Then
	 \[
	 \actdim G\le \sup \{m_\gs\mid \gs\in \cs(L)\},
	 \]
 where $m_\gs=\dim M_\gs$.
 \end{Prop} 

We turn to the case where each $M_v$ is a closed aspherical manifold not equal to a point. (For example, when $G$ is a $\raag$ each $M_v$ is a circle.) We can convert this case into the first case by the simple expedient of taking the product of each $M_v$ with $[0,1]$.
Then $M'_v=M_v\times [0,1]$ is a manifold with boundary of dimension $m_v+1$.   Proposition~\ref{p:gprdt1} then has the following corollary.

\begin{corollary}\label{c:closedgprdt}
Suppose each $G_v$ is the fundamental group of a closed aspherical manifold $M_v$. Let
$G=\prod_{L^1} G_v$ be the graph product. Then
	 \[
	 \actdim G\le \sup \{m'_\gs\mid \gs\in \cs(L)\},
	 \]
 where $M'_\gs=M_\gs\times [0,1]^{\vertex \gs}$ and $m'_\gs=\dim M'_\gs=\dim M_\gs +\dim\gs +1$.
 \end{corollary} 
 
The argument in Proposition~\ref{p:gprdt1} and Corollary~\ref{c:closedgprdt} is easily modified to give a sharp upper bound for $\actdim G$ in the case when some $M_v$ have empty boundary and some do not.
 
When each $M_v$ is  closed and of  the same dimension $m$, then  Corollary~\ref{c:closedgprdt} has the following corollary.
 
 \begin{corollary}\label{c:closedgprdt2}
 Suppose each $G_v$ is the fundamental group of a closed aspherical $m$-manifold $M_v$, that the flag complex $L$ has dimension $d$, and that $G=\prod_{L^1} G_v$ is the graph product.Then
 	\[
 	\actdim G\le (m+1)(d+1)
	 \]
 \end{corollary}
 
 \begin{proof}
 If $\gs$ is a $d$-simplex, then $\dim M'_\gs=(m+1)(d+1)$.
 \end{proof}

The arguments above for the graph product complex of fundamental groups of closed aspherical $m$-manifolds can be generalized to a simple complex of fundamental groups of closed aspherical manifolds as defined below.
 
\begin{definition}\label{d:simpleclosed}
A \emph{simple complex of closed manifolds} over $\cs(L)$ is a collection of connected, closed manifolds $\{M_\gs\}_{\gs\in \cs(L)}$ and for each $\gt<\gs$ a $\pi_1$-injective embedding $i_{\gs\gt}:M_\gt \hookrightarrow M_\gs$ as a submanifold. There are a few more requirements:
\begin{itemize}
\item
$M_\emptyset$ is a point.
\item
If $G_\gs=\pi_1(M_\gs)$ and $\phi_{\gs\gt}:G_\gt\to G_\gs$ is the homomorphism defined by $i_{\gs\gt}$, then 
the system $\{G_\gs, \phi_{\gs\gt}\}$ is a simple complex of groups $G\cs(L)$. 
\item
If $\gt_1,\dots , \gt_k$ are faces of $\gs$, then the $M_{\gt_i}$ intersect transversely in $M_\gs$.
\item
If $\gt$ is a face of $\gs$, then $M_\gt$ has trivial normal bundle in $M_\gs$. 
\end{itemize}
For $G_\gs=\pi_1(M_\gs)$, let  $G\cs(L)=\{G_\gs\}_{\gs\in \cs(L)}$ be the associated simple complex of groups.  If each $M_\gs$ is a model for $BG_\gs$, then $\{M_\gs\}_{\gs\in \cs(L)}$ is a \emph{simple complex of closed aspherical manifolds}.
\end{definition}

\begin{definition}\label{d:codim_m}
The system $\{M_\gs, i_{\gs\gt}\}$ satisfies \emph{codimension-$m$ conditions} if $M_\gt$ has codimension $m$ in $M_\gs$ whenever $\gt$ is a codimension one face of $\gs$.  (This implies that for any face $\gt$ of $\gs$, $\codim (M_\gt, M_\gs) = m\codim (\gt,\gs)$.)
\end{definition}

Note that since $\dim M_\emptyset =0$, the codimension-$m$ conditions imply that $\dim M_\gs=m(\dim\gs +1)$.

For the graph product complex,  if each vertex manifold $M_v$ is a closed aspherical $m$-manifold, then $\{M_\gs\}_{\gs\in \cs(L)}$ is a simple complex of closed aspherical manifolds satisfying the codimension-$m$ conditions. 
  
In the more general case of a complex of closed manifolds, follow the same procedure as for graph products and replace each manifold $M_\sigma$ by the manifold with boundary $M_\sigma' = M_\sigma \times [0,1]^{\vertex \gs}$. Then $\{M_\gs' \}_{\sigma \in \cs(L)}$ is a complex of manifolds with boundary as in Subsection \ref{ss:systems}.  Theorem \ref{t:gluing} allows us glue together the $\{M_\gs'\}$ to get a manifold $M'$. 

Suppose $G\cs(L)$ is developable, that the intersection of local groups condition \eqref{e:simplexinter} holds and   that $L$ is a flag complex (so that the $K(\pi,1)$-Question for $G\cs(L)$ has a positive answer). Then we get an upper bound for the action dimension as follows.

\begin{Theorem}\label{t:complexaspherical}
Suppose that $\{M_\gs\}_{\gs \in \cs(L)}$ is a simple complex of closed aspherical manifolds over $\cs(L)$ satisfying the codimension-$m$ conditions. Let $G = \lim G\cs(L)$ be the direct limit of the $\pi_1(M_\gs)$. Also suppose as above that $L$ is a $d$-dimensional flag complex and that the $K(\pi,1)$-Question has a positive answer for $G\cs(L)$. Then $$\actdim(G) \le (m+1)(d+1).$$
\end{Theorem}

In section \ref{s:obstructors}, we will prove that if $H_d(L, \zz_2) \ne 0$, then $\actdim(G) = (m+1)(d+1).$

\paragraph{Intersecting submanifolds.}
One way to get examples of simple complexes of closed manifolds is to start with an ambient manifold $M$ together with a collection of connected, smooth submanifolds (e.g., hypersurfaces), $\{M(i)\}_{i\in I}$, so that the $M(i)$ intersect transversely and so that for each subset $J$ of $I$, the intersection $M(J)=\bigcap_{i\in J}\,M(i)$ is nonempty.  In order to make the indexing compatible with previous notation we must replace subsets of $I$ by their complements.  So, we will write $v(i)$ for $I-\{i\}$ and $V$ for $\{v(i)\}_{i\in I}$, the set of complements of singletons. A subset $J\subset I$ corresponds to the complementary subset $\gs(J):=\{v(i)\}_{i\in I-J}$ of $V$.  For example, in the case of the graph product complex, all  manifolds are subspaces of the manifold $M=\prod_{v\in V}\,M_v$, while $M(i)=\prod_{j\in I-\{i\}}\,M_{v(j)}$.  The case of a $\raag$ is a further specialization:  $M$ is the $T^I$, the torus on $I$, the $T(i)$ are coordinate subtori of codimension one, the $T_{v(i)}$ are coordinate circles and the $T_\gs$ are coordinate subtori.

To simplify the discussion,  suppose that $I=\{1,\dots, p\}$, that $\dim M= mp$ and that each $M(i)$ has the same codimension $m$ in $M$.  Since the intersections are transverse and nonempty, the intersection $M(I)$ of all the $M(i)$ is a nonempty finite set of points.  Choose one, $x_0$, as the basepoint.  For each $\gs=\gs(J)\subset V$, let $M_\gs$ be the component of $M(J)$ containing $x_0$.   Thus, $\dim M_\gs=m(\dim \gs +1)$.

If $M$ is locally $\cat(0)$ and each $M(i)$ is totally geodesic, then each $M_\gs$ is aspherical with fundamental group $G_\gs:=\pi_1(M_\gs)$. By uniqueness of geodesics, for each $\gt<\gs$ the inclusion $M_\gt\to M_\gs$ is $\pi_1$-injective.   For the same reason, the intersection of local subgroups condition \eqref{e:simplexinter} holds.  Thus, for any simplicial complex $L$, $G\cs(L)$ is a simple complex of groups and when $L$ is a flag complex, the $K(\pi,1)$-Question has a positive answer.

\begin{example}\label{ex:sym} (\emph{Complexes of hyperbolic manifolds}).
First we use the above technique to construct a system of hyperbolic manifolds satisfying the codimension-$1$ conditions.
Suppose $K=\QQ(\sqrt{d})$ is a totally real quadratic extension of $\QQ$ and $A$ is the ring of algebraic integers in $K$. Choose $\geps\in A$ so that $\geps >0$ and $\ol{\geps}<0$. Define a quadratic form $\gf:A^{p+1}\times A^{p+1} \to A$ by
	\[
	\gf(e_i,e_j)=
	\begin{cases}
	\gd_{ij},	&\text{if $(i,j)\neq (0,0$),}\\
	-\geps,	&\text{if $(i,j)= (0,0)$,}
	\end{cases}
	\]
where $\{e_i\}_{0\le i\le p}$ is the standard basis. Let $O(\gf)\subset GL_{p+1}(A)$ be the subgroup which preserves the quadratic form $\gf$. The signature of $\gf$ on $A^{p+1}\otimes _A \rr \cong \rr^{p+1}$ is $(p,1)$; so, over $\rr$, the group of isometries of $\gf$ is identified with $O(p,1)$ and $O(\gf)$ is a uniform lattice in $O(p,1)$. Let $\gG$ be a normal torsion-free subgroup of $O(\gf)$ (for example, for almost any prime ideal, we could take $\gG$ to be the corresponding congruence subgroup of $O(\gf)$). Then $M^p=\hh^p/\gG$ is a closed hyperbolic manifold. The image of $(\sqrt {\geps}, 0, \dots, 0)$ in $M^p$ is the basepoint. For $1\le i\le p$, let $r_i\in O(\gf)$ be the reflection which sends $e_i$ to $-e_i$. There is an induced involution of $\ol{r}_i$ on $M^p$. Its fixed point set is a totally geodesic submanifold of codimension one and the component containing the basepoint is the manifold $M(i)$.  If we require $\gG$ to be a subgroup of $SO(\gf)$, then $M(i)$ will be orientable and have trivial normal bundle in $M$.  For any flag complex $L$ with vertex set in $V$ we then get a simple complex of closed hyperbolic manifolds $\{M_\gs\}_{\cs(L)}$.  This system behaves  similarly to a complex of coordinate tori in $T_p$.  It is easy to modify this construction to get a subsystem of $\{M_\gs\}$ satisfying the codimension-$m$ conditions.  Let $p=\ol{p}m$.  Group the first $m$ commuting involutions together to get an involution $s_1=r_1r_2\dots r_m$.  Continue in this fashion, defining $s_j=r_{m(j-1)+1}\cdots  r_{mj}$.  The fixed set of $\ol{s}_j$ on $M$ is a hyperbolic submanifold $M(j)$ of codimension-$m$.  If $L$ is a flag complex with $\ol{p}$ vertices we can use the  intersections of the $M(j)$ as above to get a system of hyperbolic manifolds satisfying codimension-$m$ conditions.

We have that $\bigcup M_\gs$ is a model for $BG\cs(L)$.  The union is a piecewise hyperbolic space.  Since $L$ is a flag complex, its all right, piecewise spherical metric is $\cat(1)$.  Since all intersections in $M$ were orthogonal, the relevant links of the $M_\gs$ in the union are full subcomplexes of $L$.  It follows that $\bigcup M_\gs$ is a $\cat(-1)$ space.  Hence, $G=\lim G\cs(L)$ is word hyperbolic.
\end{example}

\subsection{Models for spherical Artin groups}\label{ss:hyp}
Suppose $L$ is a nerve of a Coxeter system $(W,S)$ and let $A\cs(L)$ be the associated Artin complex as in Example~\ref{ex:art}.  
For any simplex $\gs$ of $L$, there is a  corresponding spherical Coxeter group $W_\gs$ and a spherical Artin group $A_\gs$. In this subsection we construct a model for $BA_\gs$ by a manifold with boundary.
The group $W_\gs$ is an orthogonal linear reflection group on $\rr^{d(\gs)+1}$, where $d(\gs)=\dim\gs$.   By complexification, $W_\gs$ acts on $\cc^{d(\gs)+1}$. If $\ca_\gs$ denotes the arrangement of reflecting hyperplanes in $\cc^{d(\gs)+1}$, then by Deligne's Theorem in \cite{deligne}, the arrangement complement $M(\ca_\gs)$ is a manifold model for the classifying space of the pure Artin group $PA_\gs$ (where $PA_\gs$ denotes the kernel of $A_\gs\to W_\gs$). Hence, $M(\ca_\gs)/W_\gs$ is a manifold model for $BA_\gs$. Similarly, if $S(\ca_\gs):= S^{2d(\gs)+1}\cap M(\ca_\gs)$, where $S^{2d(\gs)+1}$ denotes the unit sphere in $\cc^{d(\gs)+1}$, then $S(\ca_\gs)/W_\gs$ is a model for $BA_\gs$ by a manifold of dimension $2d(\gs)+1$. Our actual approach will be to define a certain $W_\gs$-invariant bordification of $S(\ca_\gs)$ (which we will denote by the same symbol) and then use $N_\gs:=S(\ca_\gs)/W_\gs$. \nota{$N$?}

\paragraph{Hyperplane arrangements.}
A \emph{hyperplane arrangement} $\ca$ is a finite collection of affine hyperplanes in $\cc^n$. The arrangement is \emph{central} if $\bigcap_{H \in \ca} H$ is nonempty. Its \emph{rank}, $\rk(\ca)$ is the maximum codimension of any nonempty intersection of hyperplanes in $\ca$. An arrangement is \emph{essential} if its rank is $n$. A \emph{subspace} of $\ca$ is either the ambient space $\cc^n$ or a nonempty intersection of hyperplanes. The set of subspaces of $\ca$, partially ordered by reverse inclusion, is denoted $\cq(\ca)$ and is called the \emph{intersection poset}. So, if $F,E\in \cq(\ca)$, then $F<E$ $\iff$ $F\supset E$.

Suppose $\ca'$ and $\ca''$ are arrangements in $\cc^{n'}$ and $\cc^{n''}$, respectively. Define $\ca' \times \ca''$ to be the  arrangement  in $\cc^{n'} \times \cc^{n''}$  consisting of all hyperplanes of the form $H'\times \cc^{n''}$ or $\cc^{n'}\times H''$ for $H'\in \ca$, $H'' \in \ca''$. An arrangement $\ca$ is \emph{reducible} if it is isomorphic to one which admits a nontrivial product decomposition as above. Otherwise, it is \emph{irreducible}.  Note  that the product of two central arrangements is central.

The \emph{codimension}  $c(E)$ of a subspace $E$ in $\cq(\ca)$ is the complex dimension of a complementary subspace $E^\perp\subset \cc^n$. Given $E\in \cq(\ca)$ its \emph{normal arrangement} $\ca_E$  is defined by
	\[
	\ca_E:= \{H\mid H\in \ca \text{ and } H\le E \}.
	\]
(Often we will identify $\ca_E$ with the essential, central arrangement in $E^\perp$ obtained by intersecting the hyperplanes with the orthogonal complement $E^\perp$ of $E$ in $\cc^n$.) There is also an arrangement $\ca^E$ in $E$, called the \emph{restriction} of $\ca$ to $E$, defined by 
	\[
	\ca^E:=\{H\cap E\mid H\cap E \text{ is a hyperplane in $E$}\}. 
	\]

\paragraph{Arrangement complements and their bordifications.} The \emph{arrangement complement} is
\[
M(\ca) := \cc^n - \bigcup \ca.
\]
Suppose $\ca$ is an essential, central arrangement in $\cc^n$ with $\{0\}$ the maximum element of $\cq(\ca)$. 

Put
	\[
	S(\ca):= S^{2n-1} \cap M(\ca)\quad \text{and } D(\ca):=D^{2n} \cap M(\ca),
	\]
where $S^{2n-1}$ and $D^{2n}$ denote the unit sphere and disk in $\cc^n$.
Next we want to attach boundaries to each of these manifolds to obtain a manifold with corners. The idea is to remove a tubular neighborhood of each $E\in \cq(\ca)$, starting with the $E$ of smallest dimension. There is a canonical way to do this which we describe below.

Suppose $V\to E$ is a vector bundle over a manifold $E$. If $s:E\to V$ denotes the $0$-section, define the associated \emph{sphere bundle} $S(V)$ and \emph{cylinder bundle} $C(V)$ by
	\[
	S(V):= (V-s(E))/ \rr_+, \quad\text{and}\quad C(V):=(V-s(E))\times_{\rr_{+}} [0,\infty),
	\]
where $\rr_+$ is the group of positive real numbers. 
So, if $V\to E$ has fiber $\rr^m$, then the fiber of $C(V)\to E$ is the cylinder $S^{m-1}\times [0,\infty)$. Note that $C(V)$ is a manifold with boundary with interior $V-s(E)$. Next, suppose that $X$ is a manifold, that $E \subset X$ is a submanifold and that $V_E$ is the normal bundle. Let $f:V_E\to X$ be a tubular neighborhood. Define $X\odot E:=C(V)\cup X-E$ where $C(V)$ is glued onto $X-E$ via the restriction of the tubular map to the open subset $V_E - s(E)$. The manifold with boundary $X\odot E$ is called the \emph{blowup} of $X$ along $E$; it is formed from $X-E$ by adding the sphere bundle $S(V_E)$ as boundary.

Next we define a bordification of $M(\ca)$, called the \emph{blowup} of $\cc^n$ along $\ca$. Start with $\cc^n$ and then blowup the subspaces $E$ of minimum dimension to obtain a manifold with boundary. Each element of $F\in \cq(\ca)_{<E}$ is also blown up to a submanifold with boundary. We continue by blowing up subspaces of increasing dimension to obtain a manifold with corners, which we continue to denote by $M(\ca)$. (A similar procedure is described in \cite{d78}*{Ch.\ IV}.) Bordifications of $S(\ca)$ and $D(\ca)$ are defined in the same fashion.  The boundary of $M(\ca)$ is a union of manifolds with boundary, $\partial _EM(\ca)$, indexed by the proper subspaces $E\in \cq(\ca)$, where $\partial_EM(\ca)$ means the part of the boundary which results from blowing up $E$, i.e.,
	\begin{equation}\label{e:partialE}
	\partial_EM(\ca) = S(\ca_E) \times \hat{E},
	\end{equation}
where $S(\ca_E)$ is the blowup of the normal arrangement in the unit sphere of $E^\perp$ and $\hat{E}= M(\ca^E)$ is the blowup of $E$ along $\ca^E$. (Note that the right hand side of \eqref{e:partialE} is a product since the normal bundle of $E$ in $\cc^n$ is trivial.) Thus, $S(\ca_E)$ is a submanifold of $\partial M(\ca)$ and for each $F<E$, $S(\ca_F) \cap E$ is also a submanifold of $\partial_F(\hat{E})$.

If the central arrangement $\ca$ decomposes as $\ca=\ca_1\times \cdots \times\ca_k$, then there is an obvious homotopy equivalence $S(\ca)\sim S(\ca_1)\times\cdots \times S(\ca_k)$.

\paragraph{Reflection arrangements.}
Now suppose $(W_\gs, S_\gs)$ is a spherical Coxeter system, where $S_\gs=\vertex(\gs)$. Let $\ca^\rr_\gs$ be the associated real hyperplane arrangement in $\rr^{d(\gs)+1}$ and $\ca_\gs$  its complexification. The intersection posets $\cq(\ca^{\rr}_\gs)$ and $\cq(\ca)$ are canonically identified. As usual, $\cs(\gs)$ denotes the face poset of $\gs$. Since $\gs$ is the dual of the fundamental simplex for $W_\gs$ on $S^{d(\gs)+1}$, each face $\gt\le\gs$ is dual to a face which corresponds to a subspace $E^\rr(\gt)\subseteq \rr^{d(\gs)+1}$ or equally well to $E(\gt)\subseteq \cc^{d(\gs)+1}$. This gives an order-preserving injection $\gt\mapsto E(\gt)$ from $\cs(\gs)$ to $\cq(\ca_\gs)$. The subspace $E(\gt)$ is the subspace of $\cc^{d(\gs)+1}$ fixed by $W_\gt$; hence, $W_\gt$ acts as a reflection group on the normal space $E(\gt)^\perp$. By Deligne's Theorem, $S(\ca_\gs)/W_\gs$ is a model for $BA_\gs$ by a manifold with boundary of dimension $2d(\gs)+1$.  Hence, $\actdim A_\gs\le 2d(\gs)+1$.  When $(W_\gs, S_\gs)$ is reducible this estimate can be improved.  So, suppose its irreducible components are $(W_i, S_i)=(W_{\gs_i}, S_{\gs_i})$, where $1\le i\le k$.  Put $\ca_i=\ca_{\gs_i}$ and $n_i=\dim\gs_i +1$.  Then $W=W_1\times \cdots \times W_k$ and $S(\ca_1)/W_1\times\cdots \times S(\ca_k)/W_k$ is a manifold model for $BA_\gs$.  This gives the following.

\begin{Prop}\label{p:sphericalactdim} Suppose the decomposition of a spherical Coxeter group $W_\gs$ into irreducibles is given by $W_\gs=W_1\times \cdots \times W_k$.  Then the spherical Artin group $A_\gs$ has action dimension $\le 2n-k$ ($=\sum 2d(\gs_i) +1$).
\end{Prop}

\paragraph{The Artin complex.}
Consider an Artin complex $A\cs(L)$ for which the $K(\pi,1)$-Question has a positive answer.  Suppose the dimension of $L$ is $d$. As explained in Example~\ref{ex:art}, since the Salvetti complex for $A$ is of dimension $d+1$, $\gdim A=d+1$.  So, on general principles, $\actdim A\le 2d+2$.  As we explained below, one can use the gluing technique of Subsection~\ref{ss:systems} to obtain the same estimate.
To this end define a complex of aspherical manifolds with boundary $\{M_\gs\}_{\gs\in \cs(L)}$ by putting
	\begin{equation}\label{e:ntau}
	M_\gs:= M(\ca_\gs)/W_\gs .
	\end{equation}
It is a manifold with boundary of dimension $2d(\gs)+2$ and a model for $BA_\gs$.  Next we need to define embeddings $i_{\gs\gt}:M_\gt\hookrightarrow \partial  M_\gs$, whenever $\gt<\gs$.  For a fixed $\gs$, $\ca_\gs$ is an arrangement in $\cc^{d(\gs)+1}$. If  $\gt<\gs$, then $E(\gt)$ is a subspace of $\ca_\gs$ and its normal arrangement $\ca_{E(\gt)}$ ($\cong \ca_\gt$) is an arrangement in $E(\gt)^\perp$. Choose a basepoint $x\in \hat{E}(\gs)$ and identify $S(\ca_{E(\gt)})$ with $S(\ca_{E(\gt)})\times x\subset \partial_{E(\gt)} S(\ca_\gs)$. ($S(\ca_{E(\gt)})$ is the fiber of the sphere bundle of the normal bundle of $S(\hat{E}(\gt))$ in $S(\ca_\gs)$.)  Hence, 
we get an embedding $S(\ca_{E(\gt)})\hookrightarrow \partial S(\ca_\gs)$. Since $M(\ca_\gt)=S(\ca_{E(\gt)})\times [0,\infty)$ and $M(\ca_\gs)=S(\ca_\gs)\times [0,\infty)$ by taking the product with the identity map on $[0,\infty)$, we get the embedding $M(\ca_\gt)\hookrightarrow \partial M(\ca_\gs)$. Taking  quotients by $W_\gs$ and $W_\gt$ ($=$ the stabilizer of $E(\gt)$ in $W_\gs$), this induces an embedding $i_{\gs\gt}:M_\gt\hookrightarrow \partial M_\gs$. Here a small remark is needed: if $W_\gt$ and $W_{\gt'}$ are conjugate subgroups of $W_\gs$, then, as defined, $M_\gt$ and $M_{\gt'}$ are the same submanifold  of $\partial M_\gs$; however, if $\gt\neq\gt'$, we want their images to be disjoint. This is easily arranged by picking a different point $x'\in \hat{E}(\gt')$ for $\gt'$, so that $M_\gt$ and $M_\gt'$ will be parallel submanifolds in $\partial M_\gs$. It is then straightforward to define the dual disk $D_\gt$ to $M_\gt$ and isotopies as in \eqref{e:J} so that $\{M_\gt\}_{\gt\in \cs(L)}$ becomes a complex of manifolds with boundary.  Applying Theorem~\ref{t:gluing} we get a manifold with boundary $M$ of dimension $2d+2$ which is a model for $BA$.  This gives an alternate proof for the following.

\begin{Prop}\label{p:artactdim}
Suppose the $K(\pi,1)$-Question has a positive answer for $A\cs(L)$.  Let $A=\lim A\cs(L)$ and $d=\dim L$.  Then
$\actdim A\le 2d+2$.
\end{Prop}

\paragraph{Permutohedra.} When $L$ is EDCE and the manifolds over the vertices of $L$ are closed, it is necessary to use a trick with permutohedra in order to apply the method of Subsection~\ref{ss:systems0} to get sharp upper bounds for the action dimensions.
This trick is already needed in the case of Artin groups, indeed for $\raag$s. In the case of a $\raag$, $A_L$, our complex of groups is given by the fundamental groups of a complex of tori $\{M_\gt\}_{\gt \in \cs(L)}$, where $M_\gt =(S^1)^{\vertex\gt}$. In Corollary~\ref{c:closedgprdt2} we produced a manifold model for $A_L$ of dimension $2d+2$ (where $d=\dim L$) by using as a complex of manifolds with boundary $M'_\gt:=M_\gt\times I^{\vertex \gt}$. When $L$ is EDCE the dimension can be decreased by one using the complex of manifolds with boundary, $N_\gt=M_\gt\times P(\gt)$ where $P(\gt)$ is a certain $d(\gt)$-dimensional polytope called a ``permutohedron". We shall see in Subsection~\ref{ss:EDCE} below that the same trick works for any complex of closed aspherical manifolds.
 
 \begin{definition}\label{d:permutohedron}
 Suppose $\gs$ is a $d$-simplex. The \emph{permutohedron on $\gs$}, denoted by $P(\gs)$, is the $d$-dimensional convex simple polytope obtained by truncating the proper nonempty faces of $\gs$. The facets of $P(\gs)$ are indexed by the elements of the interval $(\emptyset, \gs)$ in $\cs(\gs)$. (A \emph{facet} is a face of codimension one). Denote the facet corresponding to $\gt$ by $\partial_\gt P(\gs)$. Alternatively, $P(\gs)$ can be defined by blowing up the proper faces of $\gs$ by a procedure similar to that described in the paragraph on bordifications in Subsection~\ref{ss:hyp}.
\end{definition}

Whenever $\gt<\gs$, there is a natural inclusion $i_{\gs\gt}: P(\gt)\to \partial_\gt P(\gs)$. Permutohedra arise naturally in blowups of hyperplane complements.

\begin{example}\label{ex:coord}
Suppose $\ca$ denotes the coordinate hyperplane arrangement in $\cc^{d+1}$, defined by $z_i=0$, $0\le i\le d$. Let $\gD$ be the spherical $d$-simplex defined by
$\gD=[0,\infty)^{d+1}\cap S^{d}$ and define $p:S^{2d+1}\to \gD$ by $(z_0,\dots, z_k)\mapsto (|z_1|^2,\dots, |z_k|^2)$. As in subsection \ref{ss:hyp}, let $S(\ca)$ denote the blowup of the coordinate hyperplane arrangement in the unit sphere of $\cc^{d+1}$. The map $p$ induces a map $\hat{p}: S(\ca)\to P(\gD)$ which is the projection map of a trivial bundle with fiber $T^{d+1}$. Hence, the coordinate hyperplane complement $S(\ca)$ is diffeomorphic as a manifold with corners to $T^{d+1}\times P(\gD)$.
\end{example}

\begin{example}\label{ex:gen} More generally, suppose $\ca$ is a central hyperplane arrangement in $\cc^n$ and that $\ca=\ca_0\times \dots \times \ca_k$ is its decomposition into irreducible components. Since $M(\ca)\cong M(\ca_0)\times \cdots \times M(\ca_k)$, we see that $S(\ca)$ is homotopy equivalent to $S(\ca_0)\times \cdots \times S(\ca_k)$.  In fact, by using a permutahedron we get a corresponding diffeomorphism of manifolds with corners,  
	\begin{equation}\label{e:permprod}
	S(\ca)\,\mapright{\cong}\,[S(\ca_0)\times \cdots \times S(\ca_k)]\times P(\gD),
	\end{equation}
where $P(\gD)$ is a permutohedron of dimension $k$.
If $\ca_i$ is an arrangement in $\cc^{n(i)}$, then $\cc^n=\cc^{n(0)}\times \cdots \times \cc^{n(k)}$. A vector in $\cc^n$ is given by $(z(0),\dots, z(k))$, where $z(i)\in \cc^{n(i)}$.
Define $p:S^{2n-1}\to \gD$ by $(z(0),\dots, z(k)) \mapsto (|z(0)|^2,\dots, |z(k)|^2)$. There is an induced map on blowups, $\hat{p}:S(\ca)\to P(\gD)$. Using $\hat{p}$ and the various projections we get \eqref{e:permprod}. 
For example, suppose $W_\gs$ is a spherical Coxeter group and that $W_\gs=W_0\times \cdots \times W_k$ is its decomposition into irreducibles. As in the paragraph above on reflection arrangements, let $\ca_i$ be the reflection arrangement corresponding to $W_i$ and $\ca_\gs$ the arrangement for $W$. As in \eqref{e:ntau}, put $M_\gs=S(\ca_\gs)/W_\gs$ and $M_i=S(\ca_i)/W_i$.  By \eqref{e:permprod}  we have a diffeomorphism of manifolds with corners:
	\[
	M_\gs\,\mapright{\cong}\, [M_0\times \cdots \times M_k]\times P(\gD).
	\]
\end{example}

\subsection{$L$ is EDCE}\label{ss:EDCE}
When $L$ is EDCE (cf.\ Definition~\ref{d:EDCE}) we can decrease by $1$ our estimates of upper bounds for the action dimension in Propositions~\ref{p:gprdt1} and \ref{p:artactdim}, Theorem~\ref{t:complexaspherical} and Corollaries~\ref{c:closedgprdt} and \ref{c:closedgprdt2}.  In each case we apply Theorem~\ref{t:gluing2} of Subsection~\ref{ss:systems0}. We assume throughout this subsection that $\dim L=d$.

We first consider the Artin group case.  When $L$ is EDCE we can improve the upper bound in Proposition~\ref{p:artactdim} to get the following result, first proved in \cite{giang}.

\begin{Prop}\label{p:artinEDCE}
Suppose the $K(\pi,1)$-Question has a positive answer for $A\cs(L)$.  If $L$ is EDCE, then
$\actdim A\le 2d+1$ ($= 2\gdim A - 1$).
\end{Prop}

\begin{proof}[Sketch of Proof]
To prove this we essentially use the same complex of manifolds with boundary as in the paragraph on the Artin complex in Subsection~\ref{ss:hyp}, except that $M(\ca_\gs)/W_\gs$ is replaced by $S(\ca_\gs)/W_\gs$.  So, we start with the complex of aspherical manifolds $\{M^o_\gs\}_{\gs\in \cs^o(L)}$, where $M^o_\gs=S(\ca_\gs)/W_\gs$.
By taking products with other disks, we can  arrange that for all maximal simplices $\gs$, each  $M^o_\gs$ has dimension $2d+1$.  Next embed $L$ in a contractible simplicial complex $L_c$ of the same dimension $d$.
As in Proposition~\ref{p:ext}, there is a trivial extension of this to a complex of manifolds with boundary over $\cs^o(L_c)$. Finally, apply Theorem~\ref{t:gluing2} to get the result.
\end{proof}

Essentially the same argument gives the following improvement of Theorem~\ref{t:complexaspherical}.

\begin{Theorem}\label{t:complexEDCE}
Suppose that $\{M_\gs\}_{\gs \in \cs(L)}$ is a simple complex of closed aspherical manifolds over $\cs(L)$ satisfying the codimension-$m$ conditions. Let $G = \lim G\cs(L)$ be the direct limit of the $\pi_1(M_\gs)$. Also suppose as before that $L$ is a $d$-dimensional flag complex and that the $K(\pi,1)$-Question for $G\cs(L)$ has a positive answer. If $L$ is EDCE, then $$\actdim(G) \le (m+1)(d+1)-1.$$
\end{Theorem}

\begin{proof}[Sketch of Proof]
For each $\gs\in \cs^o(L)$, let $N^o_\gs=M_\gs\times P(\gs)$, where $P(\gs)$ is a permutohedron (cf. Definition~\ref{d:permutohedron}).  By taking products with other disks we can  arrange that for all maximal simplices $\gs$ the dimensions of the $N^o_\gs$ are equal. Then the complex of manifolds with boundary $\{N^o_\gt\}_{\gt \in \cs^o(L)}$ is a complex of manifolds with boundary satisfying the conditions in Subsection~\ref{ss:systems0}.  The proof is finished exactly as in the proof of the previous proposition:  embed $L$ in a contractible simplicial complex $L_c$ of the same dimension $d$; there is a trivial extension of $\{N^o_\gs\}$ to a complex of manifolds with boundary over $\cs^o(L_c)$; then use Theorem~\ref{t:gluing2} to complete the proof.
\end{proof}

For the graph product complex of fundamental groups, Theorem~\ref{t:complexEDCE} has the following corollary (cf. Corollary~\ref{c:closedgprdt2}), which gives part (ii) of Theorem B in the Introduction. 

 \begin{corollary}\label{c:gprdt2EDCE}
Suppose  $L$ is a flag complex, that for each $v\in \vertex L$, $G_v$ is the fundamental group of a closed aspherical $m$-manifold $M_v$ and that $G=\prod_{L^1} G_v$ is the graph product. If $L$ is EDCE, then
 	\[
 	\actdim G\le (m+1)(d+1)-1
	 \]
 \end{corollary}

The case of a RAAG is where $m = 1$.  Then, it is proved in \cite{ados} that for $d \ne 2$, the EDCE condition, can be replaced with the weaker condition that $H_d(L, \zz_2) = 0$. 

\begin{remark}
On the other hand, Corollaries~\ref{c:closedgprdt2} and \ref{c:gprdt2EDCE} have advantages  over the results of \cite{ados}.  For example, suppose $\dim L =d$ and that $G$ is the graph product over $L^1$ of free abelian groups of rank $m$.  Then $\gdim G= m(d+1)$.  The group $G$ also is the $\raag$ associated to the simplicial complex $\bar{L}$, which is the polyhedral join over $L$ of $(m-1)$-simplices.  (The notion of a ``polyhedral join'' is defined in Definition~\ref{d:polyjoin} below.)  So, $\dim \bar{L} =m(d+1)-1$.  Corollary~\ref{c:closedgprdt2}  yields, $\actdim G\le (m+1)(d+1)$, while \cite{ados} only gives, $\actdim G\le 2\dim \bar{L} +2 = 2m(d+1)$.  (When $L$ is EDCE both estimates can be improved by $1$:  Corollary~\ref{c:gprdt2EDCE} gives $\actdim G\le (m+1)(d+1)-1$ while \cite{ados} gives $\actdim G\le 2m(d+1)-1$).
\end{remark}

\section{Obstructors}\label{s:obstructors}

The \emph{ordered $2$-point configuration space} $\ti\cac(X)$ of a space $X$ is the space of ordered pairs of distinct points in $X$, i.e., $$\ti\cac(X)= (X\times X)-D,$$ where $D$ denotes the diagonal (in many situations it will be better to remove a neighborhood of the diagonal).
The \emph{$2$-point configuration space} $\cac(X)$ is the space of unordered pairs of distinct points:
	\[
	\cac(X):=\ti\cac(X)/\zz_2,
	\]
where $\zz_2$ acts by switching the factors. The double cover $\ti\cac(X)\to \cac(X)$ is classified by a map $c:\cac(X)\to \rr P^\infty=B\zz_2$. Let $w_1$ denote the nontrivial element of $H^1(\rr P^\infty, \zz_2)$, so that $(w_1)^n$ is the nontrivial element of $H^n(\rr P^\infty,\zz_2)$.

\begin{definition}\label{d:vk1}
The \emph{$\zz_2$-valued van Kampen obstruction for $X$ in dimension $n$},  $\vktwo^n(X) \in H^n(\cac(X);\zz_2)$, is defined by $\vktwo^n(X)=c^*(w_1)^n$.
\end{definition}

If $\vktwo^n(X)\ne 0$, then we say $X$ is an \emph{$n$-obstructor}.

\begin{Remark}
One can also define a $\zz$-valued van Kampen obstruction by replacing $w_1$ by the nontrivial element of $H^1(\rr P^\infty; \zz^-)$ where $\zz^-$ means twisted integer coefficients  (cf.\ \cite{ados}). We will not use this refinement in this paper.
\end{Remark}

\begin{definition}\label{d:coarse}
A map $f:X \rightarrow Y$ between metric spaces is a \emph{coarse embedding} if there exist two nondecreasing functions  $\rho_+, \rho_-: \rr^+ \rightarrow \rr^+$ such that $\lim_{t\rightarrow \infty} \rho_-(t) = \infty$ and
 $$\rho_-(d(x,y)) \le d(f(x), f(y)) \le \rho_+(d(x,y)).$$
\end{definition}

Given a finite simplicial complex $K$, let 
	\[
	\cone_\infty K := K \times [0,\infty)/{K \times \{0\}}
	\]
denote the cone of infinite radius on $K$. Equip $\cone_\infty K$ with a proper metric so that for each pair of disjoint simplices $\gs$ and $\gt$, the distance between $\gs \times [t,\infty)$ and $\gt \times [t,\infty)$ goes to infinity as $t\to \infty$. 

Suppose, for simplicity, that the group $G$ is of type $F$. The method of \cite{bkk} consists of the following two steps. 
	\begin{enumerate1}
	\item
	Find a coarse embedding $\cone_\infty K\to EG$ for a suitable complex $K$.
	\item
	Compute the van Kampen obstruction of $K$ in degree $n$. (By the Linking Lemma of \cite{bkk}*{p.\ 223}, this is an obstruction to coarsely embedding $\cone_\infty K$ in a contractible $(n+1)$-manifold).
	\end{enumerate1}
	
	Putting this together, we have the following definition from \cite{bkk}. 
\begin{definition}
The obstructor dimension of $G$, denoted $\obdim G$, is the maximal $n+2$ such that there exists a complex $K$ with nonzero van Kampen obstruction in degree $n$ and a coarse embedding of $\text{Cone}_\infty(K) \rightarrow EG$. 
\end{definition}	
	
\begin{remarks}\label{remarks} 
(i) Since Bestvina-Kapovich-Kleiner \cite{bkk} want to define obstructors for a general group $G$, instead of using (1) they consider  proper, expanding, Lipschitz maps from the $0$-skeleton of $\cone_\infty K$ to $G$.  When $G$ acts cocompactly on $EG$ this amounts to finding a coarse embedding $\cone_\infty(K) \to EG$.  

(ii)  For a finite complex $K$, the cohomology class $\vktwo^n(K)$ is the classical obstruction for finding a PL embedding of $K$ into $\rr^n$.  Indeed, if $K$ embeds in $\rr^n$, then $\cac(K)$ embeds in $\cac(\rr^n)\sim \rr P^{n-1}$, i.e., $\cac(K)$ classifies into $\rr P^{n-1}$ and hence, $\vk^n(K)=c^*(w_1)^n=0$.  In the converse direction, suppose $F:K\to \rr^n$ is a PL map in general position and that $\gS$ is a mod two $n$-cycle in $\cac(K)$.  Then the result of evaluating the cohomology class on the cycle, $\langle \vktwo^n(K),\gS\rangle$, counts the self-intersections of $F$ which lie in $\gS$. (This is explained in \cite{bkk}*{\S 2.1}.)  In Subsection~\ref{ss:gp} we will use this method of calculating self-intersections of general position maps to prove certain van Kampen obstructions are nonzero.

(iii) If $G$ is a word-hyperbolic group or if $EG$ is a $\cat(0)$-space, then $EG$ has a $\mathcal{Z}$-set compactification $\ol{EG}$.  Put $\partial G:=\ol{EG}-EG$. If a finite simplicial complex $K$ is a subspace of $\partial G$,  we get a coarse embedding $\cone_\infty K\to EG$ by choosing a basepoint in $EG$ and coning off $K$.
\end{remarks}

\paragraph{The coarse van Kampen obstruction.}  A generalization of (1) is used by Yoon in \cite{yoon}. He considers instead coarse embeddings $f:T\to EG$ for some contractible CW complex $T$ with a proper metric into $EG$. We are assuming $BG$ is a finite complex and that $EG$ has a length metric induced from a length metric on $BG$.
As in \cite{bkk}, Yoon considers the van Kampen obstruction for embedding $T$ in a contractible $(n+1)$-manifold. 
The $2$-point configuration space of a contractible $(n+1)$-manifold $W$ is homotopy equivalent to $\rr P^{n}$. If $T$ embeds in $W$, then $\ti\cac(T)\subset \ti\cac(W)$ and hence, $c:\cac(T)\to \rr P^\infty$ factors through $\rr P^{n}$. So, if $\vk^{n+1}(T)\neq 0$, then $T$ does not embed in a contractible $(n+1)$-manifold.  The following lemma of \cite{bkk}*{Lemma 8} is important, at least psychologically, to understanding the case $T=\cone_\infty K$.

\begin{lemma}\label{l:cone}\textup{(The Cone Lemma of \cite{bkk}).}
If $K$ is an $n$-obstructor, then $\cone K$ is a $(n+1)$-obstructor, i.e., $\vk^{n+1}(\cone K)=\vk^n(K)$.
\end{lemma}

One actually needs to show the stronger result that $T$ does not coarsely embed in any contractible $(n+1)$-manifold.  For $T=\cone_\infty K$ this is proved in \cite{bkk} by using the Linking Lemma.  In the case of a more general contractible complex $T$, Yoon considers the ``deleted configuration spaces'' $\ti{\cac}_r(T):=[(T\times T)-N_r(D)]$, where $N_r(D)$ means a $\zz_2$-stable $r$-neighborhood of the diagonal.  One then needs to show that the van Kampen obstruction remains nonzero on the quotients, $\cac_r(T):=\ti{\cac}_r(T)/\zz_2$.  So, we define the \emph{coarse van Kampen obstruction} to be the image of $\vk^{n+1}(T)$ in $\lim_{r\to \infty} H^{n+1}(\cac_r(T);\zz_2)$. We define the \emph{proper obstructor dimension} of $G$, denoted $\pobdim(G)$, to be the maximal $n +1$ such that there is a coarse embedding of a contractible complex $T$ into $EG$ such that the coarse van Kampen obstruction of $T$ in dimension $n$ is nonzero. Yoon shows that $\pobdim(G) \le \actdim(G)$. 

\subsection{Configurations of subgroups and sheets}\label{ss:config}

\begin{definition}\label{d:polyjoin}
Suppose $\{X_v\}_{v\in V}$ is a collection of spaces indexed by the vertex set $V$ of a simplicial complex $L$.  For $\gs \in \cs^o(L)$, let $X(\sigma)$ denote the join $\aster_{v \in \vertex \sigma} X_v$. Define the \emph{polyhedral join over $L$} of the $\{X_v\}_{v \in V}$ by 
	\begin{equation}\label{e:poly}
	\aster_L X_v := \bigcup_{\gs \in \cs^o(L)} X(\sigma).
	\end{equation}	
\end{definition}

As in Subsection~\ref{ss:posetgps}, suppose that $G\cq = \{G_\gs\}_{\gs\in \cq}$ is a developable complex of groups with $\vert\cq\vert$ simply connected, that $BG\cq$ is its aspherical realization and that $G=\pi_1(BG\cq)=\lim G\cq$.  Put 
	\begin{equation}\label{e:cg}
	\cg:=\bigcup_{\gs\in \cq} G_\gs
	\end{equation}
and call it a \emph{configuration of standard subgroups}. Suppose each finite subset $\{\gs_1,\dots \gs_k\}$ has a greatest lower bound, $\bigcap\gs_i$. With the word metrics, each inclusion $G_\gs\hookrightarrow G$ is a coarse embedding. Furthermore, given any two subgroups $H_1$ and $H_2$ of $G$, the coarse intersection of $H_1$ and $H_2$ in $G$ is coarsely equivalent to their actual intersection $H_1 \cap H_2$ (see \cite{msw03} for precise definitions of the coarse intersection). Since $\bigcap G_{\gs_i}= G_{\bigcap \gs_i}$, this implies the inclusion of the configuration $\cg$ into $G$ is also a coarse embedding. 

Let $EG\cq$ denote the universal cover of $BG\cq$. For each $\gs$, choose a basepoint $b'_\gs\in BG_\gs$ and a path connecting it to the basepoint $b'\in BG\cq$ ($BG_\gs$ is a subcomplex of $BG\cq$). Choose a lift of $b'$ to a basepoint $b\in EG\cq$. The path from $b'$ to $b'_\gs$ lifts to a path from $b$ to a point $b_\gs$. Let $EG_\gs$ denote the component of the inverse image of $BG_\gs$ in $EG\cq$ containing the basepoint $b_\gs$ (so, $EG_\gs$ is a copy of the universal cover of $BG_\gs$). Put
	\begin{equation}\label{e:ecg}
	E\cg:= \bigcup_{\gs\in \cq} EG_\gs
	\end{equation}
and call it a \emph{standard configuration of sheets} in $EG\cq$. Identify the $G_\gs$-orbit of $b_\gs$ with the group $G_\gs$. If each $G_\gs$ is type $F$, then $E\cg$ is quasi-isometric to the union of orbits $\cg b$. Hence, $E\cg\hookrightarrow EG\cq$ also is a coarse embedding in the type $F$ case.

\begin{example}\label{ex:raags} (\emph{$\raag$s}). As in Example~\ref{ex:art}, suppose $A_L$ is the $\raag$ associated to a flag complex $L$ and that $A\cs(L)$ is the Artin complex.  For each $\gs\in \cs(L)$, let $\rr^\gs$ denote the Euclidean space with basis $\vertex \gs$, let $\zz^\gs\subset \rr^\gs$ be the integer lattice and let $T^\gs=\rr^\gs/\zz^\gs$ be the torus.  The local group $A_\gs$ is $\zz^\gs$. The spaces $BA_\gs$ and $EA_\gs$ are identified with $T^\gs$ and $\rr^\gs$, respectively.   Let $d(\gs)=\dim \gs$. The \emph{octahedron on $\gs$}, denoted $O\gs$, is the $(d(\gs)+1)$-fold join of $S^0$'s, where the copies of $S^0$ are indexed by $\vertex \gs$. So, each \emph{standard sheet} $\rr^\gs$ is identified with $\cone_\infty O\gs$ and the configuration of standard sheets $E\cg$ is identified with $\cone_\infty OL$, where, as in Definition~\ref{d:polyjoin}, $OL$ denotes the polyhedral join of $0$-spheres, i.e.,
	\begin{equation}\label{e:ol}
	OL:=\aster_L S^0 = \bigcup_{\gs\in\cs^o(L)} O\gs.
	\end{equation}
\end{example}

Sometimes we shall consider  finer configurations of subgroups.  Suppose each local group $G_\gs$ contains a simple complex of groups $\cq_\oslash(\gs)$, where the poset $\cq_\oslash(\gs)$ has  geometric realization homeomorphic to that of $\cq_{\le \gs}$.  Moreover, the union of these posets will be the poset of simplices in a simplicial complex $\cq_\oslash$ with the same geometric realization as $\cq$. For each $\ga\in \cq_\oslash(\gs)$,  the corresponding local group $H_\ga$ is a subgroup of $G_\gs$. Then $\ch(\gs) :=\bigcup H_\ga$ is a configuration of subgroups in $G_\gs$ and $E\ch(\gs):=\bigcup EH_\ga$ is a configuration of sheets in $EG_\gs$.  Taking the union over all $\gs\in\cq$, we get a configuration of subgroups $\ch:=\bigcup \ch(\gs) \subset G$ and a configuration of sheets $E\ch:=\bigcup E\ch(\gs) \subset EG$. The archetype is the case of a general Artin group $A$ considered in \cite{dh16}. This is explained in the next example. 
\begin{example}\label{ex:abelartin} (\emph{The configuration  of standard free abelian subgroups in an Artin group}).
As in Example~\ref{ex:art}, let $A\cs(L)$ denote the Artin complex associated to a Coxeter system $(W,S)$ with nerve $L$. Let $A_L$ be the associated Artin group. For each $\gs\in \cs(L)$, $W_\gs$ is the corresponding spherical Coxeter group and $\mathbf{D}_\gs$ is its Coxeter diagram.  As explained in \cite{dh16},  there is a subdivision $\gs_\oslash$ of $\gs$ whose vertices correspond to the connected  subdiagrams of $\mathbf{D}_\gs$.  In other words the vertices of $\gs_\oslash$ are the irreducible special subgroups of $W_\gs$.  Corresponding to each vertex we have the infinite cyclic group generated by the element $\gD^2_\gs$ in the pure spherical Artin group $PA_\gs$. So, the group corresponding to a vertex $s$ of $L$ is the square of the corresponding Artin generator $x_s$. The simplices of $\gs_\oslash$  index the \emph{standard free abelian subgroups} of  $A_\gs$. These subdivisions of simplices fit together to give a subdivision $L_\oslash$ of $L$. For example, the edge  $\{s,t\}$ of $L$ is subdivided into two edges exactly when $3\le m_{st}<\infty$. This leads to the configuration of standard free abelian subgroups in $A_L$ and a configuration of \emph{standard flats} in $A\cs(L)$:
	\[
	\bigcup_{\ga\in \cs(L_\oslash)} \zz^\ga \subset A_L \quad\text{and}\quad \bigcup_{\ga\in \cs(L_\oslash)} \rr^\ga \subset EA\cs(L).
	\]
Moreover, the configuration of standard flats is isometric to $\cone_\infty OL_\oslash$, where, as in \eqref{e:ol}, $OL_\oslash$, is a polyhedral join of $0$-spheres. So, when the $K(\pi,1)$-Question for $A\cs(L)$ has a positive answer, $\cone_\infty OL_\oslash$ is coarsely embedded in $EA_L$. 
\end{example}

It follows from \cite{dh16} that the inclusion of standard infinite subgroups corresponding to the vertices of $L_\oslash$ defines a homomorphism $\Phi:A_{L_\oslash}\to A_L$ from the $\raag$ $A_{L_\oslash}$ to the Artin group $A_L$.  This leads us to the following.

\begin{conjecture}\label{conj}
The homomorphism $\Phi:A_{L_\oslash}\to A_L$ is injective.
\end{conjecture}

Conjecture~\ref{conj} is related to a conjecture of Tits which was proved by proved by Paris and Crisp in \cite{cp}.  Let $\gG$ denote the subgraph of $L^1$ consisting of the edges labeled  $2$.  Let $A_\gG$ be the $\raag$ defined by $\gG$.  There is a homomorphism $A_\gG\to A_L$ which sends each generator for $A_\gG$ to the square of the corresponding generator for $A_L$.  Since the flag complex determined by $\gG$ is a full subcomplex of $L_\oslash$,  the Tits Conjecture amounts to the conjecture that the homomorphism $\Phi$ in Conjecture~\ref{conj} is injective.  Conjecture~\ref{conj} is still open even for spherical Artin groups. For example, if $A$ is the braid group $B_4$, the Tits Conjecture states that there is an injective homomorphism from $\zz^2 \ast \zz$ into $B_4$, whereas our conjecture predicts an injective homomorphism from the RAAG,  $A_{\text{Cone}(C_5)}$ to $B_4$, where $C_5$ is a five-cycle. 

There are similar configurations of abelian subgroups for affine hyperplane complements (cf.\ section~\ref{hypcomp}). The configurations of free abelian and nilpotent groups in \cite{bf} as well as the configurations of free abelian subgroups generated by Dehn twists and``Mess subgroups'' in the mapping class group in \cite{despotovic} follow similar lines.

\subsection{The complex $O_m L$}\label{ss:hol}
Given a $d$-dimensional flag complex $L$, let $O_mL$ denote the polyhedral join of $(m-1)$-spheres, as in Definition~\ref{d:polyjoin}, 
i.e., 
	\[
	O_mL:= \aster_L S^{m-1} := \bigcup_{\gs \in \cs^o(L)} O_m\gs,
	\]
where $O_m\gs$ denotes the $(d(\gs)+1)$-fold join of $(m-1)$-spheres, $S_v$, indexed by $\vertex \gs$ (eventually each of these $(m-1)$-spheres will be given a simplicial structure). Thus, $O_m\gs$ is a sphere of dimension $m(d(\gs)+1) -1$. So, the dimension of $O_mL$ is $m(d+1)-1$. We denote this dimension by $\gd_m(L)$ (or simply by $\gd$):
	\begin{equation}\label{e:gd}
	\gd=\gd_m(L):= \dim O_mL=m(d+1)-1
	\end{equation}
Let $OL=O_1L$. In \cite{ados} the van Kampen obstruction of $OL$ was computed in many cases.

\begin{Theorem}\label{ados}\textup{(\cite{ados}).} Let $L$ be any $d$-dimensional flag complex. If $H_d(L;\zz_2)\neq 0$, then $\vktwo^{2d}(OL)\neq 0$. Therefore, $ \actdim A_L = 2\gdim A_L = 2d+2$. 
\end{Theorem}
The idea for the proof of this in \cite{ados} was to construct a specific $2d$-cycle $\gO$ in the $2$-point configuration space $\cac(OL)$ so that $\vk^{2d}(OL)$ evaluates nontrivially on $\gO$. In the following subsections, we will generalize this to $O_mL$. 

The main result in this section is the following.
	
\begin{Theorem}\label{t:vk}
Suppose $L$ is a $d$-dimensional flag complex and that $\gd=\dim O_mL = m(d+1)-1$. If $H_d(L;\zz_2)\neq 0$, then $O_mL$ is a $(\gd+d)$-obstructor.
\end{Theorem}

The proof of Theorem~\ref{t:vk} will occupy Subsection~\ref{ss:omega}. A similar result holds for polyhedral joins of spheres over $L$ when the spheres are allowed be of different dimensions, and the proof is essentially the same as the proof of Theorem~\ref{t:vk}. Note that $\gd+d=m(d+1)-1 + d = (m+1)(d+1)-2$. 

Before proving Theorem~\ref{t:vk}, we note in the following proposition that the vanishing of the van Kampen obstruction of $O_mL$ in degrees higher than $\gd+d$ follows from the gluing constructions as in Corollaries~\ref{c:closedgprdt2} and \ref{c:gprdt2EDCE}.

\begin{Prop}\label{p:vk}
Suppose $L$ is a $d$-dimensional flag complex and that $\gd=\dim O_mL$.
\begin{enumeratei}
\item
$\vk^{\gd+d+1}(O_mL)=0$.
\item
If $L$ is EDCE, then $\vk^{\gd+d}(O_mL)=0$.
\end{enumeratei}
\end{Prop}

\begin{proof}
For each $v\in \vertex L$, choose a closed, nonpositively curved $m$-manifold $M_v$. Put $G_v=\pi_1(M_v)$ and let $\prod_L G_v$ be the graph product of the $\{G_v\}$. As before, $BG\cs(L)$ is the polyhedral product of the $M_v$. Then $\cone_\infty O_mL$ can be identified with the configuration of standard sheets in $EG$. By Proposition~\ref{p:gprdt1}, $BG$ thickens to a manifold of dimension $(m+1)(d+1)= \gd+d+2$. So, $\cone_\infty O_mL$ coarsely embeds into a contractible $(\gd+d+2)$-manifold and hence, $\vk^n(O_mL)=0$ for $n\ge \gd+d+1$, giving statement (i). Similarly, if $L$ is EDCE, then
$BG$ thickens to a $(\gd+d+1)$-manifold, so that $\vk^{\gd+d}(O_mL)=0$.
\end{proof}

\begin{remark}
One can give a different argument for 
Proposition~\ref{p:vk}~(i) by showing directly that $O_mL$ has an embedding of codimension $(d+1)$ in Euclidean space and hence, that $\vk^{\gd+d+1}(O_mL)=0$.
\end{remark}

\subsection{Sheets of contractible manifolds}\label{ss:cmfld}
As in Definitions~\ref{d:simpleclosed} and \ref{d:codim_m}, suppose $\{M_\gs\}_{\gs\in \cs(L)}$ is a simple complex of closed aspherical manifolds over $\cs(L)$ satisfying the codimension-$m$ conditions. Let $G\cs(L)$ be the associated simple complex of groups.  Put $G=\lim G\cs(L)$. Assume the $K(\pi,1)$-Question has a positive answer for $G\cs(L)$.

Each $M_\gs$ is a subspace of $BG$ ($=BG\cs(L)$). Let $\widetilde M_\gs$ be the copy of its universal cover containing a given basepoint $b\in EG$. Let $E\cg=\bigcup \widetilde M_\gs$ be the union of sheets. If each $\widetilde M_\gs$ is $\cat(0)$ and if its visual boundary, $\partial_\infty \widetilde M_\gs$, is homeomorphic to the sphere $O_m\gs$, then $\widetilde M_\gs$ is homeomorphic to $\cone_\infty (O_m\gs) = \rr^{m(d(\gs)+1)}$. Hence, $E\cg$ is homeomorphic to $\cone_\infty (O_mL)$. This uses the assumption that the $\widetilde M_\gt$ intersect transversely; so that for $\gt<\gs$, the visual sphere of $\widetilde M_\gt$ is identified with the standard subsphere $O_m\gt\subset O_m\gs$. By Theorem~\ref{t:vk}, if $H_d(L;\zz_2)\neq 0$, then $\vk^{\gd+d}(O_mL)\neq 0$ and hence, $\obdim G \ge \gd+d +2$. 

Our goal in this subsection is to show how to generalize this without the assumption that  $\widetilde M_\gs$ has a $\mathcal{Z}$-set compactification with boundary a sphere.  In particular,  we do not need to assume that $\widetilde M_\gs$ is simply connected at infinity. So, assuming Theorem~\ref{t:vk} (which will be proved in the next subsection), we here prove the following.

\begin{Theorem}\label{t:pobdimcm}
Suppose that $\{M_\gs\}_{\gs\in \cs(L)}$ is a simple complex of closed aspherical manifolds over $\cs(L)$, where $L$ is a $d$-dimensional flag complex. Take hypotheses and notation as above. If $H_d(L,\zz_2)\neq 0$, then the coarse van Kampen obstruction of $E\cg$ in degree $\gd+d+1$ is nonzero. So, $\pobdim G\ge \gd+d+2$. Hence, $\actdim G=\gd+d+2$.
\end{Theorem}

This theorem applies, for example, when $G$ is the graph product of fundamental groups of closed aspherical manifolds.

To simplify notation write $E_\gs$ instead of $\widetilde M_\gs$ and $C_\gs$ instead of $\cone (O_m\gs)$. Also, write $E$ for $E\cg= \bigcup E_\gs$ and $C$ for $\bigcup C_\gs$. We can identify $C_\gs$ with an open neighborhood of the basepoint in $E_\gs$. Thus, $C$ is an open neighborhood of the basepoint $b$ in $E$ and $C\times C$ is an open neighborhood of $(b,b)$ in $E\times E$.  The inclusion $C\times C\hookrightarrow E\times E$ takes the diagonal to the diagonal so we have a $\zz_2$-equivariant inclusion $\ti{\cac}(C)\hookrightarrow \ti{\cac}(E)$ inducing an inclusion of $2$-point configuration spaces,
	\(
	i:\cac(C)\hookrightarrow \cac(E).
	\)

\begin{lemma}\label{l:CE}
The inclusions $\ti{\cac}(C)\hookrightarrow \ti{\cac}(E)$ and $i:\cac(C)\hookrightarrow \cac(E)$ are both homotopy equivalences.
\end{lemma}

We will assume that each $E_\gs$ comes with a proper $G_\gs$-invariant metric so that the inclusions $E_\gt\hookrightarrow E_\gs$ are isometries. Extend these metrics to a metric on $E$  by taking the induced path metric. This implies that if $x$ and $y$ are points in $E$ such that $x \in E_\sigma$, $y \in E_\tau$, then there is a point $z \in E_{\gs \cap \gt}$ with $d(x,z) + d(y,z) = d(x,y)$.
As in \cite{yoon}, let $N_r$ be an $r$-neighborhood of the diagonal $D$ in $E\times E$. Write $\ti{\cac}_r(E)$ for $(E\times E)-N_r$ and $\cac_r(E)$ for its quotient by the free $\zz_2$-action.

\begin{lemma}\label{l:ENr}
The inclusions $\cac_r (E)\hookrightarrow \cac (E)$ induce an isomorphism on cohomology:
	\[
	H^*(\cac(E))\mapright{\cong}\lim_{r\to \infty} H^*(\cac_r(E)),
	\]
where $\lim$ means direct limit.
\end{lemma}

Both Lemmas~\ref{l:CE} and \ref{l:ENr} have similar proofs. 

Note that $E=\bigcup E_\gs$ is a poset of contractible manifolds over $\cs(L)$. In particular, that $E_\gs\cap E_\gt = E_{\gs\cap \gt}$. Similarly, $E\times E=\bigcup_{(\gs,\gt)} E_\gs\times E_\gt$ is a poset of contractible manifolds over $\cs(L)\times \cs(L)$. The diagonal $D(E)$ intersects the terms in this decomposition as follows: 
\(
(E_\gs\times E_\gt)\cap D(E) = D_{\gs\cap\gt},
\)
where $D_{\gs\cap\gt}$ denotes the diagonal in $E_{\gs\cap \gt} \times E_{\gs\cap \gt}$. It is a properly embedded, contractible submanifold with trivial normal bundle in the contractible manifold $E_{\gs\cap \gt} \times E_{\gs\cap \gt}$. Thus,
	\[
	\ti{\cac}(E):=(E\times E)-D(E) = \bigcup_{(\gs,\gt)}(E_\gs\times E_\gt)-D_{\gs\cap\gt}.
	\] 
The manifold $(E_\gs\times E_\gt)-D_{\gs\cap\gt}$ is homotopy equivalent to a normal sphere $S^{c(\gs,\gt) -1}$, where $c(\gs,\gt)$ is the codimension of $D_{\gs\cap\gt}$ in $E_\gs\times E_\gt$. The normal vector space of $D_{\gs\cap\gt}$ in $E_\gs\times E_\gt$ decomposes as $V^a+V^b+V^d$, where $V^a$ is the normal space of $E_{\gs\cap\gt}$ in $E_\gs$, $V^b$ is the normal space of $E_{\gs\cap\gt}$ in $E_\gt$, and $V^d$ is the normal space of $D_{\gs\cap\gt}$ in $E_{\gs\cap\gt}$. Thus, $S^{c(\gs,\gt)-1}$ has a join decomposition as $S^{a-1}*S^{b-1} *S^{d-1}$. The involution (switching factors) maps $V^a$ and $V^b$ into different factors and acts by the antipodal map on $V^d$. It follows that the image of the normal sphere to $D_{\gs\cap\gt}$ in $(E_\gs\times E_\gt)$ in the $2$-point configuration space $\cac(E)$ is homotopy equivalent to $S^{a-1}*S^{b-1}*\rr P^{d-1}$, i.e., to a suspension of projective space. 

The previous paragraph goes through, mutatis mutandis, for the union of cones $C=\bigcup C_\gs$, as well as, for its ordered $2$-point configuration space, $\ti{\cac}(C)$.

\begin{proof}[Proof of Lemma~\ref{l:CE}]
Both $\ti{\cac}(C)$ and $\ti{\cac}(E)$ are posets of spaces over $(\cs(L)\times\cs(L))_{>(\emptyset,\emptyset)}$. The geometric realization of this poset is the join, $L*L$. The relative homology groups of $(\ti{\cac}(E),\ti{\cac}(C))$ can be computed from a spectral sequence for the poset of spaces, e.g., see \cite {do12}. Its $E^1$-page is
	\[
	E^1_{pq}= C_p(L*L;H_q((E_\gs\times E_\gt) -D, (C_\gs\times C_\gt)-D)).
	\]
(The coefficients in this spectral sequence are not locally constant).
Since $C_\gs\times C_\gt-D \hookrightarrow E_\gs\times E_\gt -D$ is a homotopy equivalence, $H_q(E_\gs\times E_\gt -D, C_\gs\times C_\gt-D)=0$ for all $q$. Hence, $H_*(\ti{\cac}(E),\ti{\cac}(C))$ vanishes in all degrees. When $L$ is connected, $L*L$ is simply connected; so, by van Kampen's Theorem, $\ti{\cac}(C)$ and $\ti{\cac}(E)$ are simply connected. Hence, by Whitehead's Lemma, $\ti{\cac}(C)\hookrightarrow \ti{\cac}(E)$ is a homotopy equivalence. A more careful analysis yields the same statement even when $L$ is not connected. We shall not give the argument since all we need is that the map induces an isomorphism on homology. Since $\ti{\cac}(C)\hookrightarrow \ti{\cac}(E)$ is $\zz_2$-equivariant it induces a homotopy equivalence $\cac(C)\to \cac(E)$.
\end{proof}

To prove Lemma~\ref{l:ENr}, first note that $\ti{\cac}_r(E)$ also has a decomposition as a poset of contractible manifolds:
	\[
	\ti{\cac}_r(E)= \bigcup_{(\gs,\gt)} (E_\gs \times E_\gt)- N_r(D_{\gs\cap \gt}),
	\]
where $N_r(D_{\gs\cap \gt})$ means the $r$-neighborhood of $D_{\gs\cap\gt}$ in $E_\gs\times E_\gt$. By coarse Alexander duality (e.g., see \cite{yoon}), $(E_\gs \times E_\gt)- N_r(D_{\gs\cap \gt})$ has the same pro-homology type as $(E_\gs \times E_\gt)- D_{\gs\cap \gt}$ (which is homotopy equivalent to $S^{c(\gs,\gt)-1}$); so, the inclusions induce an isomorphisms:
	\begin{equation}\label{e:alexdual}
	H^*(S^{c(\gs,\gt)-1})=H^*((E_\gs \times E_\gt)- D_{\gs\cap \gt})\to \lim_{r\to \infty} H^*((E_\gs \times E_\gt)- N_r(D_{\gs\cap \gt})) .
	\end{equation}

\begin{proof}[Proof of Lemma~\ref{l:ENr}] The cohomology spectral sequences for the poset of spaces give spectral sequences with $E_1$-pages:
	\begin{align*}
	E_1^{pq} &= C^p(L*L;H^q ((E_\gs \times E_\gt)- D_{\gs\cap \gt}))  \\
	E_1^{pq}(r) &=C^p(L*L;H^q ((E_\gs \times E_\gt)- N_r(D_{\gs\cap \gt}))) 
	\end{align*}
By \eqref{e:alexdual}, the inclusions induce an isomorphism,
$E_1^{pq} \to \lim E_1^{pq}(r)$
and hence, by the comparison theorem for spectral sequences, an isomorphism, 
	\[
	H^*(\ti{\cac}(E))\mapright{\cong} \lim_{r\to \infty} H^*(\ti{\cac}_r(E)).
	\] 
There is a similar isomorphism for $H^*(\cac_r(E))$.
\end{proof}

\begin{proof}[Proof of Theorem~\ref{t:pobdimcm}]
By Theorem~\ref{t:vk}, $\vk^{\gd+d}(O_mL) \neq 0$.  By Lemma~\ref{l:cone}, this implies $\vk^{\gd+d+1} (\cone(O_mL))\neq 0$. By Lemma~\ref{l:CE}, its image $\vk^{\gd+d+1}(E)\in H^{\gd+d+1}(\cac(E))$ also is not zero. Finally, by Lemma~\ref{l:ENr}, the coarse van Kampen obstruction (i.e., the image of this class in $\lim H^{\gd+d+1}(\cac_r(E))$) is $\neq 0$. Therefore, $\pobdim G\ge \gd+d+2$. Since, by Proposition~\ref{p:gprdt1}, $\actdim G\le \gd+d+2$, the last sentence of the theorem follows.
\end{proof}

\subsection{The van Kampen obstruction and general position}\label{ss:gp}

For a finite simplicial complex $K$ there is an equivalent definition of the van Kampen obstruction in terms of a general position map of $K$ into Euclidean space which we now describe. First, replace $\cac(K)$ by the \emph{simplicial $2$-point configuration space of $K$}:	
	\begin{equation}\label{e:config}
	\mathcal{C}(K) = [(K \times K) - D] / \zz_2,
	\end{equation}
where 
	\(
	D = \{(\sigma,\tau) \in K \times K \mid \sigma \cap \tau \ne \emptyset\}
	\)
is a simplicial thickening of the diagonal.  Since the $2$-point configuration space and the simplicial $2$-point configuration space are homotopy equivalent, we can denote both $\cac(K)$ without risking confusion.

\begin{definition}\label{d:vk2}
Let $K$ be a $k$-dimensional simplicial complex, and let $f: K \rightarrow \mathbb{R}^n$ be a general position map. This means, in particular, that if $\sigma$ and $\tau$ are two disjoint simplices of $K$ with $\dim \gs + \dim \gt = n$, then the images of $\sigma$ and $\tau$ intersect in a finite number of points. The \emph{van Kampen obstruction} $\vktwo^n(K) \in H^n(\mathcal{C}(K); \mathbb{Z}_2)$ is the cohomology class of the cocycle $\cvk$ ($=\cvk^n(K)$) defined by $$\langle \cvk,\{\sigma,\tau\}\rangle = |f(\sigma) \cap f(\tau)| \mod 2,$$where $\{\gs,\gt\}$ means an unordered pair of disjoint simplices in $K$ (we are using $|X|$ to denote the cardinality of a finite set $X$).
\end{definition}

Generalizing \cite{mtw11}*{Appendix D}, we give the following description of a cocycle representing 
$\vktwo^n(K)$.
Given any total ordering of the vertices of $K$, there is a general position map $f$ from $K$ to $\rr^n$ defined by sending the $i^{th}$ vertex in $K$ to $\lambda(i)$, where $\lambda = (t, t^2, \ldots, t^n) \in \rr^n$ is the moment curve, and extending linearly.
Suppose $\gs,\gt \in K$ with $\dim \gs + \dim \gt = n$. The convex hull of the union of  vertices of $f(\gs)$ and $f(\gt)$ is the cyclic polytope $C(n+2,n)$. If $\gs$ and $\gt$ intersect, then neither of them can be contained in faces of $C(n+2,n)$. The faces of $C(n+2,n)$ are completely determined by \emph{Gale's Evenness Condition}, which in this case says that a set $T$ of $n$ vertices of $C(n+2,2)$ spans a face if and only if the two missing elements in Vert$(C(n+2,n)) - T$ are separated by an even number of elements of $T$.

Two simplices $\gs$ and $\gt$ with $\dim(\gs) + \dim(\gt) = n$   are said to be \emph{meshed} if the order on their vertices is either 

\begin{align*}
 v_0& < w_0 < v_1 < w_1 < \dots < v_{n/2} < w_{n/2} \text{\ \ , or} \\
v_0 &< w_0 < v_1 < w_1 < \dots < w_{(n-1)/2} < v_{(n+1)/2}.
\end{align*}

\begin{lemma}
Two simplices $\gs$ and $\gt$ with $\dim(\gs) + \dim(\gt) = n$ intersect under the map $f$ if and only if they are meshed. 
\end{lemma}
\begin{proof}
If $\gs$ and $\gt$ are not meshed, there are two vertices $v_i$ and $v_{i+1}$ of $\gs$ with no vertex of $\gt$ between them. 
In $C(n+2,n)$ the union of the vertices of $f(\gt)$ and all other vertices of $f(\gs)$ except $v_i$ and $v_{i+1}$ spans a face by the evenness condition. 
Therefore, $f(\gt)$ is contained in a face of the cyclic polytope and so cannot intersect $f(\gs)$.

If $\gs$ and $\gt$ are meshed, then the evenness condition implies that $f(\gs)$ and $f(\gt)$ are not proper faces of $C(n+2,n)$. 

Suppose that $f(\gs)$ and $f(\gt)$ do not intersect. Let $H$ be a hyperplane separating $f(\sigma)$ and $f(\tau)$, so that $H$ partitions the vertices of $C(n+2,n)$ into $\vertex f(\sigma)$ and $\vertex f(\tau)$. If $f(\sigma)$ (or $f(\tau)$) is in the interior of $C(n+2,n)$, then another vertex of $C(n+2,n)$ is on the same side of $H$; hence, $f(\sigma) \subset \partial C(n+2,n)$,  a contradiction.
\end{proof}

Note that if the difference between the dimensions of $\sigma$ and $\tau$ is greater than one, then $f(\sigma)$ and $f(\tau)$ are disjoint.

\subsection{Proof of Theorem~\ref{t:vk}}\label{ss:omega}
We recall the construction in \cite{ados}. Suppose that $L$ is a $d$-dimensional complex with $H_d(L; \zz_2) \ne 0$, and suppose $C$ is a $d$-cycle in $L$ with coefficients in $\zz_2$. Identify $C$ with its support. 
Choose a $d$-simplex $\gD\in C$ with vertices $v_0,\dots, v_d$.
Let $v_i^{\pm}$ denote the two vertices in $OC$ lying above $v_i$.
Let $D^C(\gD)$ be the full subcomplex of $OL$ containing $C^-$ and the vertices $v_0^+,\dots, v_d^+$ of $\gD^+$. 
We say $D^C(\gD)$ is $C$ \emph{doubled over the simplex} $\gD$. Suppose $\ga, \gb$ are disjoint $d$-simplices in $D^C(\gD)$.
Define a chain $\gO \in C_{2d}(\cac(D^C(\gD));\zz_2)$ by declaring the $2d$-cell $\{\ga,\gb\}$ of $\cac(D^C(\gD))$ to be in $\gO$ if and only if
	\begin{itemize}
	\item $\ga \cap \gb = \emptyset$, and
	\item $\vertex \Delta \subset p(\vertex \ga) \cup p(\vertex \gb)$.
	(Here $p:\vertex OL \to \vertex L$ is the natural projection.)
	\end{itemize}

It is proved in \cite{ados} that $\gO$ is a cycle and that $\cvk^{2d}(D^C(\gD))$ evaluates nontrivially on $\gO$. 

Next, we define a subcomplex $D_m^C(\gD)$ of $O_mL$. We assume that each sphere in $O_mL$ is triangulated as the boundary of an $m$-simplex. Let $D_m^C(\gD)$ be the full subcomplex of $O_mL$ containing $\vertex(C) \times \{1\} \hspace{.5mm} \bigcup \hspace{.5mm} \vertex(\Delta) \times \{ 2, 3, \dots m\}$. So, $D_m^C(\gD)$ is constructed by replacing each vertex $v \in \Delta$ with the boundary of a $m$-simplex. As before, let $\delta = \delta_m(L)$  denote the dimension of $O_mL$.

Surprisingly, the above definition of $\gO$  works in the case of $O_mL$.  Define a chain $\gO_m$ in $C_{d+\delta}(\cac(D_m^C(\gD));\zz_2)$ to be the union of all cells $\{\sigma, \tau\}$ such that 
	\begin{itemize}
	\item $\sigma \cap \tau= \emptyset$ 
	
	\item $\vertex \Delta \subset p(\vertex \sigma) \cup p(\vertex \tau)$.
	\end{itemize}
 
It  will be shown in Theorem~\ref{l:cycle} below that $\gO_m$ is a cycle. We first need a few lemmas which restrict the possible $(d+\delta)$-cells in $\cac(D^C_m(\gD))$.

\begin{definition}
For any $w \in \vertex \Delta$ and $\sigma, \tau \in O_mL$, let $M_w^{\sigma\tau}$ be the collection of \emph{missing} vertices in $p^{-1}(w)$, i.e., $M_w^{\gs\gt}$ is the set of vertices in $p^{-1}(w)$ that are not contained in $\sigma \cup \tau$. 
\end{definition} 

Note that if $p(\sigma)$ misses a vertex $w$ in $\vertex \Delta$, then $|M^{\gs\gt}_w| \ge 1$ for any $\tau$, since the preimage $p^{-1}(w)$ does not span a simplex in $O_m(L)$. 

\begin{lemma}\label{l:missing}
If $\{\sigma,\tau\} \in \gO_m$, then for any $w \in \vertex \Delta, |M_w^{\gs\gt}| \le 1$.
\end{lemma}

\begin{proof}
The cardinality of $\vertex (\sigma \cup \tau)$ is the same as that of $\vertex p^{-1}(\Delta)$. Assume $\vertex (p(\sigma) \cup p (\tau))$ includes $l$ vertices not contained in $\Delta$, so that $\vertex (\sigma \cup \tau)$ contains $d+ \delta
 - l$ vertices in $p^{-1}(\Delta)$. If $\vertex \sigma$ contains a vertex outside$\Delta$, then $p(\sigma)$ misses a vertex $w$ of $\Delta$, which implies that $|M_w^{\gs\gt}| \ge 1$. Furthermore, if $v$ and $v'$ are distinct vertices in $L - \Delta$ which are contained in $\vertex \sigma \cup \tau$, then there are distinct vertices $w$ and $w'$ so that $|M_w^{\gs\gt}|$ and $|M_{w'}^{\gs\gt}|$ are both $\ge 1$. Otherwise, $\sigma$ and $\tau$ would both miss a vertex $w \in \Delta$ and $\{\sigma,\tau\}$ would not be contained in $\gO_m$. 
Similarly, if there are $l$ vertices which are not contained in $\Delta$, then there are $l$ vertices $w_1, \dots, w_l$ of $\Delta$ with $|M_{w_i}^{\gs\gt}|\ge 1$. Therefore, if $|M_w^{\gs\gt}| > 1$ for any $w$, the number of total missing vertices in $p^{-1}(\Delta)$ is greater than $l$, a contradiction.
\end{proof}

\begin{lemma}\label{verts} 
If $\{\sigma,\tau\} \in \gO_m$ then $p(\sigma)$ and $p(\tau)$ are in $L^{(d)}$. 
\end{lemma}

\begin{proof}
Assume $\vertex (p(\sigma) \cup p(\tau))$ includes $l$ vertices not contained in $\Delta$, so that $p(\gs) \cup p(\gt)$ contains $d+l+1$ vertices. The proof of Lemma \ref{l:missing} implies there are $l$ vertices $w_1, \ldots, w_l$ of $\Delta$ such that $|M_{w_i}^{\gs\gt}| = 1$. For the other $d +1 - l$ vertices of $\vertex \Delta$, $|M_w^{\gs\gt}| = 0$. Neither $\vertex \sigma$ nor $\vertex \tau$ can contain $\vertex p^{-1}(w)$, so each such $w$ is contained in $p(\sigma)$ and $p(\tau)$. Thus, $p(\gs) \cap p(\gt)$ contains at least $d + 1 - l$ vertices. Since $|\vertex p(\gs)|$ and $|\vertex p(\gt)|$ is bounded above by $d+1$, the equality $$|\vertex p(\gs) \cup \vertex p(\gt)|= |\vertex (p(\sigma) \cup p(\tau))| + |\vertex (p(\gs) \cap p(\gt))|$$ implies that $|\vertex p(\gs)|$ and $|\vertex p(\gt)|$ both equal $d+1$.
\end{proof}

The next two theorems are our computation of the van Kampen obstruction of $O_mL$. 

\begin{Theorem}\label{l:cycle} 
Let $L$ be a $d$-dimensional flag complex with $H_d(L, \zz_2) \ne 0$. Let $C$ be a $d$-dimensional cycle contained in $L$, and $\Delta \subset C$ a $d$-simplex. Then $\gO_m \in C_{d+\delta}(\cac(D_m^C(\gD));\zz_2)$ is a $(d + \delta)$-cycle.
\end{Theorem}

\begin{proof}
We assume that $m > 1$, since the $m = 1$ case was proved in \cite{ados} (this slightly simplifies the argument).
Let $\{\sigma,\alpha\}$ be a $(d + \delta - 1)$-cell in $\cac(D_m^C(\gD))$.
 We claim the sum of the cardinality of the sets
 \begin{align*}
 \mathcal{V}_1 &:= \{v \in \vertex D_m^C(\gD) | \{\sigma \ast v, \alpha\} \in \gO_m\}\\ 
 \mathcal{V}_2 &:= \{v \in \vertex D_m^C(\gD) | \{\sigma, \alpha \ast v\} \in \gO_m\}
 \end{align*}
 is even.  Note that some vertices of $D_m^C(\gD)$ may be contained in both $\mathcal{V}_1$ and $\mathcal{V}_2$.

First, suppose $p(\sigma)$ and $p(\alpha)$ are in $C^{(d)}$. In this case, if $v \in \mathcal{V}_i$, then $p(v) \in \Delta$. By Lemma \ref{verts}, we can assume $|M_w^{\gs\ga}| = 0,1$ or $2$ for all $w \in \Delta$; otherwise, $\mathcal{V}_1$ and $\mathcal{V}_2$ would be empty. 

If $\vertex \sigma\cap \vertex p^{-1}(w)\ne \emptyset$ and $\vertex \alpha \cap \vertex p^{-1}(w) \ne \emptyset$, then each vertex of $M_w^{\gs\ga}$ is in $\mathcal{V}_1$ and $\mathcal{V}_2$, hence $w$ contributes an even number to the sum of $|\mathcal{V}_1|$ and$|\mathcal{V}_2|$. If $\vertex \sigma \cap \vertex p^{-1}(w) =\emptyset$ and $\vertex \alpha \cap \vertex p^{-1}(w) = \emptyset$, then $\mathcal{V}_1$ and $\mathcal{V}_2$ are empty. 

Next, suppose that $\vertex \sigma \cap \vertex p^{-1}(w) \ne \emptyset$ and $\vertex \alpha \cap \vertex p^{-1}(w) = \emptyset$, so that $M_w^{\gs\ga}$ makes no contribution to $\mathcal{V}_2$.
If $|M_w^{\gs\ga}| = 1$, then the missing vertex is not in $\mathcal{V}_1$ since $p^{-1}(w)$ does not span a simplex in $O_mL$. If $|M_w^{\gs\ga}| = 2$, then each vertex in $M_w^{\gs\ga}$ is contained in $\mathcal{V}_1$, so again $M_w^{\gs\ga}$ contributes an even number to $|\mathcal{V}_1|$. The same argument works if $\vertex \sigma \cap \vertex  p^{-1}(w)= \emptyset$ and $\vertex \alpha \cap \vertex p^{-1}(w) \ne \emptyset$.

Now, assume $p(\sigma)$ is a $d$-simplex of $L$ and $p(\alpha)$ is a $(d-1)$-simplex of $L$. In this case, $\mathcal{V}_1$ is empty by Lemma \ref{verts}.
Again, we consider the sets $M_w^{\gs\ga}$ for $w \in \vertex \Delta$.

Note that $|M_w^{\gs\ga}| = 2$ for at most one $w \in \vertex \Delta$ by Lemma \ref{l:missing}, and if $|M_w^{\gs\ga}| = 2$ for some $w \in \Delta$, then $\mathcal{V}_2$ is contained in $p^{-1}(w)$. In this case, there are $0$ or $2$ vertices in $\mathcal{V}_2$, depending on whether or not $w$ is in the link of $\alpha$. 
Suppose $|M_w^{\gs\ga}| \ne 2$ for all $w \in \Delta$. Since $C$ is a cycle and $p(\alpha)$ is $(d-1)$-dimensional, the link $\Lk_C(p(\alpha))$ is an even number of vertices. For each $w \in \Lk_C(p(\alpha)) \cap \Delta$, by assumption there is precisely one vertex in $\mathcal{V}_2$ which is in $M_w^{\ga\gs}$ (if there were zero vertices, then $\gs$ would contain all of $p^{-1}(w)$). 

Now, we claim that all the vertices in $\Lk_C(p(\alpha)) - \Delta$ are in $\mathcal{V}_2$. Such a vertex $v$ is not in $\mathcal{V}_2$ if and only if it is contained in $\gs$. Since $p(\gs) \cup p(\ga)$ contains $\Delta$, if $v$ were in $\gs$ this would imply that $\Delta \cup v$ is a simplex in $L$, which contradicts $L$ being $d$-dimensional and flag. Therefore, each vertex in $\Lk_C(p(\alpha)) - \Delta$ is in $\mathcal{V}_2$, and since $|M_w^{\ga\gs}| = 1$ for all $w \in  \Lk_C(p(\alpha)) \cap \Delta$, the total cardinality of the set $\mathcal{V}_2$ is even (and equal to the cardinality of $\Lk_C(p(\alpha)))$. 
 \end{proof}

\begin{Theorem}\label{t:neq0}
	Let $C$ be the support of a cycle in $H_d(L;\zz_2)$, let $D_m^C(\Delta)$ be $C$ doubled over a $d$-simplex $\gD \subset C$ and let $\gO_m\in Z_{d + \delta }(\cac(O_mL);\zz_2)$ be as above.
	Then $\cvk^{d+\delta}$ evaluates nontrivially on $\gO_m$.	
\end{Theorem}

\begin{proof}
	Order the vertices of $\gD$ so that $v_0<\cdots < v_{d+1}$ and then order the other vertices of $L$ so that each vertex of $\gD$ is $<$ each vertex of $L-\gD$. Extend this to an ordering on the vertex set of $O_mL$, by 
$$v_0^1 < v_0^2 < \dots < v_0^{m} < v_1^1 < v_1^2 < \cdots < v_1^{m} < \dots < v_{d+1}^{1} < \dots < v_{d+1}^{m}.$$

We have the following decomposition of the obstruction cocycle $\cvk^{d+\delta}$ evaluated on $\gO_m$:
	\[
		\sum_{\{\sigma,\tau\} \in \gO_m}\cvk^{d+\delta}(\{\sigma,\tau\}) = \sum_{\substack{\{a,b\} \in M \\
		\vertex \Delta\subset a \cup b}} \ \sum_{\substack{\{\sigma,\tau\} \in \gO_m \\
		p(\sigma) = a\\
		p(\tau) = b}}\cvk^{d+\delta}(\{\sigma,\tau\}).
	\]
If $a = b = \Delta$, then there is exactly one meshed pair, because the union of vertices of $\gs$ and $\gt$ is precisely the set of vertices of $O_m\Delta$, and for each $m$-simplex in $D_m(\Delta)$ there is a unique pair of meshed faces. Now, suppose $a \ne \Delta$.
	Let $b$ such that $\vertex \Delta\subset a \cup b$, and let $\sigma \in p^{-1}(a)$.
	If $p(\tau) = b$, then $\{\sigma,\tau\} \in \gO_m$ if and only if $\sigma \cap \tau = \emptyset$. For all $w \in \Delta - p(\sigma)$, if $\tau$ contains more than two vertices of $p^{-1}(w)$, then $\sigma$ and $\tau$ do not mesh by our choice of ordering. On the other hand, the vertices of $\sigma \cup \tau$ can omit at most $1$ vertex of $p^{-1}(w)$ by Lemma \ref{l:missing}; hence, no meshing can occur.
	Therefore, the only contribution comes from the unique meshed pair with $a = b = \Delta$, and hence, $\cvk^{d+\delta}$ evaluates nontrivially on $ \gO_m$. 
	\end{proof}

Putting this all together, we get the following.

\begin{Theorem}\label{hatvk}
 Let $L$ be a $d$-dimensional flag complex and let $O_mL$ be a polyhedral join over $L$ of $(m-1)$-spheres. Let $\delta = \dim O_mL$. If $H_d(L;\zz_2)\neq 0$, then $\vktwo^{d+\delta}(O_mL)\neq 0$.
\end{Theorem}

Theorems~\ref{hatvk} and \ref{t:complexaspherical} have the following corollary.

\begin{corollary}\label{c:gp}
Suppose $L$ is a $d$-dimensional flag complex, and $G = \prod_{L^1}G_v$ is a graph product over $L$, where each $G_v$ is the fundamental group of a closed aspherical $m$-manifold. If $H_d(L; \zz_2) \ne 0$, then $$\obdim G = \actdim G = (m+1)(d+1) = \gdim G + (d+1)$$
\end{corollary}

Combining this corollary with Theorem~\ref{t:complexEDCE}, gives Theorem B in the Introduction.

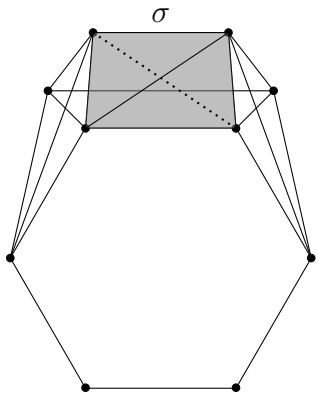
\begin{figure}
\centering

\begin{tikzpicture}[scale = .5]

 \draw (0:4) \foreach \x in {60,120,...,359} {
        -- (\x:4)
      }-- cycle (90:4);
      
  \draw [fill] (4,0) circle [radius=0.1];
  \draw [fill] (-4,0) circle [radius=0.1];
   \draw [fill] (2,3.45) circle [radius=0.1];
  \draw [fill] (-2,3.45) circle [radius=0.1];
  \draw [fill] (2,-3.45) circle [radius=0.1];
  \draw [fill] (-2,-3.45) circle [radius=0.1];
  
    \draw [fill] (3,4.45) circle [radius=0.1];
  \draw [fill] (-3,4.45) circle [radius=0.1];
    \draw [fill] (1.8, 6) circle [radius=0.1];
  \draw [fill] (-1.8, 6) circle [radius=0.1];
  \draw (3, 4.45) -- (-3, 4.45);
  
  \draw (2, 3.45) -- (1.8, 6) -- (3, 4.45) -- (2, 3.45);
   \draw (-2, 3.45) -- (-1.8, 6) -- (-3, 4.45) -- (-2, 3.45);
   
   \draw (4,0) -- (1.8, 6);   
    \draw (4,0) -- (3, 4.45);  
     \draw (-4,0) -- (-1.8, 6);   
    \draw (-4,0) -- (-3, 4.45);  
    
    \draw[fill = gray, opacity = .5] (-2,3.45) -- (2,3.45) -- (1.8, 6) -- (-1.8, 6) -- (-2, 3.45);
         \draw[dotted, thick] (2, 3.45) -- (-1.8, 6);
          \draw (-2, 3.45) -- (1.8, 6);
 \draw (-1.8,6) -- (1.8, 6);
 
 \node[above] at (0, 6) {$\sigma$};
\end{tikzpicture}
\caption{$D_2(\Delta)$ when $L$ is a $1$-cycle and each sphere is $1$-dimensional.}
\label{f:cycle}
\end{figure}
\begin{Remark}
It may be confusing why we chose to replace each vertex in $O_mL$ with the boundary of an $m$-simplex. 
In fact, at first it seemed more natural to us to replace each vertex with an $(m-1)$-octahedron, as this would give $O_mL$ a simple flag triangulation. 
However, we could not find a way to extend the definition of \cite{ados} to this case. 
We'll illustrate this with a simple example, see Figure~\ref{f:cycle}. 

Suppose $L$ is a $1$-cycle, and replace two of the vertices with cellulated $S^1$'s. 
We need to construct a $4$-cycle $\beta \in C_4(\cac(D_m(\Delta)); \zz_2)$. 
For this to be the case, then for any $3$-cell $\sigma$ in $\cac(D_2(\Delta))$, the collection of $1$-cells $\{ \tau | \{\sigma, \tau\} \in \beta\}$ must form a $1$-cycle. 
If we triangulate the $S^1$'s as the boundary of a $2$-simplex, then for each such $\sigma$ there is a natural $1$-cycle containing all the vertices not contained in $\sigma$. On the other hand, if we replace each vertex with a $1$-octahedron, then we have to make a choice of $1$-cycle to pair with $\sigma$.  We could not find a way to do this consistently.
\end{Remark}

\remark{(\emph{Homology below the top dimension}).}\label{s:ast}
In Section \ref{hypcomp}, we will need a slight generalization of the previous arguments, where we consider complexes with homology below the top dimension. 

Suppose $L$ is a $d$-dimensional flag complex, and let $C$ be the support of a cycle in $H_{k}(L, \zz_2)$ for $k < d$. If $C$ is a full subcomplex, then the arguments in Theorem \ref{l:cycle} and Theorem \ref{t:neq0} generalize to show that $\vktwo^{k + \delta}(O_mL) \ne 0$.

If $C$ is not full, then $ \gO_m$ may not be a cycle if we choose $\Delta$ incorrectly. However, the argument generalizes if the following $\ast$-condition is satisfied, see \cite{ados}:

\[
	\tag{$*$} \text{For all $\gs,\gt \in C$ with $\Delta^0 \subset \sigma \cup \tau$ we have $\sigma \cap \tau \subset \Delta$.}
\]

We do not know an example of a $d$-dimensional complex $L$ and a class $\phi \in H_k(L;\mathbb{Z}_2)$ such that the $\ast$-condition fails for the support of every representative $C$ for $\phi$. 

One instance where the $\ast$-condition is automatically satisfied is if the cycle is in the top dimension. 

\begin{Theorem}\label{act=2n}
If $L$ is a $d$-dimensional flag complex and $H_d(L;\zz_2) \ne 0$, then $\vk^{d+\delta} (O_mL) \ne 0$.
\end{Theorem}

\section{Obstructors for hyperplane complements}\label{hypcomp}

In Subsection \ref{ss:config} we defined various configurations of standard subgroups for  simple complexes of groups. In this section, we will show that for any finite arrangement $\ca$ of affine hyperplanes in $\cc^n$ there is a configuration of abelian subgroups in the fundamental group of the complement $\pi_1(M(\ca))$, indexed by the simplices in  a certain simplicial complex, which is 
homeomorphic to the geometric realization $|\cq(\ca)|$ of the intersection poset.
If $\ca$ satisfies certain conditions, this simplicial complex will satisfy the $\ast$-condition in Remark \ref{s:ast} below. When these conditions hold, the obstructor dimension method will imply that if $\ca$ is irreducible, essential, and not central, then $\actdim(\pi_1(M(\ca))) \ge 2n$.  In particular, when $M(\ca)$ is aspherical, $\actdim(\pi_1(M(\ca))) = 2n$.

\subsection{Free abelian subgroups}
Many of the terms which we use in this subsection were defined earlier in Subsection \ref{ss:hyp}.  Recall that the intersection poset $\cq(\ca)$ is ordered by reverse inclusion.
Given $G \in \cq(\ca)$, let $\ca_G := \{H \in \ca \mid H\le G\}$ be the induced central subarrangement of hyperplanes containing $G$. We will use these central subarrangements to construct free abelian subgroups of $\pi_1(M(\ca))$. The same construction of these free abelian groups is given in \cite{dsy}.

\begin{lemma}\label{subgp} 
For any $H \le G \in \cq(\ca)$, $\pi_1(M(\ca_H))$ injects into $\pi_1(M(\ca_G))$.
\end{lemma}

\begin{proof}
There is a natural inclusion $j: M(\ca_G) \rightarrow M(\ca_H)$. We define a map $f: M(\ca_H) \rightarrow M(\ca_G)$ by first choosing a point $x$ in $H$ and a small ball $B_x$ that only intersects hyperplanes in $M(\ca_H)$. We can deformation retract $M(\ca_H)$ to $B_x$ and then compose with the inclusion $B_x \rightarrow M(\ca_G)$. Clearly, $j \circ f$ is homotopic to the identity, and therefore $f_\ast: \pi_1(M(\ca_H)) \rightarrow \pi_1(M(\ca_G))$ is injective.
\end{proof}

\begin{lemma}\label{l:center}
For any central arrangement $\ca$, $\pi_1(M(\ca))$ has an infinite center. 
\end{lemma}

\begin{proof}
There is a projectivization map $p: \mathbb{C}^n - \{0\} \rightarrow \mathbb{C}\mathbb{P}^{n-1}$ with fiber $\cc^\ast$. The restriction to $M(\ca)$ is a trivial bundle (\cite{ot}*{Proposition 5.1}). Then $\gamma = i(\pi_1(\cc^\ast)) \subset M(\ca)$ is in the center of $\pi_1(M(\ca))$. 
\end{proof} 

Of course, if the central arrangement $\ca$ is  reducible, then the center of $\pi_1(M(\ca))$ has rank greater than one (take central elements from each factor). If $\ca$ is irreducible, central, and $M(\ca)$ is aspherical, it turns out that the center is infinite cyclic. 

Next, suppose $G \in \cq(\ca)$ is such that $\ca_G$ is irreducible. By the previous lemma, there is an  element $\gamma_G$ of infinite order in the center of $\pi_1(M({\ca}_G))$. Furthermore, if $G_1 < G_2 < \dots < G_n$ is a chain in $\cq(\ca)$ with each $\ca_{G_i}$ irreducible, then since 
	\[
	\pi_1(M(\ca_{G_1})) \subset \pi_1(M(\ca_{G_2})) \subset \dots \subset \pi_1(M(\ca_{G_n})),
	\]
we obtain a free abelian group of rank $n$ generated by $\gamma_{G_1}, \gamma_{G_2}, \dots, \gamma_{G_n}$ (its rank is $n$ by  Theorem \ref{homology} below).
On the other hand, suppose $\ca_G$ decomposes as a product of irreducibles: 
	\[
	\ca(G)\cong\ca_{G_1}\times \cdots\times \ca_{G_k}.
	\]
Then
	\[
	M(\ca_G) \cong M(\ca_{G_1}) \times M(\ca_{G_2}) \times \dots \times M(\ca_{G_k}) 
	\] 
and we obtain further free abelian groups as products of the free abelian subgroups in the fundamental groups of the factors. 

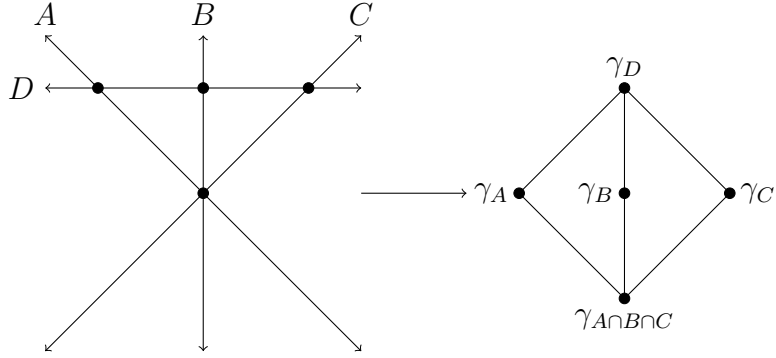
\begin{figure}\label{fig:cxarr}
   \centering
   \begin{tikzpicture}[scale = .7]
   
   \draw[<->] (-3,3) -- (-3,-3);
   \draw[<->](-6,3) -- (0,-3);
   \draw[<->] (-6,-3) -- (0,3);
   \draw[<->] (-6,2) -- (0,2);
   \draw[->] (0,0) -- (2,0);
   
   \node[above] at (-6,3) {$A$};
   \node[above] at (-3,3) {$B$};
   \node[left] at (-6,2) {$D$};
   \node[above] at (0,3) {$C$};
   
   \draw (3,0) -- (5,-2);
   \draw (5,0) -- (5,-2);
   \draw (7,0) -- (5,-2);
   \draw (3,0) -- (5,2);
   \draw (5,0) -- (5,2);
   \draw (7,0) -- (5,2);
   
   \node[above] at (5,2) {$\gamma_{D}$};
     \node[below] at (5,-2) {$\gamma_{A \cap B \cap C}$};
      \node[left] at (3,0) {$\gamma_{A}$};
      
      \draw [fill] (-3,2) circle [radius=0.1];
      \draw [fill] (-3,0) circle [radius=0.1];
      \draw [fill] (-5,2) circle [radius=0.1];
      \draw [fill] (-1,2) circle [radius=0.1];
      \draw [fill] (5,2) circle [radius=0.1];
      \draw [fill] (5,-2) circle [radius=0.1];
      \draw [fill] (3,0) circle [radius=0.1];
      \draw [fill] (7,0) circle [radius=0.1];
      \draw [fill] (5,0) circle [radius=0.1];

        \node[left] at (5,0) {$\gamma_{B}$};
          \node[right] at (7,0) {$\gamma_{C}$};
   \end{tikzpicture}
   
 \caption{A real line arrangement and an associated configuration of abelian groups for its complexification. }
 \end{figure}

We will produce a configuration of abelian groups based on a simplicial complex with vertex set  $\{G \in \cq(\ca)\mid \ca_G\text{ is irreducible}\}$. Two elements $G$ and $H$ of $\cq(\ca)$ are \emph{comparable} if $H<G$ or $G<H$.  Distinct vertices $v_G$ and $v_H$ are connected by an edge if and only if 1) $G$ and $H$ are comparable or 2) $\ca_{G\, \cap\, H} = \ca_G \times \ca_H$ (see Figure \ref{fig:cxarr}).

To properly describe the simplicial complex and the corresponding configuration of abelian groups, we need the notion, introduced by De Concini and Procesi  \cite{dp}, of a \emph{building set} for $\cq(\ca)$.

\subsection{Building sets}
Given a collection of subspaces $\cg$ in $\cq(\ca)$ and an element $X \in \cq(\ca)$, let $\cg_{\le X}$ denote the set of  elements in $\cg$ that contain $X$.  Let $\max \cg_{\le X}$ be the set of maximal elements of $\cg_{\le X}$.

\begin{definition}
A collection $\cg$ of subspaces in $\cq(\ca)$ is a \emph{building set} if for any $X \in \cq(\ca)$ such that $\max \cg_{\le X} = \{G_1,G_2, \dots G_k\}$, we have $$\ca_X \cong \ca_{G_1} \times \ca_{G_2} \times \dots \times \ca_{G_k}.$$
\end{definition}
There are two canonical choices for a building set. 
Note that $\cq(\ca)$ is itself a building set, since $\max \cq(\ca)_{\le X} = X$. 
Also, note that the poset $$\viq = \{G \in \cq(\ca) \mid \ca_G \text{ is irreducible}\}$$
 is a building set. The set $\viq$ is called the \emph{set of irreducibles} in $\cq(\ca)$; eventually, $\viq$ will be the vertex set of a simplicial complex $\iq$ called the \emph{irreducible complex}.
In the case of $\viq$, if we decompose $\ca_X$ into irreducible components 
	\[
	\ca_X \cong \ca_{G_1} \times \ca_{G_2} \times \dots \times \ca_{G_k},
	\] 
then $\max \viq_{\le X} = \{G_1,G_2,\dots, G_k\}$.
In fact, by considering $X \in \cq(\ca)$ such that $\ca_X$ is irreducible, it is obvious that every building set must contain $\viq$. 

Any building set determines a collection of ``nested subsets''. These subsets will be the simplices of a simplicial complex on which our configuration of abelian groups will be based. 

\begin{definition}
Let $\cg$ be a building set for an affine arrangement $\ca$.  A subset $\alpha \subset \cg$ is \emph{$\cg$-nested} if for any subset $\{X_i\}$ of $\ga$ consisting of pairwise incomparable elements $X_i \in \alpha$, the intersection $\cap X_i$ is nonempty and does \emph{not} belong to $\cg$.

\end{definition}

Note that if  $X_1 > X_2 > \dots > X_n$ is any chain in $\cg$, then $\{X_1,\dots,X_n\}$ is a $\cg$-nested subset. It is also obvious from the definition that the nested subsets form a simplicial complex, that is, if $\beta \subset \alpha$ and $\alpha$ is $\cg$-nested, then $\beta$ is $\cg$-nested. For example, if the building set is $\cq(\ca)$, then the nested set complex is the barycentric subdivision $|\cq(\ca)|$, as the nested subsets are precisely chains in $\cq(\ca)$. Let $\iq$ denote the simplicial realization of the $\viq$-nested subsets.

\begin{lemma}\label{abelnested}
Let $\alpha$ be a simplex in $\iq$. For each vertex of $\alpha$ corresponding to $G \in \viq$, let $\gamma_G$ denote the central element in $\pi_1(M(\ca_G))$ corresponding to the fiber of the Hopf fibration as defined in Lemma \ref{subgp}. Then the subgroup generated by $\{\gamma_G\}_{G \in \vertex \alpha}$ is free abelian. 
\end{lemma}

\begin{proof}
We claim if $G$ and $H$ are connected by an edge in $\viq$, then $\gamma_G$ and $\gamma_H$ commute. Two elements $G$ and $H$ of $\viq$ are \emph{comparable} if $G<H$ or $H<G$.  A two element subset  $\{G,H\}$ with two elements is nested if and only if $H$ and $G$ are comparable or if $G \cap H \ne \emptyset$ and $G \cap H \notin \viq$. In the first case, $\gamma_G$ and $\gamma_H$ commute by Lemma \ref{subgp}. For the second case, we have that $\ca_{G \cap H}$ is reducible, hence $\ca_{G \cap H} = \ca_{G_1} \times \ca_{G_2} \times \dots \times\ca_{G_k}$. If $\ca_G$ and $\ca_H$ are contained in a single $\ca_{G_i}$, this would imply that $G_i \in G \cap H$, which is a contradiction. This implies $\gamma_G$ and $\gamma_H$ commute. Since each $\gamma_G \cong \zz$ we are done.
\end{proof}

\begin{Theorem}[\cite{dsy}*{Corollary 4.5}]\label{homology}
For each $G \in \viq$, the image $\overline{\gamma_G}$ of $\gamma_G$ in $H_1(M(\ca),\zz)$ satisfies the relation 
	\[
	\overline{\gamma_G} = \sum_{\substack{H\in \ca\\ H\le G}} \overline{\gamma_H}.
	\] 

Furthermore, for any simplex $\ga$ of $\iq$, $\{\overline{\gamma_G} \mid G\in \vertex \ga\}$ is linearly independent. 
\end{Theorem}

\begin{corollary}
For any simplex $\alpha \in I\cq(\ca)$, the free abelian subgroup constructed in Lemma \ref{abelnested} is free abelian of rank $\dim \ga + 1$.
\end{corollary}

In general, the simplicial complex $\iq$ formed from the $\viq$-nested subsets is not a flag complex. For example, if $\ca$ is the complexification of the real arrangement which consists of $n$ general position lines in $\rr^2$, then $\iq$ for the complexification is the one-skeleton of an $n$-simplex. Our configuration of abelian groups will always be based on the flag completion of $\iq$, where the flag completion is the unique flag complex with the same one-skeleton as $I\cq(\ca)$. 

\begin{Theorem}\label{t:abelcomplex}
Let $\ca$ be an affine arrangement, and let $\fiq$ be the flag completion of $\iq$. Then $\pi_1M(\ca)$ admits a configuration of free abelian groups based on $\fiq$.
\end{Theorem}

\begin{proof}
Assume that $\sigma$ and $\tau$ are two simplices in $\fiq$, and $\zz^\gs$ and $\zz^\gt$ the corresponding free abelian subgroups we have constructed. 
We have seen above that the subgroups $\zz^\gs$ and $\zz^\gt$ have ranks equal to  $\dim\gs + 1$ and $\dim\gt +1$, respectively. 
We must show that these subgroups intersect in the subgroup corresponding to $\zz^{\gs \cap \tau}$. 
Since each of these subgroups maps isomorphically onto its image in $H_1(M(\ca),\zz)$, it suffices to prove that the images of the subgroups intersect correctly. 
However, this follows immediately from Theorem \ref{homology}: it is obvious that $\zz^{\gs \cap \gt} \subset \zz^{\gs} \cap \zz^\gt$ and since each of the $\gamma_G$ are linearly independent we have $\zz^\gs \cap \zz^\gt \subset \zz^{\gs \cap \gt}$.
\end{proof}

\begin{remark}
Note that the simplicial complex we construct does not capture all of the commuting relations between the standard abelian subgroups. For example, for the arrangement depicted in Figure, \ref{fig:cxarr}, the elements $\gamma_D$ and $\gamma_{A \cap B \cap C}$ commute, since $\gamma_D$ commutes with $\gamma_A$, $\gamma_B$ and $\gamma_C$. In this configuration of abelian groups it is only important that $\zz^\gs \cap \zz^\gt = \zz^{\gs \cap \gt}$.
\end{remark}

\subsection{The homotopy type of $\iq$} 

Suppose $\ca$ is a finite  arrangement of affine hyperplanes in $\cc^n$ of rank $l$.  By a theorem of Folkman \cite{folkman}, $\vert\cq(\ca)\vert$, the geometric realization of the intersection poset, is homotopy equivalent to a wedge of spheres, each of dimension $(l-1)$.  The number of spheres in this wedge is an integer, $\gb_\ca$.  Alternatively, $\gb_\ca = \vert \chi(M(\ca))\vert$, where $\chi(M(\ca))$ is the Euler characterisitic of the complement (see \cite{ot}).  So, if $$p_\ca(t) = \sum_i b_i(M(\ca))t^i$$ is the Poincar\'e polynomial of $M(\ca)$, then $\chi(M(\ca))=p_\ca(-1)$.

If $\ca$ is central, then $\vert\cq(\ca)\vert$ is a cone and hence, $\gb_\ca=0$ and $p_\ca(t)$ has $(1+t)$ as a factor.  Conversely, by a theorem of Crapo \cite{crapo}, if $\gb_\ca=0$, then $\ca$ decomposes as $\ca\cong \ca'\times \ca_1$ where $\ca_1$ is a nontrivial central arrangement.  It follows that $p_\ca(t)$ has $(1+t)^k$ as a factor if and only if $\ca$ can be decomposed as $\ca'\times \ca_1\times \dots \times \ca_k$, where $\ca_i$ are nontrivial central arrangements.

Feichtner and M{\"u}ller showed that if $\cg$ and $\cg'$ are two building sets for a hyperplane arrangement, the simplicial complexes corresponding to the nested subsets are homeomorphic via a series of stellar subdivisions \cite{fm}. 
In particular, this implies that $\iq$ is a simplicial complex with $H_{\rk(\ca)-1}(\iq, \mathbb{Z}_2) \ne 0$ if $\ca$ is irreducible and has no central factor.

If $\ca$ is a central arrangement, there is an affine arrangement $d(\ca) \subset \cc^{n-1}$, called the ``deconing'' of $\ca$, which is obtained by projectivizing $\ca$ and then removing any projectivized hyperplane. The realization of the intersection poset $|\cq(\ca)|$ is isomorphic to the cone on $|\cq(d(\ca))|$, so therefore the Poincar\'e polynomials satisfy $p_\ca(t) = (1+t)p_{d(\ca)}(t)$. If $\ca$ is essential and irreducible, this implies by Crapo's theorem that $d(\ca)$ is essential and has no central factor.

\subsection{Obstructor dimensions of arrangement complements}

Suppose that $\ca$ is an irreducible, essential, and noncentral arrangement in $\cc^n$. We constructed a simplicial complex $\iq$ with $H_{n-1}(\iq, \zz) \ne 0$. The configuration of abelian groups is obtained by taking the flag completion $\fiq$. Since the dimension of $|\cq(\ca)|$ is always $\rk(\ca) -1$, it is a flag simplicial complex with top dimensional homology if $\ca$ is essential and has no central factors.

Therefore, we can apply Theorem~\ref{t:vk} if $\fiq$ is $(n-1)$-dimensional. We first show this holds for the case that $M(\ca)$ is aspherical. 

\begin{lemma}
Let $\ca$ be a complex hyperplane arrangement such that $M(\ca)$ is aspherical. Then $\dim (\fiq) \le n-1$.
\end{lemma}

\begin{proof}
It is a standard fact that $M(\ca)$ is homotopy equivalent to an $n$-dimensional complex.
If $M(\ca)$ is aspherical, then  $\gdim(\pi_1(M(\ca))) \le n$. Therefore, the rank of any abelian subgroup of $M$ is $\le n$, which implies $\dim (\fiq) \le n-1$.
\end{proof}

This along with Theorem \ref{t:vk} and Theorem \ref{t:abelcomplex} immediately implies Theorem C from the introduction.

\begin{Theorem}\label{t:nothe}
If $M(\ca)$ is an essential, aspherical arrangement with no central factors, then $\actdim(\pi_1(M(\ca))) = 2n$. In particular, $M(\ca)$ is not homotopy equivalent to a $(2n-1)$-manifold.
\end{Theorem}

A similar argument shows a general result for central arrangements. 

\begin{Theorem}
Let $\ca$ be an aspherical, irreducible, essential, central arrangement in $\cc^n$. Then $\actdim(\pi_1(M(\ca))) = 2n-1$. 
\end{Theorem}
\begin{proof}
 If $\ca$ is essential and irreducible, then Theorem \ref{mainhyper} implies that if $d(\ca)$ is the deconing of $\ca$, then $\obdim(\pi_1(M(d(\ca)))) = 2n-2$.  Then since $M(\ca) \cong M(d(\ca)) \times S^1$, we have $\pi_1(M(\ca)) = \mathbb{Z} \times \pi_1(M(d(\ca)))$ and hence $\obdim ( \pi_1(M(d(\ca)))) = \actdim( \pi_1(M(d(\ca)))) = 2n-1$.
\end{proof}

The product formula for obstructor dimension gives the following corollary which is the general answer for the obstructor dimension of aspherical hyperplane arrangements.
\begin{corollary}
Let $\ca$ be an affine aspherical arrangement in $\cc^n$ and suppose that $$\ca \cong \ca_1 \times \ca_2 \times \dots \times \ca_k \times \ca'$$ where each $\ca_i$ is irreducible and central, and $\ca'$ has no central factor. Then $$\obdim(\pi_1(M(\ca))) = \actdim(\pi_1(M(\ca))) = 2n - k.$$
\end{corollary}

Using Theorem \ref{t:nothe}, this gives the following computation of the action dimension of any spherical Artin group (cf.\ Proposition~\ref{p:sphericalactdim}). 
\begin{corollary}[Le \cite{giang}]\label{l:sphericalactdim}
Suppose $A_\gs$ is an irreducible spherical Artin group of rank $d+1$. Then $\obdim A_\gs=\actdim A_\gs=2d+1$. Therefore, the action dimension of a spherical Artin group is the sum of the action dimensions of its irreducible factors. 
\end{corollary}

When $M(\ca)$ is not aspherical, $\fiq$ could have much larger dimension than $\iq$. For example, we can take any arrangement and add new hyperplanes in general position. Then $\iq$ for the new arrangement will have arbitrarily high dimensional flag completion, as the new hyperplanes themselves induce the $1$-skeleton of a $n$-simplex in $\iq$. We now give a condition on our arrangement that guarantees that $\fiq$ contains a simplex which satisfies the $\ast$-condition from Subsection \ref{s:ast}.

\begin{definition}\label{d:ic}

Let $\ca$ be a complex hyperplane arrangement. A \emph{complete chain of irreducibles} is a chain of subspaces $G_0 > G_1 > \dots > G_n$ such that each $\ca_{G_i}$ is irreducible.
\end{definition}

We claim that the simplex $\sigma := [G_0,G_1, \dots G_n]$ in $\fiq$ satisfies the $\ast$-condition. Note that if a simplex does not satisfy the $\ast$-condition, then there is a vertex in $\fiq$ that is connected to each of the vertices of that simplex by an edge. So, suppose to the contrary that $H$ is a subspace which is connected to each $G_i$ by an edge in $\fiq$.
Now, since $H$ is connected to $G_n$ in $\fiq$ and $G_n$ is $0$-dimensional, it must contain $G_n$. 
Let $G_i$ be the maximal subspace in the chain that is not contained in $H$. Since $H$ and $G_i$ are connected by an edge in $\iq$, we must have $\ca_{H \cap G_i} \cong \ca_H \times \ca_{G_i}$. Since $H$ contains $G_{i+1}$ we must have $H \cap G_i = G_{i+1}$. Therefore, the splitting $\ca_{H \cap G_i} \cong \ca_H \times \ca_{G_i}$ would contradict $G_i \in \viq$. 

Therefore, we have constructed a $d$-dimensional flag complex $\fiq$ such that $H_{n-1}(\fiq, \zz_2) \ne 0$. 
Since $(\fiq, \sigma)$ satisfies the $\ast$-condition, we can apply \cite{ados} to the complex $O(\fiq)$ to get the following:

\begin{Theorem}\label{mainhyper} 
Let $\ca$ be an arrangement in $\cc^n$ that is essential and not central. If $\ca$ contains a complete chain of irreducibles, then $$\actdim(\pi_1(M(\ca))) \ge \obdim(\pi_1(M(\ca))) \ge 2n.$$
\end{Theorem}

If $\ca$ is an inessential arrangement, then we can still compute lower bounds for the action dimension. This is because $M(\ca)$ splits as $\mathbb{C}^k \times \ca'$, where $\ca'$ is an essential arrangement in $\mathbb{C}^{n-k}$, and if $\ca$ is not central then $\ca'$ is not central. Conversely, our results say nothing about general position hyperplane arrangements, though this is not very interesting in this context. Hattori showed in \cite{hat} that the complement of a general position arrangement has free abelian fundamental group and that the arrangement is homotopy equivalent to a certain skeleton of a $k$-torus.

\begin{example}\label{fiber}
We now describe another type of arrangement whose complement is always aspherical. 
First, an arrangement $\ca$ is said to be \emph{strictly linearly fibered} over $\ca_G$ if $G$ is a line and the restriction of the projection $\pi_G : \cc^n \rightarrow \cc^n/G$ to $M(\ca)$ is a fiber bundle projection. 
A fiber-type arrangement is defined inductively: $\ca \subset \cc^n$ is fiber-type if there is a line $G$ such that $\ca$ is strictly linearly fibered over $G$ and the induced arrangement $\pi(\ca) \subset \cc^{n-1}$ is fiber-type. 
For example, the braid arrangement is fiber-type. 
The action dimension of these examples was already known by work of Falk and Randell in \cite{fr} and results in \cite{bkk}. 
Indeed, Falk and Randell showed that the fundamental group of the complement of a fiber-type arrangement is an iterated semidirect product of free groups, so the computation of action dimension followed from Corollary 27 in \cite{bkk}.
\end{example}

\begin{example}
A \emph{complex reflection} is a periodic affine automorphism of $\cc^n$ whose fixed point set  a complex hyperplane.
A \emph{complex reflection group} is a finite group acting on $\cc^n$ by complex reflections. For example, every finite Coxeter group is a complex reflection group by complexification of the action on $\mathbb{R}^n$. There groups were completely classified by Shepard and Todd \cite{shepardtodd1954}, who showed that they fit into several infinite families depending on $3$ parameters and 34 exceptional cases. The complement of the fixed hyperplanes is a central arrangement, and the fundamental groups of such hyperplane complements can be thought of as generalizations of spherical Artin groups. It is known that all such hyperplane complements are aspherical (the remaining exceptional cases were resolved in \cite{bessis}). Therefore, if the arrangement for a finite reflection group is essential and irreducible, then the action dimension of $\pi_1(M(\ca))$ is precisely $2n-1$. 
\end{example}

\section{Questions}

Here are four questions that came up during our work. 
When the   $d$-dimensional flag complex $L$ is EDCE and $G$ is the fundamental group for graph product complex of closed aspherical $m$-manifolds (or more generally, the group associated to  a complex of closed aspherical manifolds), we showed that $\actdim G \le (m+1)(d+1)-1$.  On the other hand, in order to show that the corresponding van Kampen obstruction is $0$, we only need the weaker assumption $H_d(L,\zz_2)=0$.

\begin{Question} \label{q0}
Is our upper bound for $\actdim G$ still valid when the hypothesis that $L$ is EDCE is replaced by the hypothesis $H_d(L,\zz_2)=0$?
\end{Question}

\begin{Question}\label{q1}
If $OL$ piecewise linearly embeds in a sphere of codimension $k$, does $O_mL$ piecewise linearly embed in a sphere of codimension $k$? Together with the main theorem of \cite{ados} this would imply  that if $L$ is a flag complex with $H_n(L,\zz_2) = 0$, then $\embdim(O_mL) < d + \delta$.  (Here $\embdim(O_mL)$ means the minimum dimension of a sphere into which there is a PL embedding of $O_mL$.)
\end{Question} 

\begin{Question}\label{q2}
 Let $K_L$ be a polyhedral join of simplicial complexes $K_s$ over $L$. Is there a formula for the van Kampen obstruction of $K_L$ in terms of the van Kampen obstructions of the $K_s$ and the homology of $L$?
 
 \end{Question}

\begin{Question}\label{q3}

Suppose that $\ca$ is an essential, noncentral arrangement which admits a complete chain of irreducibles (see Definition~\ref{d:ic}). Is it possible for $M(\ca)$ to be homotopy equivalent to a $(2n-1)$-manifold?
\end{Question}


\bibliography{1}

Michael W. Davis, Department of Mathematics, The Ohio State University, 231 W. 18th Ave., Columbus, Ohio 43210, \url{davis.12@osu.edu}

Giang Le, Department of Mathematics, Oregon State University, 368 Kidder Hall, Corvallis, OR 97331, \url{giangl@oregonstate.edu}

Kevin Schreve, University of Michigan-Ann Arbor, Department of Mathematical Sciences, PO Box 413, Ann Arbor, MI 48109, \url{schreve@umich.edu} 

\obeylines
\end{document}